\documentclass[a4paper,10pt,twoside,openright]{article}
\usepackage[T1]{fontenc}
\usepackage[utf8]{inputenc}
\usepackage{courier}
\usepackage{amsmath}
\usepackage{amsthm}
\usepackage{amssymb}
\usepackage{mathrsfs}
\usepackage{bbm}
\usepackage{graphicx}
\usepackage{epstopdf}
\usepackage{subfigure}
\usepackage{enumerate}
\usepackage{enumitem}
\usepackage{fancyhdr}
\usepackage{tikz}
\usepackage{setspace}
\usepackage{makeidx}
\usepackage{thmtools}
\usepackage{thm-restate}
\usepackage{hyperref}
\usepackage{cleveref}
\usepackage[margin=35mm]{geometry}
\setlist{nolistsep}
\makeindex
\usetikzlibrary{arrows}
\hypersetup{colorlinks=true,linkcolor=blue}
\hypersetup{
	citecolor = purple,
}

\usepackage{stackengine}

\newtheorem{theorem}{Theorem}
\newtheorem{proposition}[theorem]{Proposition}
\newtheorem{lemma}[theorem]{Lemma}

\begin{document}

\begin{titlepage}

\vskip0.2truecm

\begin{center}

{\LARGE \bf A condition that implies full homotopical complexity of orbits}

\end{center}

\vskip  0.4truecm

\centerline {{\large Salvador Addas-Zanata}}
\vskip 0.1truecm
\centerline {{\large Bruno de Paula Jacoia}}

\vskip 0.2truecm

\centerline { {\sl Instituto de Matem\'atica e Estat\'\i stica }}
\centerline {{\sl Universidade de S\~ao Paulo}}
\centerline {{\sl Rua do Mat\~ao 1010, Cidade Universit\'aria,}} 
\centerline {{\sl 05508-090 S\~ao Paulo, SP, Brazil}}
 
\vskip 0.7truecm

\begin{abstract}

 We consider closed orientable surfaces $S$ of genus $g>1$ 
and  homeomorphisms $f:S\rightarrow S$ homotopic to the identity. A set 
of hypotheses is presented, called fully essential system of curves $\mathscr{C}$ 
and it is shown that under these hypotheses, the natural lift of $f$ to the universal 
cover of $S$ (the Poincar\'e disk $\mathbb{D}),$ denoted $\widetilde{f},$ has 
complicated and rich dynamics. In this context we generalize results that hold 
for homeomorphisms of the torus homotopic to the identity 
when their rotation sets contain zero in the interior. 
In particular, we prove that if $f$ is a $C^{1+\epsilon }$ diffeomorphism for some 
$\epsilon >0$ and $\pi :\mathbb{D}\rightarrow S$ is the covering map, then there 
exists a contractible hyperbolic $f$-periodic saddle point $p\in S$ such that for 
any $\widetilde{p}\in \pi ^{-1}(p),$

 $$W^u(\widetilde{p}) \pitchfork W^s(g(\widetilde{p})) $$
for all deck transformations $g\in Deck(\pi ).$ By $\pitchfork,$ we mean a 
topologically transverse intersection between the manifolds, see the 
precise definition in subsection 1.1. 
We also show that the homological rotation set of such a $f$ is a compact convex 
subset of $\mathbb{R}^{2g}$ with maximal dimension and all points in its interior 
are realized by compact $f$-invariant sets, periodic orbits in the rational case, and  
 $f$ has uniformly bounded displacement with respect to 
rotation vectors in the boundary of the rotation set. Something that  implies, in case $f$ is 
area-preserving, that the rotation vector of  
Lebesgue measure belongs to the interior of the rotation set.

\end{abstract} 

\vskip 1.3truecm

\noindent{\bf e-mails:} sazanata@ime.usp.br and bpjacoia@ime.usp.br

\vskip 2.3truecm


\vskip 0.3truecm

\vfill
\hrule
\noindent{\footnotesize{The first author is partially supported 
by CNPq, grant: 306348/2015-2 and the second author was supported by CNPq, grant: 159853/2013-3.}}

\end{titlepage}

\baselineskip=6.2mm

\section{Introduction}

\subsection{Preliminaries}

The main motivation for this work is to generalize some results that hold
for homeomorphisms and diffeomorphisms of the torus homotopic to the
identity, to homeomorphisms and diffeomorphisms of closed surfaces of higher
genus (for us, higher genus means larger than $1$), also in the homotopic to
the identity class.

In the study of torus homeomorphisms, a useful concept inherited from
Poincar\'e's work on circle homeomorphisms is that of rotation number, or in
the two-dimensional case, rotation vectors. Actually in the two-dimensional
setting, one usually does not have a single rotation vector, but a rotation
set, which is most precisely defined as follows: Given a homeomorphism $f:%
\mathbb{T}^2\rightarrow \mathbb{T}^2$ homotopic to the identity and a lift
of $f$ to $\mathbb{R}^2$, $\widetilde f:\mathbb{R}^2\to \mathbb{R}^2$, the
Misiurewicz-Ziemian rotation set $\rho (\widetilde f)$ is defined as (see 
\cite{mmkz89}):

\begin{equation}
\label{rotsetident} \rho(\widetilde{f}) = \bigcap_{i\geq 1}\overline{%
\bigcup_{n\geq i}\bigg\{ \frac{\widetilde{f}^n(\widetilde{p})-\widetilde{p}}{%
n} : \widetilde{p} \in \mathbb{R}^2 \bigg\}} 
\end{equation}

This set is a compact convex subset of $\mathbb{R}^2$ (see \cite{mmkz89}),
and it was proved in \cite{jf89} and \cite{MZ} that all points in its
interior are realized by compact $f$-invariant subsets of $\mathbb{T}^2$,
which can be chosen as periodic orbits in the rational case. By saying that
some vector $v\in \rho (\widetilde f)$ is realized by a compact $f$%
-invariant set, we mean that there exists a compact $f$-invariant subset $%
K\subset \mathbb{T}^2$ such that for all $p\in K$ and any $\widetilde p\in
\pi ^{-1}(p),$ where $\pi :\mathbb{R}^2\to \mathbb{T}^2$ is the associated
covering map, the following holds

\begin{equation}
\label{deffrotvect} \lim_{n \to \infty}\frac{\widetilde{f}^n(\widetilde{p})- 
\widetilde{p} }{n} = v. 
\end{equation}

Moreover, the above limit, whenever it exists, is called the rotation vector
of the point $p$, denoted $\rho (p)$.

Before presenting the results in the torus that we want to generalize to
other surfaces, we need a definition:

\begin{description}
\item[Definition (Topologically transverse intersections):]  If $M$ is a
surface, $f:M\rightarrow M$ is a $C^1$ diffemorphism and $p,q\in M$ are $f$%
-periodic saddle points, then we say that $W^u(p)$ has a topologically
transverse intersection with $W^s(q)$ (and write $W^u(p)\pitchfork W^s(q)$),
whenever there exists a point $r\in W^s(q)\cap W^u(p)$ ($r$ clearly can be
chosen arbitrarily close to $q$ or to $p$) and an open ball $B$ centered at $%
r$, such that $B\setminus \alpha =B_1\cup B_2$, where $\alpha $ is the
connected component of $W^s(q)\cap B$ which contains $r$, with the following
property: there exists a closed connected arc $\beta \subset W^u(p),$ such
that $\beta \subset B$, $r\in \beta $, and $\beta \setminus r$ has two
connected components, one contained in $B_1\cup \alpha $ and the other
contained in $B_2\cup \alpha $, such that $\beta \cap B_1\neq \emptyset $
and $\beta \cap B_2\neq \emptyset $. Clearly a $C^1$-transverse intersection
is topologically transverse. Note that as $\beta \cap \alpha $ may contain a
connected arc containing $r$, the ball $B$ may not be chosen arbitrarily
small.
\end{description}

\vskip 0.2truecm

{\bf Remark:\ }The consequence of a topologically transverse intersection
which is more relevant to us is a $C^0$ $\lambda $-lemma: If $W^u(p)$ has a
topologically transverse intersection with $W^s(q),$ then $W^u(p)$ $C^0$%
-accumulates on $W^u(q).$

\vskip 0.2truecm

In \cite{saz15} it is proved that if $(0,0)\in int(\rho (\widetilde f))$ and 
$f$ is a $C^{1+\epsilon }$- diffeomorphism for some $\epsilon >0$, then $%
\widetilde f $ has a hyperbolic periodic saddle point $\widetilde p\in 
\mathbb{R}^2$ such that

\begin{equation}
\label{c1eptoro}W^u(\widetilde p)\pitchfork W^s(\widetilde p)+(a,b), 
\end{equation}
for all $(a,b)\in \mathbb{Z}^2$ ($W^u(\widetilde p)$ is the unstable
manifold of $\widetilde p$ and $W^s(\widetilde p)$ is its stable manifold).
Note that as $\widetilde p$ is a periodic point for $\widetilde f$, the same
holds for all integer translations of $\widetilde p$ and moreover, for any
integer vector $(a,b),$ $W^{u,s}(\widetilde p+(a,b))=W^{u,s}(\widetilde{p}%
)+(a,b)$.

In the area-preserving case, this result implies that:

\begin{itemize}
\item  $\overline{W^u(\widetilde{p})}=\overline{W^s(\widetilde{p})}$ is a $
\widetilde{f}$-invariant equivariant closed connected subset of $\mathbb{R}^2
$ and there exists $M=M(f)>0$ such that any connected component $\widetilde{D%
}$ of $\left( \overline{W^u(\widetilde{p})}\right) ^c$ is an open
topological disk, whose diameter is less than $M$ and $D\stackrel{def.}{=}%
\pi (\widetilde{D})$ is a $f$-periodic disk. Moreover, for any $f$-periodic
disk $D\subset \mathbb{T}^2$, $\pi ^{-1}(D)\subset \left( \overline{W^u(
\widetilde{p})}\right) ^c.$

\item  for any $\rho =(s/q,r/q)\in int(\rho (\widetilde{f}))\cap \mathbb{Q}^2
$, if we consider the map $\widetilde{f}^q(\bullet )-(s,r),$ then there
exists a point $\widetilde{p}_\rho $ which is a hyperbolic periodic saddle
point for $\widetilde{f}^q(\bullet )-(s,r)$, its stable and unstable
manifolds have similar intersections as in (\ref{c1eptoro}) and 
$$
\overline{W^u(\widetilde{p}_\rho )}=\overline{W^s(\widetilde{p}_\rho )}=
\overline{W^u(\widetilde{p})}=\overline{W^s(\widetilde{p})}. 
$$

So, the above set is the same for all rational vectors in the interior of
the rotation set. We denote it by $R.I.(\widetilde{f})$ (region of
instability of $\widetilde{f}$) and a similar definition can be considered
in the torus: $R.I.(f)\stackrel{def.}{=}\pi (\overline{W^u(\widetilde{p})})=
\overline{W^u(p)}$, where $p=\pi (\widetilde{p})$ is $f$-periodic. Every $f$%
-periodic open disk in $\mathbb{T}^2$ is contained in a connected component
of the complement of $R.I.(f)$ and every such connected component is a $f$%
-periodic open disk, whose diameter when lifted to the plane is smaller than 
$M$.

\item  every open ball centered at a point of $R.I.(f)$ has points with all
rational rotation vectors contained in the interior of $\rho (\widetilde{f})$%
.

\item  if $f$ is transitive, then $\widetilde{f}$ is topologically mixing in
the plane. This follows easily from the fact that if $f$ is transitive, then 
$R.I.(f)=\mathbb{T}^2$ and $R.I.(\widetilde{f})=\overline{W^u(\widetilde{p})}%
=\overline{W^s(\widetilde{p})}=\mathbb{R}^2$.
\end{itemize}

As we already said, the above results were obtained in \cite{saz15} under a $%
C^{1+\epsilon }$ condition. In \cite{strict} and \cite{nancy}, some
analogous results were proved for homeomorphisms, by completely different
methods, but the conclusions of some are weaker.

What about surfaces of higher genus?

In this setting, starting with the definition of rotation set, things are
more involved. If $S$ is a closed orientable surface of genus $g>1$, the
definition of rotation set needs to take into account the fact that $\pi
_1(S),$ the fundamental group of $S$, and $H_1(S,\mathbb{Z})$, the first
integer homology group of $S$ are different: the first is almost a free
group with $2g$ generators. There is only one relation satisfied by the
generators. While the second is $\mathbb{R}^{2g}$.

Maybe the most immediate consequence of this is the fact that in order to
define a rotation set for surfaces of higher genus, if one wants it to be
not too complicated, and have some properties similar to what happens in the
torus, a homological definition must be considered. In the following we
present the definition of homological rotation set and homological rotation
vector as it appeared in \cite{akft15}. The idea to use homology in order to
define rotation vectors goes back to the work of Schwartzman \cite
{schwartzman}.

\subsection{Rotation vectors and rotation sets}

Let $S$ be a closed orientable surface of genus $g>1$ and $I:[0,1]\times
S\to S$ an isotopy from the identity map to a homeomorphism $f:S\to S.$

For $\alpha $ a loop in $S$ (a closed curve), $[\alpha ]\in H_1(S,\mathbb{Z}%
)\subset H_1(S,\mathbb{R})$ is its homology class. Recall that $H_1(S,%
\mathbb{Z})\simeq \mathbb{Z}^{2g}$ and $H_1(S,\mathbb{R})\simeq \mathbb{R}%
^{2g}.$ We will also consider $H_1(S,\mathbb{R})$ endowed with the stable
norm as in \cite{mg}, which has the property that $||[\gamma ]||\leq
l(\gamma )$ for any rectifiable loop $\gamma $, where $l(\gamma )$ is the
length of the loop.

For any fixed base point $b\in S,$ $\mathcal{A}_b=\{\gamma _p:p\in S\}$ is a
family of rectifiable paths such that $\gamma _p$ joins $b$ to $p$ and the
length of $\gamma _p$ is bounded by a uniform constant $C_{\mathcal{A}_b}.$

For any point $p \in S$ we want to construct a path in $S$ from $p$ to $%
f^n(p)$ and then form a loop by adding $\gamma_p$ and $\gamma_{f^n(p)}$.
Consider the path $I_p$ joining $p$ to $f(p)$ given by $t\mapsto I(t, p)$ .
Also, for each $n \in \mathbb{N}$ define the path $I_p^n$ joining $p$ to $%
f^n(p)$ by 
$$
I_p^n = I_p*I_{f(p)}*\ldots*I_{f^{n-1}(p)},%
$$
where $\beta * \delta$ is the concatenation of the path $\beta$ with the
path $\delta$.

For each $p \in S$ let $\alpha^n_p$ be the closed loop based at $b$ formed
by the concatenation of $\gamma_p$, the path $I_p^n$ in $S$ from $p$ to $%
f^n(p)$ and $\gamma_{f^n(p)}$ traversed backwards, that is 
$$
\alpha^n_p = \gamma_p * I^n_p * \gamma_{f^n(p)}^{-1}.%
$$

We can now define the homological displacement function of $p$ as 
$$
\Psi_f(p) = [\alpha_p].%
$$

For the function $\Psi _f:S\to H_1(S,\mathbb{R})$ we abbreviate its Birkhoff
sums as 
$$
\Psi _f^n(p)=\sum_{k=0}^{n-1}\Psi _f(f^k(p)). 
$$

Note that since $\alpha _p^n$ is homotopic to $\alpha _p*\alpha
_{f(p)}*\ldots *\alpha _{f^{n-1}(p)},$ 
$$
[\alpha _p^n]=\sum_{k=0}^{n-1}[\alpha _{f^k(p)}]=\sum_{k=0}^{n-1}\Psi
_f(f^k(p))=\Psi _f^n(p). 
$$

Also, the path $I_p^n$ can be replaced by any path joining $p$ to $f^n(p)$
and homotopic with fixed endpoints to $I_p^n$. This implies that $\Psi _f$
depends only on $f,$ on the choice of $\mathcal{A}_b$ and on the homotopy
class of the isotopy $I.$ In particular, $\Psi _f$ is bounded. Indeed, as $S$
is compact, sup$\{d_{\mathbb{D}}(\widetilde{q},\widetilde{f}(\widetilde{q})):
\widetilde{q}\in \mathbb{D}\}=C_{\max \_f}<\infty $, and if we replace the
path $I_p$ by the projection of the geodesic segment in $\mathbb{D}$ joining 
$\widetilde{p}\in \pi ^{-1}(p)$ to $\widetilde{f}(\widetilde{p}),$ as the
length of this path is smaller than $C_{\max \_f},$ then $||\Psi _f||\leq
2C_{\mathcal{A}_b}+C_{\max \_f}.$

As we just said, $\Psi _f$ depends on the choice of the basepoint $b$ and
the family $\mathcal{A}_b$. However, given another basepoint $b^{\prime }\in
S$ and a family $\mathcal{A}_{b^{\prime }}^{\prime }=\{\gamma _p:p\in S\}$
of rectifiable paths whose lengths are uniformly bounded by $C_{\mathcal{A}%
_{b^{\prime }}^{\prime }}$ such that $\gamma _p^{\prime }$ joins $b^{\prime }
$ to $p$, defining $\alpha ^{\prime }$$_p^n$ analogously, one has

\begin{equation}
\label{computation}[\alpha ^{\prime }{}_p^n]=[\gamma _p^{\prime
}*I_p^n*\gamma ^{\prime }{}_{f^n(p)}^{-1}]=[\alpha _p^n*\delta _p^n]=[\delta
_p^n]+\Psi _f^n(p), 
\end{equation}

where $\delta _p^n=\gamma _{f^n(p)}*\gamma ^{\prime }{}_{f^n(p)}^{-1}*\gamma
_p^{\prime }*\gamma _p^{-1}$. Indeed, the loop $\alpha ^{\prime }{}_p^n$ is
freely homotopic to $I_p^n*\delta _p^n$. In particular, if $\Psi _f^{\prime
}(p)=[\alpha _p^{\prime }]$, then

\begin{equation}
\label{bound}||\Psi _f^n(p)-\Psi _f^{^{\prime }n}(p)||\leq 2C_{\mathcal{A}%
_b}+2C_{\mathcal{A}_{b^{\prime }}^{\prime }}. 
\end{equation}

Finally, if the limit 
\begin{equation}
\label{defrotvech}\rho (f,p)=\lim _{n\to \infty }\frac 1n\Psi _f^n(p)\in
H_1(S,\mathbb{R}) 
\end{equation}
exists, we say that $p$ has a well-defined (homological) rotation vector.

After all this, we are ready to present the definition of the (homological)
rotation set of $f,$ which is analogous to the definition for the torus \cite
{mmkz89}. The Misiurewicz-Ziemian rotation set of $f$ over $S$ is defined as
the set $\rho _{mz}(f)$ consisting of all limits of the form 
$$
v=\lim _{k\to \infty }\frac 1{n_k}\Psi _f^{n_k}(p_k)\in H_1(S,\mathbb{R}), 
$$
where $p_k\in S$ and $n_k\to \infty $. By (\ref{bound}), the rotation set
depends only on $f,$ but not on the choice of the isotopy, the basepoint $b$
or the arcs $\gamma _p$. This definition coincides with 
$$
\rho _{mz}(f)=\bigcap_{m\geq 0}\overline{\bigcup_{n\geq m}\bigg\{\frac{\Psi
_f^n(p)}n:p\in S\bigg\}}. 
$$

In particular, since $\Psi _f$ is bounded, the rotation set is compact.

Note that, using a computation similar to (\ref{computation}), if one
chooses a rectifiable arc $\beta$ joining $f^n(p)$ to $p$ one has

\begin{equation}
\label{aha!} [I^n_p * \beta] = [\gamma^{-1}_p * \alpha^n_p * \gamma_{f^n(p)}
* \beta] = \Psi^n_f(p) + [\gamma_{f^n(p)} * \beta * \gamma^{-1}_p]. 
\end{equation}

Thus, $||I_p^n*\beta -\Psi _f^n(p)||\leq 2C_{\mathcal{A}_b}+l(\beta ).$ As a
consequence, an alternate but equivalent definition of rotation vectors and
rotation sets is obtained by considering all limits of the form 
$$
v=\lim _{k\to \infty }\frac 1{n_k}[I_{p_k}^{n_k}*\beta _k], 
$$
where $p_k\in S$, $n_k\to \infty $ and $\beta _k$ are rectifiable arcs
joining $f^{n_k}(p_k)$ to $p_k$ such that $l(\beta _k)<\infty $.

Moreover, it is possible to choose the arcs $\gamma_p$ in the definition of $%
\Psi_f$ so that the map $p \mapsto \Psi_f$ is not only bounded, but also
Borel measurable \cite{jf96}.

This is important if one wants to define rotation vectors of invariant
measures. Let $\mathcal{M}(f)$ be the set of all $f$-invariant Borel
probability measures. The rotation vector of the measure $\mu \in 
\mathcal{M}(f)$ is defined as 
$$
\rho _m(f,\mu )=\int \Psi _fd\mu \in H_1(S,\mathbb{R}). 
$$

By the Birkhoff ergodic theorem, for $\mu $-almost every point $p\in S$ the
limit $\rho (f,p)=\lim _{n\to \infty }\frac 1n\Psi _f^n(p)$ exists and $\rho
_m(f,\mu )=\int \rho (f,p)d\mu $. Moreover, if $\mu $ is an ergodic measure,
then $\rho (f,p)=\rho _m(f,\mu )$ for $\mu $-almost every point $p$.

Due to these facts and (\ref{bound}), the rotation vector of a measure is
also independent of any choices made in the definitions. Denote by $\rho
_m(f)$ the rotation set of invariant measures, that is $\rho _m(f)=\cup
_{\mu \in \mathcal{M}(f)}\rho _m(f,\mu )$ and $\rho _{erg}(f)$ the
corresponding set for ergodic measures. The proof of theorem 2.4 of \cite
{mmkz89}, without modifications, implies that 
$$
\rho _m(f)=Conv(\rho _{erg}(f))=Conv(\rho _{mz}(f)). 
$$

In particular, every extremal point of the convex hull of $\rho _{mz}(f)$ is
the rotation vector of some ergodic measure, and therefore, it is the
rotation vector of some recurrent point.

The main problems with this definition of rotation set are the following:

\begin{itemize}
\item  although it is compact, it is does not need to be convex;

\item  it is not known if vectors in the interior of the rotation set are
always realized by invariant sets;

\item  it is also not known if when $0$ is in the interior of the rotation
set, a result analogous to (\ref{c1eptoro}) holds, not even in the Abelian
cover of $S$ (see definition below).
\end{itemize}

\begin{description}
\item[Definition (Abelian cover)]  : Let $S$ be a closed orientable surface
of genus $g>1$. The Abelian cover of $S$ is a covering space for $S$, for
which the group of deck transformations is the integer homology group of $S$.
\end{description}

\vskip 0.2truecm

\subsection{Main motivation}

The main objective of this work is to give conditions which imply
complicated and rich dynamics in the universal cover of $S$, analogous to
what happens for a homeomorphism of the torus homotopic to the identity when
its rotation set contains $(0,0)$ in its interior.

This type of problem has already been studied for surfaces of higher genus
by P. Boyland in \cite{pb09}. But in that paper he considered the Abelian
cover of $S$ instead of the universal cover. As far as we know, this is the
only published result on this kind of problem. Boyland considered
homeomorphisms $f:S\to S$ of a special type, very important for our work: $f$
is isotopic to the identity as a homeomorphism of $S$, but it is
pseudo-Anosov relative to a finite $f$-invariant set $K\subset S$, see \cite
{afflvp}. He presented some conditions equivalent to $f$ having a transitive
lift to the Abelian cover of $S$.

The hypotheses of our main results will imply, in particular, that if a
homeomorphism $f:S\to S$ homotopic to the identity satisfies these
hypotheses, then there exists a constant $C_f\geq 0$ such that if $%
\widetilde f:\mathbb{D}\to \mathbb{D}$ is the only lift of $f$ to $\mathbb{D}
$ (the Poincar\'e disk) which commutes with all deck transformations of $S,$
for any fixed $\widetilde{Q}\subset \mathbb{D},$ a fundamental domain of $S,$
and for any deck transformation $g,$ there exist an integer $N>0$ which can
be chosen arbitrarily large, a compact set $K=K(N)\subset V_{C_f}(\widetilde{%
Q})$ (the open $C_f$-neighborhood of $\widetilde{Q}$ in the metric $d_{%
\mathbb{D}}$ of $\mathbb{D}$, which is the lift of the hyperbolic metric $d$
in $S),$ such that 
$$
\widetilde f^N(K)=g(K). 
$$

In particular, this implies that $\cup _{n>0}\widetilde{f}^n(V_{C_f}(
\widetilde{Q}))$ accumulates in the whole boundary of $\mathbb{D}$ and given
any compact set $M\subset \mathbb{D}$, if 
$$
D_M\stackrel{def.}{=}\{g(V_{C_f}(\widetilde{Q})):g\text{ is some deck
transformation for which }g(V_{C_f}(\widetilde{Q}))\ \cap M\ \neq \emptyset
\}, 
$$
then there exists $N_M>0$ such that for all $n\geq N_M,$ $\widetilde{f}%
^n(V_{C_f}(\widetilde{Q}))$ intersects all expanded fundamental domains
contained in $D_M.$ By expanded fundamental domains, we mean translates of $%
V_{C_f}(\widetilde{Q})$ by deck transformations.

In case $f$ is $C^{1+\epsilon },$ theorem 2 clearly shows that $C_f=0$ is
enough.

In the torus case, if $(0,0)$ belongs to the interior of the rotation set,
an analogous property holds. Therefore, in the situation when the
fundamental group is not Abelian (surfaces of genus larger than $1$), our
hypotheses, the fully essential system of curves $\mathscr{C}$ (see
definition \ref{thesystem}), are an analog for $(0,0)$ being in the interior
of the rotation set when the surface is the torus.

As a by-product of our main results, we obtain that in the $C^{1+\epsilon }$
setting, the homological rotation set is a compact convex subset of $%
\mathbb{R}^{2g},$ $2g$-dimensional: it is equal to the rotation set of the $%
f $-invariant Borel probability measures and all rational points in its
interior are realized by periodic orbits. Non-rational points in the
interior of the rotation set are also realized by compact $f$-invariant sets.

We are indebted to Alejandro Passeggi, who pointed this consequence of
theorem 2 to us.

Moreover, as a corollary of the ideas used in this last result, we can
extend the main theorems from \cite{jlm} to our setting. This is done in
theorems 4 and 5.

\subsection{Precise statements of the main results}

In what follows, we precisely present the main results of this paper. Assume 
$S$ is a closed orientable surface of genus $g>1$ and $\pi : \widetilde{S}%
\rightarrow \widetilde{S}$ is its universal covering map. We may identify
the universal cover $\widetilde{S}$ with the Poincar\'e disk $\mathbb{D}$
and denote by $Deck(\pi )$ the groups of deck transformations of $S$.
Consider $f:S\to S$ a homeomorphism isotopic to the identity and let $%
\widetilde f:\mathbb{D}\to \mathbb{D}$ be the endpoint of the lift of the
isotopy from $Id$ to $f$ which starts at $Id:\mathbb{D}\to \mathbb{D}.$ We
call $\widetilde{f}$ the natural lift of $f$.

\begin{description}
\item[Definition]  \label{thesystem} {\bf (Fully essential system of curves $%
\mathscr{C}$) \ref{thesystem}:} We say that $f:S\rightarrow S$ is a
homeomorphism with a fully essential system of curves $\mathscr{C}$ if there
exist different closed geodesics $\gamma _1,\ldots ,\gamma _k$ in $S$, $%
k\geq 1$, such that $\left( \cup _{i=1}^k\gamma _i\right) ^c$ only has
non-essential connected components and for each $i\in \{1,\ldots ,k\}$,
there are $f$-periodic points $p_i^{-}$ and $p_i^{+}$ such that their
trajectories under the isotopy are closed curves freely homotopic to $\gamma
_i,$ or concatenations of $\gamma _i,$ with both possible orientations.
\end{description}

\begin{theorem}
\label{firstthm} Let $f:S\to S$ be a homeomorphism isotopic to the identity
with a fully essential system of curves $\mathscr{C}$ and $\widetilde{f}$ be
its natural lift. Then there exists a constant $C_f\geq 0$ such that for all 
$g\in Deck(\pi )$ and any fundamental domain $\widetilde{Q}\subset \mathbb{D}
$ of $S,$ there exist arbitrarily large natural numbers $N>0,$ a point $
\widetilde{r}=\widetilde{r}(N)\in \mathbb{D}$ and a compact set $%
K=K(N)\subset V_{C_f}(\widetilde{Q})$ such that 
$$
\widetilde{f}^N(\widetilde{r})=g(\widetilde{r})\text{ and }\widetilde{f}%
^N(K)=g(K). 
$$
\end{theorem}

\begin{description}
\item[Remark:]  If for $\widetilde{p}\in \mathbb{D},$ $\widetilde{f}^n(
\widetilde{p})=g(\widetilde{p})$ for some $g\in Deck(\pi ),$ then for every $%
h\in Deck(\pi ),$ $\widetilde{f}^n(h(\widetilde{p}))=hgh^{-1}(h(\widetilde{p}%
))$.
\end{description}

\vskip 0.2truecm

\begin{theorem}
\label{teo} For some $\epsilon >0$ let $f:S\to S$ be a $C^{1+\epsilon }$
diffeomorphism isotopic to the identity with a fully essential system of
curves $\mathscr{C}$ and let $\widetilde{f}$ be its natural lift. Then there
exists a contractible hyperbolic $f$-periodic saddle point $p\in S$ such
that for any $\widetilde{p}\in \pi ^{-1}(p)$ and for every $g\in Deck(\pi ),$
$$
W^u(\widetilde{p})\pitchfork W^s(g(\widetilde{p})). 
$$
\end{theorem}

\begin{description}
\item[Remark:]  A point $p\in S$ being contractible means that all $
\widetilde{p}\in \pi ^{-1}(p)$ are $\widetilde{f}$-periodic.
\end{description}

\vskip 0.2truecm

To prove this result we have to work with a pseudo-Anosov map $\phi $
isotopic to $f$ relative to the finite invariant set of periodic points
associated with the fully essential system of curves $\mathscr{C}$. Using
several properties of the stable and unstable foliations of this map, it is
possible to prove a result similar to theorem 2 for $\phi $ and then using
Handel's global shadowing \cite{mh85} and other technical results on Pesin
theory \cite{asmp03}, \cite{ak80} we can finally prove the theorem for the
original map $f.$ This procedure is similar to what was done in \cite{saz15}.

The main part of this paper is proving theorem 2 for relative pseudo-Anosov
maps and this is done in lemma \ref{lema2}.

The next results are consequences of theorem 2, exactly as in \cite{saz15}.
They all share the same hypotheses: Suppose for some $\epsilon >0,$ $f:S\to
S $ is a $C^{1+\epsilon }$ area preserving diffeomorphism isotopic to the
identity with a fully essential system of curves $\mathscr{C}.$\ 

\vskip0.2truecm

{\bf Corollary 1.} {\it If $f$ is transitive, then $f$ cannot have a
periodic open disk. In the general case, there exists $M=M(f)>0$ such that
if $D\subset S$ is a $f$-periodic open disk, then for any connected
component $\widetilde D$ of $\pi ^{-1}(D),$ $diam(\widetilde D)<M$ in the
metric $d_{\mathbb{D}}$.}

\vskip0.2truecm

In \cite{akft15}, it is proved that in case $f$ is just an area-preserving
homeomorphism of $S$ and the fixed
point set is inessential, then all $f$-invariant open disks have
diameter bounded by some constant $M>0.$ If moreover, for all $n>0,$ the set
of $n$-periodic points is inessential, then for each $n>0,$ the set of
$n$-periodic open disks has bounded diameter. But the bound may not be uniform
with the period. In our situation, with  much stronger hypotheses,
corollary 1 gives a uniform bound.
         
\vskip0.2truecm

{\bf Corollary 2.} {\it There exists a contractible hyperbolic $f$-periodic
saddle point $p\in S$ (the one from theorem 2) such that $R.I.(f)\stackrel{%
def.}{=}\overline{W^u(p)}=\overline{W^s(p)},$ is compact, $f$-invariant and
all connected components of the complement of $R.I.(f)$ are $f$-periodic
disks. Moreover, for all $\widetilde{p}\in \pi ^{-1}(p),$ $R.I.(\widetilde f)%
\stackrel{def.}{=}\pi ^{-1}(R.I.(f))=\overline{W^s(\widetilde p)}=\overline{%
W^u(\widetilde p)}$ is a connected, closed, $\widetilde f$-invariant,
equivariant subset of $\mathbb{D}.$}

\vskip0.2truecm

Following the notation of \cite{akft15}, under our hypotheses,
their set of all dynamically fully essential points is equal to $R.I.(f)$

\vskip0.2truecm

{\bf Corollary 3.} {\it If $f$ is transitive then there exists a
contractible hyperbolic $f$-periodic saddle point $p\in S$ (the one from
theorem 2) such that $\overline{W^u(p)}=\overline{W^s(p)}=S$ and for any $%
\widetilde p\in \pi ^{-1}(p),$ $\overline{W^u(\widetilde p)}=\overline{%
W^s(\widetilde p)}=\mathbb{D},$ something that implies that $\widetilde{f}$
is topologically mixing.} \vskip0.2truecm

Finally, in the third theorem we study the homological rotation set $\rho
_{mz}(f).$

\begin{theorem}
\label{rot} Let $f:S\to S$ be a $C^{1+\epsilon }$ diffeomorphism isotopic to
the identity with a fully essential system of curves $\mathscr{C}$. Then the
(homological) rotation set $\rho _{mz}(f)$ is a $2g$-dimensional compact
convex subset of $H_1(S,\mathbb{R})\simeq \mathbb{R}^{2g}$. Moreover, if $%
v\in int(\rho _{mz}(f)),$ then there exists a compact set $K\subset S$ such
that for all $q\in K,$ $\rho (f,q)=v.$ In case $v$ is a rational point, $K$
can be chosen as a periodic orbit.
\end{theorem}

The last two results generalize the main theorems of \cite{jlm} to the
context of this paper.

\begin{theorem}
\label{bounddisp} Let $f:S\to S$ be a $C^{1+\epsilon }$ diffeomorphism
isotopic to the identity with a fully essential system of curves $\mathscr{C}
$. Then there exists $M(f)>0,$ such that for any $\omega \in \partial \rho
_{mz}(f),$ any hyperplane $\omega \in H\subset \mathbb{R}^{2g},$ which does
not intersect $interior(\rho _{mz}(f))$ ($H$ is called a supporting
hyperplane), any $p\in S$ and $n>0,$ 
$$
\left( [\alpha _p^n]-n.\omega \right) .\overrightarrow{v_H}<M(f), 
$$
where $\overrightarrow{v_H}$ is the unitary normal to $H,$ which points
towards the connected component of $H^c$ that does not intersect $\rho
_{mz}(f).$
\end{theorem}

\begin{theorem}
\label{boyland} Let $f:S\to S$ be a $C^{1+\epsilon }$ area-preserving
diffeomorphism isotopic to the identity with a fully essential system of
curves $\mathscr{C}$. Then the rotation vector of Lebesgue measure belongs
to $interior(\rho _{mz}(f)).$
\end{theorem}

\section{Some background, auxiliary results and their proofs}

In this section we present some important results we use, along with some
definitions and a short digression on hyperbolic surfaces, Thurston
classification of homeomorphisms of surfaces up to isotopy and a little of
Pesin theory. We also prove some auxiliary results and using them, in the
next sections, we prove theorems 1, 2, 3, 4 and 5.

\subsection{Properties of hyperbolic surfaces}

Let $S$ be a closed orientable surface of genus $g>1$ and $\pi : \widetilde{S%
}\to \widetilde{S}$ its universal covering map. As we said before, the
universal cover $\widetilde{S}$ is identified with the Poincar\'e disk $%
\mathbb{D}$ endowed with the hyperbolic metric. Hence, we assume that $S=%
\mathbb{D}\slash\Gamma $, where $\Gamma $ is a cocompact freely acting group
of Moebius transformations. Any nontrivial deck transformation $g\in
Deck(\pi )=\Gamma $ is a hyperbolic isometry and extends to the ''boundary
at infinity'' $\partial \mathbb{D}$ as a homeomorphism which has exactly two
fixed points, one attractor and one repeller. These fixed points are the
endpoints of some $g$-invariant geodesic $\delta _g$ of $\mathbb{D}$, called
the axis of $g$. For any point $\widetilde p\in \overline{\mathbb{D}}$, the
sequence $g^n(\widetilde p)$ converges to one endpoint of $\delta _g$ as $%
n\to -\infty $ and to the other one as $n\to \infty $. Any subarc of $\delta
_g$ joining a point $\widetilde p$ to $g(\widetilde p),$ when projected to $%
S,$ becomes an essential loop $\gamma _g,$ which is the unique geodesic in
its free homotopy class.

Given an essential loop $\gamma : [0, 1] \to S$, an extended lift of $\gamma$
is an arc $\widetilde{\gamma} : \mathbb{R} \to \mathbb{D}$ obtained by the
concatenation of arcs that are the translation of a lift of $\gamma$ by all
iterates of some deck transformation. Two extended lifts of an essential
loop coincides if and only if they share the same endpoints in $\partial%
\mathbb{D}$. If $h$ is a deck transformation that commutes with $g$, then
the axis of $g$ is equal the axis of $h$, and the group of all deck
transformations that commute with $g$ is cyclic, generated by $g$ if $%
\gamma_g$ is in the free homotopy class of a simple loop.

Let $f:S\to S$ be a homeomorphism isotopic to the identity and let $%
I:[0,1]\times S\to S$ be an isotopy from the identity map to $f$. The
isotopy $\widetilde{I}:[0,1]\times \mathbb{D}\to \mathbb{D}$ obtained by
lifting $I$ with basepoint $Id:\mathbb{D}\to \mathbb{D}$ is called the
natural lift of $I$. As we already defined in subsection 1.4, the map $
\widetilde{f}:\mathbb{D}\to \mathbb{D},$ given by $\widetilde{f}(\widetilde{p%
})=\widetilde{I}(1,\widetilde{p}),$ is called the natural lift of $f$
associated with the isotopy $I$. Natural lifts of a homeomorphism are
characterized by the property of commuting with all deck transformations,
and moreover, $\widetilde{f}$ can be extended to a homeomorphism of $
\overline{\mathbb{D}},$ as the identity on the ''boundary at infinity'' $%
\partial \mathbb{D}$, see \cite{bfdm}.

\subsection{On the fully essential system of curves $\mathscr{C}$}

Let us recall the definition of a fully essential system of curves $%
\mathscr{C}$ as in \ref{thesystem}. It is equivalent to the following: for
some $k\geq 1$, $f$ has a set of periodic points $P=\{p_1^{+},p_1^{-},\ldots
,p_k^{+},p_k^{-}\}$ satisfying:

\begin{enumerate}
\item  For each pair of points $p_i^{+}$, $p_i^{-}$, $1\leq i\leq k$, there
exist lifts $\widetilde{p}_i^{+}\in \pi ^{-1}(p_i^{+})$, $\widetilde{p}%
_i^{-}\in \pi ^{-1}(p_i^{-})$, and deck transformations $g_i^{+}$, $g_i^{-}$
such that $\widetilde{f}^{n_i^{+}}(\widetilde{p}_i^{+})=g_i^{+}(\widetilde{p}%
_i^{+})$ and $\widetilde{f}^{n_i^{-}}(\widetilde{p}_i^{-})=g_i^{-}(
\widetilde{p}_i^{-})$, for certain natural numbers $n_i^{+}$ and $n_i^{-};$

\item  For all $\widetilde{p}\in \overline{\mathbb{D}}$, $\lim _{n\to \infty
}(g_i^{+})^n(\widetilde{p})=\lim _{n\to \infty }(g_i^{-})^{-n}(\widetilde{p})
$ and \\ $\lim _{n\to \infty }(g_i^{-})^n(\widetilde{p})=\lim _{n\to \infty
}(g_i^{+})^{-n}(\widetilde{p});$

\item  Defining $\gamma _i=\pi (\delta _{g_i^{+}})=\pi (\delta _{g_i^{-}}),$
then $\gamma _i\neq \gamma _{i^{\prime }}$ when $i\neq i^{\prime };$

\item  If $\mathscr{C}=\cup _{i=1}^k\gamma _i$, then $S\setminus \mathscr{C}$
is a union of open topological disks;
\end{enumerate}

\begin{proposition}
\label{connectedC} The lift $\pi ^{-1}(\mathscr{C})$ is a closed connected
subset of $\mathbb{D}$.
\end{proposition}

{\it Proof: }First observe that $\mathscr{C}$ is the union of a finite
number of closed geodesics in $S$, therefore $\mathscr{C}$ is closed. Since $%
\pi :\mathbb{D}\to S$ is continuous, $\pi ^{-1}(\mathscr{C})$ is closed. To
see that $\pi ^{-1}(\mathscr{C})$ is connected we just observe that $%
S\setminus \mathscr{C}$ is a union of open topological disks, and therefore
all connected components of $\mathbb{D}\setminus \pi ^{-1}(\mathscr{C})$ are
bounded topological open disks. \qed

\vskip 0.2truecm

\begin{proposition}
\label{finitepath} For every $\widetilde{p},\widetilde{r}\in \pi ^{-1}(%
\mathscr{C})$, there exists a path $\gamma $ in $\pi ^{-1}(\mathscr{C})$
joining these two points, contained in the union of a finite number of
subarcs of extended lifts of geodesics in $\mathscr{C}$. \ 
\end{proposition}

{\it Proof: }Fix a point $\widetilde p\in \pi ^{-1}(\mathscr{C})$ and let $%
\mathcal{P}_{\widetilde p}$ be the set of all points $\widetilde q\in \pi
^{-1}(\mathscr{C})$ such that there exists a path joining $\widetilde p$ to $%
\widetilde q$ formed by subarcs of a finite number of extended lifts of
geodesics in $\mathscr{C}$. We will show that $\mathcal{P}_{\widetilde p}$
is an open and closed subset of $\pi ^{-1}(\mathscr{C})$.

\begin{figure}[!h]
	\centering
	\includegraphics[scale=0.23]{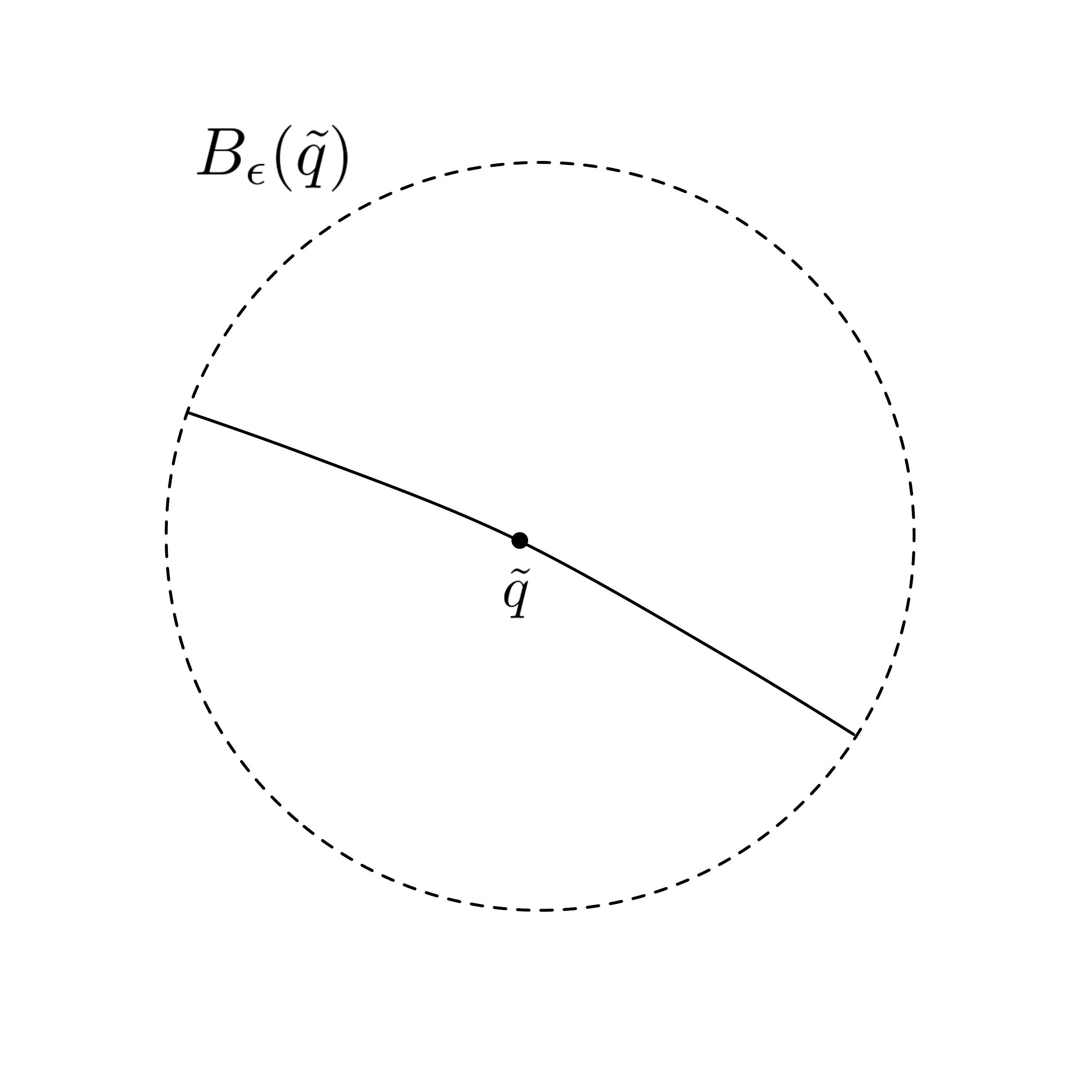}
	\includegraphics[scale=0.23]{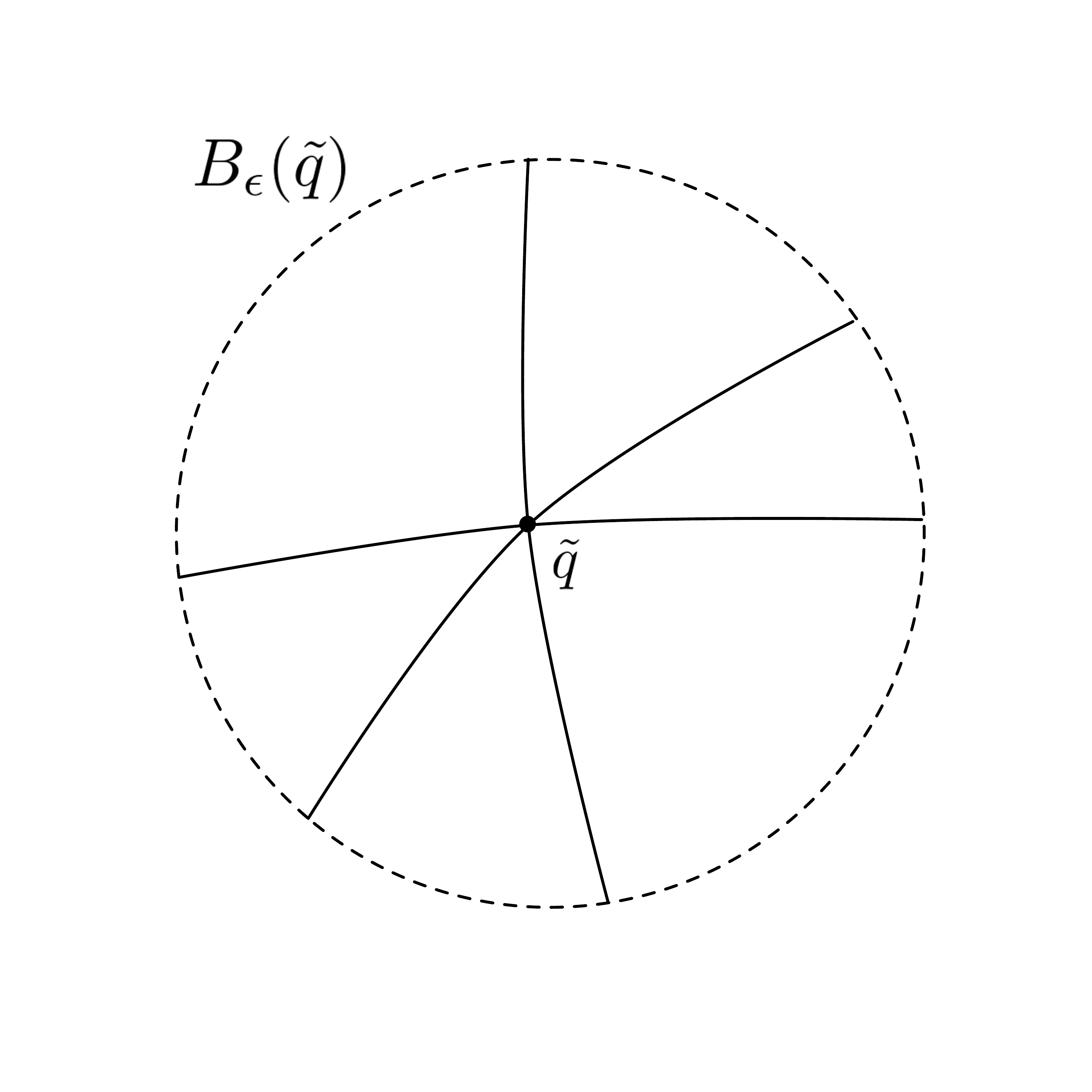}
	\caption{Possibilities for a neighborhood of $\widetilde q.$}
	\label{clopen}
\end{figure}

Let $\widetilde q$ be a point in $\mathcal{P}_{\widetilde p}$. As the set $%
\mathscr{C}$ is equal the union of a finite number of closed geodesics,
there exists $\epsilon >0$ small enough so that $B_\epsilon (\widetilde
q)\cap \pi ^{-1}(\mathscr{C})$ satisfies one of the possibilities in figure 
\ref{clopen}.

In the first case $\widetilde q$ belongs to just one extended lift of a
geodesic in $\mathscr{C}$. If $\gamma $ is the path joining $\widetilde p$
to $\widetilde q$ and it is formed by $k>0$ subarcs of extended lifts of
geodesics, it is clear that for all points in $B_\epsilon (\widetilde q)\cap
\pi ^{-1}(\mathscr{C}),$ there is a path $\gamma ^{\prime }$ joining $%
\widetilde p$ to this point formed by the same number of subarcs of extended
lifts of geodesics. In the second case $\widetilde q$ belongs to the
intersection of a finite number of extended lifts of geodesics, and again,
if the path $\gamma $ is formed by $k>0$ subarcs, then for all points in $%
B_\epsilon (\widetilde q)\cap \pi ^{-1}(\mathscr{C})$ there is a path $%
\gamma ^{\prime }$ joining $\widetilde p$ to this point formed by at most $%
k+1$ subarcs of extended lifts of geodesics. So $\mathcal{P}_{\widetilde p}$
is open.

We will prove now that $\mathcal{P}_{\widetilde{p}}^c=\pi ^{-1}(\mathscr{C}%
)\setminus \mathcal{P}_{\widetilde{p}}$ is open. Again, if $\widetilde{q}$
is a point in $\mathcal{P}_{\widetilde{p}}^c,$ there exists $\epsilon >0$
small enough such that $B_\epsilon (\widetilde{q})\cap \pi ^{-1}(\mathscr{C})
$ satisfies one of the possibilities in figure \ref{clopen}. In both cases,
if $\widetilde{q}^{\prime }\in B_\epsilon (\widetilde{q})\cap \pi ^{-1}(%
\mathscr{C})\cap \mathcal{P}_{\widetilde{p}}$, then by the same argument as
above, there is a path $\gamma ^{\prime }$ joining $\widetilde{p}$ to $
\widetilde{q}$ with a finite number of subarcs of extended lifts of
geodesics. But this is a contradiction because $\widetilde{q}\in \mathcal{P}%
_{\widetilde{p}}^c,$ so all points in $B_\epsilon (\widetilde{q})\cap \pi
^{-1}(\mathscr{C})$ are points of $\mathcal{P}_{\widetilde{p}}^c$. Hence $%
\mathcal{P}_{\widetilde{p}}^c$ is open. Since $\mathcal{P}_{\widetilde{p}}$
is an open and closed subset of the connected set $\pi ^{-1}(\mathscr{C}),$ $%
\mathcal{P}_{\widetilde{p}}=\pi ^{-1}(\mathscr{C}).$ \qed

\subsection{{}Nielsen Thurston classification of homeomorphisms of surfaces}

In this subsection we present a brief overview of Thurston's classification
of homeomorphisms of surfaces and prove a result analogous to one due to
Llibre and Mackay \cite{jlrm91}.

\subsubsection{Some definitions and the classification theorem}

Let $M$ be a compact, connected, orientable surface, possibly with boundary,
and let $f:M\to M$ be a homeomorphism. There are two basic types of
homeomorphisms which appear in the Nielsen-Thurston classification: the
finite order homeomorphisms and the pseudo-Anosov ones.

A homeomorphism $f$ is said to be of finite order if $f^n = Id$ for some $n
\in \mathbb{N}$. The least such $n$ is called the order of $f$. Finite order
homeomorphisms have zero topological entropy.

A homeomorphism $f$ is said to be pseudo-Anosov if there is a real number $%
\lambda > 1$ and a pair of transverse measured foliations $\mathscr{F}^u$
and $\mathscr{F}^s$ such that $f(\mathscr{F}^s) = \lambda^{-1} \mathscr{F}^s$
and $f(\mathscr{F}^u) = \lambda \mathscr{F}^u$. Pseudo-Anosov homeomorphisms
are topologically transitive, have positive topological entropy and Markov
partitions \cite{afflvp}.

A homeomorphism $f$ is said to be reducible by a system 
$$
C = \bigcup^{n}_{i = 1} C_i%
$$
of disjoint simple closed curves $C_1, \ldots, C_n$, called reducing curves,
if:

\begin{itemize}
\item  $\forall i$, $C_i$ is not homotopic to a point, nor to a component of 
$\partial M$,

\item  $\forall i\neq j$, $C_i$ is not homotopic to $C_j$,

\item  $C$ is invariant under $f$.
\end{itemize}

\begin{theorem}[Nielsen-Thurston]
\label{nielsenthurston} If the Euler characteristic $\chi (M)<0$, then every
homeomorphism $f:M\to M$ is isotopic to a homeomorphism $\phi :M\to M$ such
that either:

\begin{enumerate}
\item  $\phi $ is of finite order; \label{finiteorder}

\item  $\phi $ is pseudo-Anosov; \label{pseudoanosov}

\item  $\phi $ is reducible by a system of curves $C$, and there exist
disjoint open annular neighborhoods $U_i$ of $C_i$ such that 
$$
U=\bigcup_iU_i 
$$
is $\phi $-invariant. Each component $S_i$ of $M\setminus U$ is mapped to
itself by some least positive iterate $n_i$ of $\phi $, and each $\phi
^{n_i}|_{S_i}$ satisfies $\ref{finiteorder}$ or $\ref{pseudoanosov}$. Each $%
U_i$ is mapped to itself by some least positive iterate $m_i$ of $\phi $
fixing the boundary components, and each $\phi ^{m_i}|_{U_i}$ is a
generalized twist.
\end{enumerate}
\end{theorem}

Homeomorphisms $\phi$ as in theorem \ref{nielsenthurston} are called
Thurston canonical forms for $f$.

\begin{figure}[!h]
	\centering
	\includegraphics[scale=0.13]{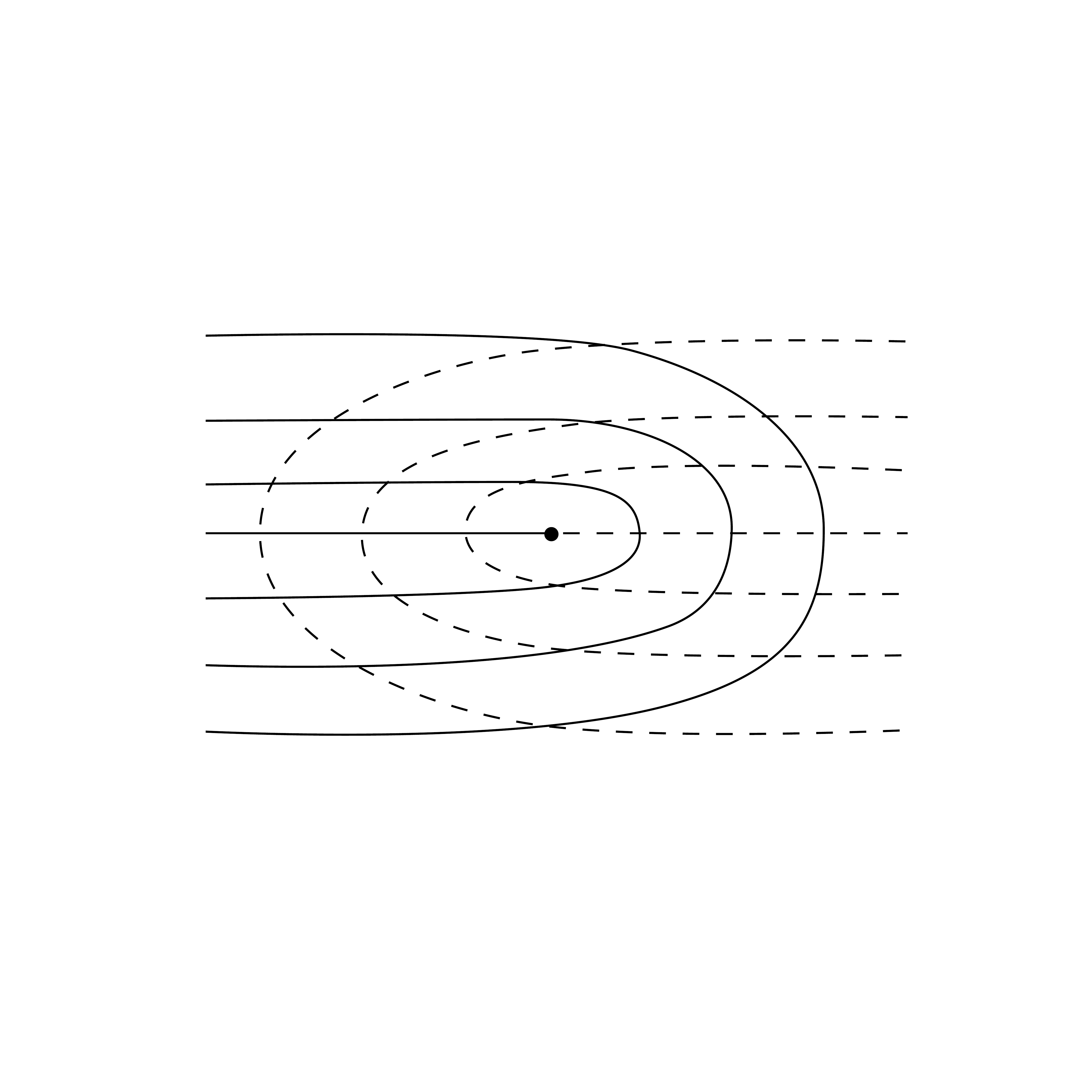}
	\includegraphics[scale=0.13]{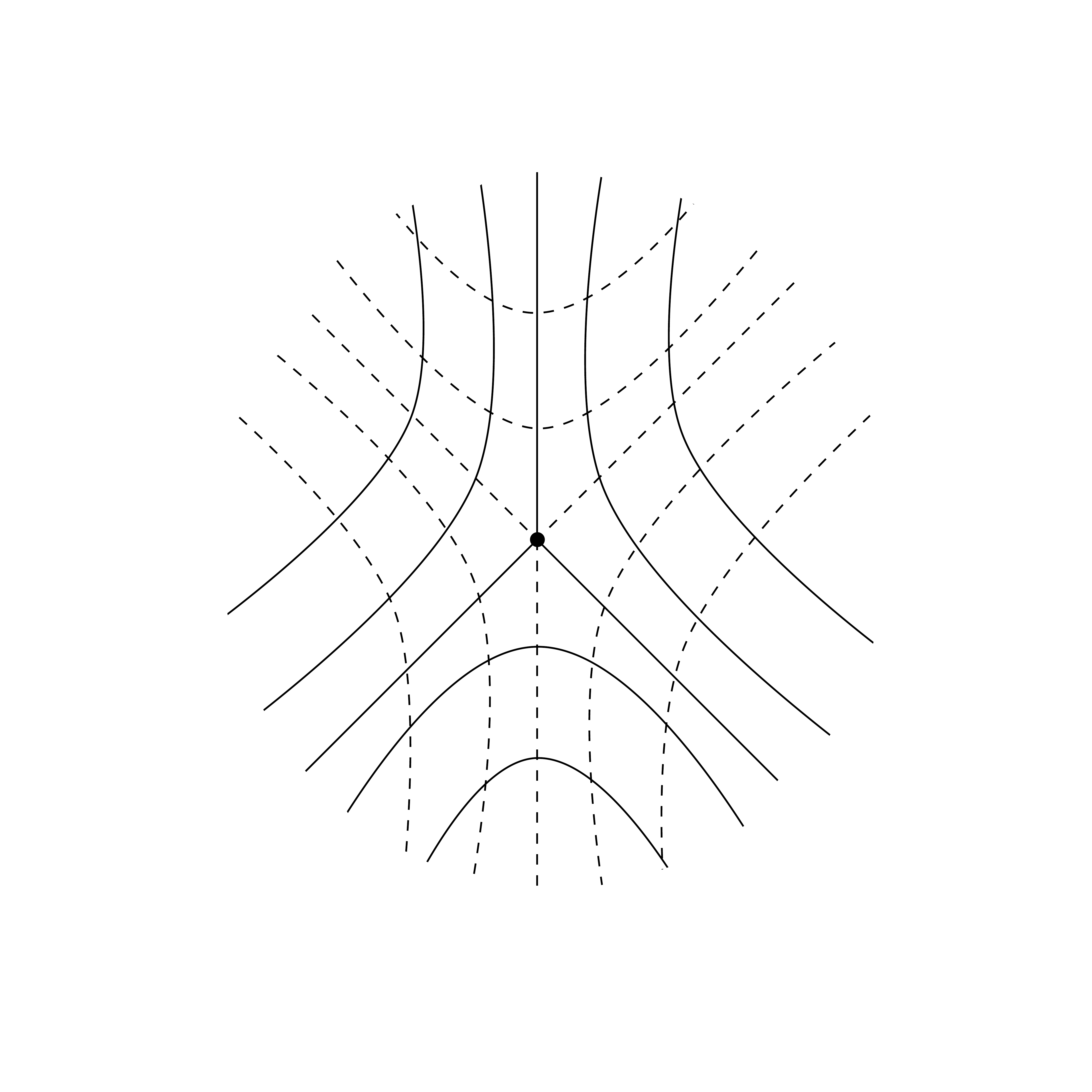}
	\caption{Examples of a $1$-prong and a $3$-prong singularity respectively.}
\end{figure}

We say that $\phi : M \to M$ is pseudo-Anosov relative to a finite invariant
set $K$ if it satisfies all of the properties of a pseudo-Anosov
homeomorphism except that the associated stable and unstable foliations may
have $1$-pronged singularities at points in $K$ \cite{mh90}. Equivalently,
let $N$ be the compact surface obtained from $M \setminus K$ by
compactifying each puncture with a boundary circle, let $p : N \to M$ be the
map that collapses these boundary circles to points. Then $\phi$ is
pseudo-Anosov relative to $K$, if and only if there is a pseudo-Anosov
homeomorphism $\Phi : N \to N$ such that $\phi \circ p = p \circ \Phi.$ %

\subsubsection{The beginning of the work}

The following result is the first step towards the proof of the main
theorems.

\begin{lemma}
\label{parel} Let $f:S\to S$ be a homeomorphism isotopic to the identity
with a fully essential system of curves $\mathscr{C}$ and let $P$ be the set
of periodic points associated with the geodesics in $\mathscr{C}$. Then,
there exists an integer $m_0>0$ such that $f^{m_0}$ is isotopic relative to $%
P$ to $\phi :S\to S$, a homeomorphism which is pseudo-Anosov relative to $P$.
\end{lemma}

{\it Proof: }Let $f$ be a homeomorphism with a fully essential system of
curves $\mathscr{C}$ and $P$ be the set of all periodic points associated
with the geodesics in $\mathscr{C}$. We will write $P=\{p_1^{+},p_1^{-},%
\ldots ,p_k^{+},p_k^{-}\}$. For each $1\leq i\leq k$, there exists integers $%
n_i^{+},n_i^{-}>0$ such that $f^{n_i^{+}}(p_i^{+})=p_i^{+}$ and $%
f^{n_i^{-}}(p_i^{-})=p_i^{-}$. Take $m_0>0$ an integer such that all points
in $P$ are fixed points for $f^{m_0}$.

We will follow the same ideas used by Llibre and MacKay in \cite{jlrm91}.
Let $\phi :S\rightarrow S$ be Thurston canonical form associated to $%
f^{m_0}. $ Of course we are considering $f^{m_0}:S\backslash P\rightarrow
S\backslash P$ and so, $\phi $ is also a homeomorphism from $S\backslash P$
into itself. But it can be extended in a standard way to the set $P,$ giving
a homeomorphism of $S$ into itself, also homotopic to the identity as a
homeomorphism of $S,$ which we still call $\phi .$

Let us show that $\phi $ is pseudo-Anosov relative to $P.$ First note that $%
\phi $ can not be of finite order, since points in $P$ move in different
homotopical directions. As $k\geq 1,$ there are at least 2 fixed points for $%
\phi,$ $p_1^{+}$ and $p_1^{-},$ whose trajectories under the natural lift $
\widetilde{\phi}: \mathbb{D} \rightarrow \mathbb{D}$ move under different
deck transformations.

Now, suppose $\phi $ is reducible by a system of curves $C$. As in \cite
{jlrm91}, say a simple closed curve $\gamma $ on a surface of genus $g$ with
punctures is non-rotational if after closing the punctures $\gamma $ is
homotopically trivial. If $\gamma $ is a non-rotational reducing curve, then
it must surround at least two punctures. Without loss of generality, we may
suppose that $\gamma $ surrounds $p_i^{+}$ and $p_j^{+}$, $i\neq j$. Since $%
\gamma $ is a reducing curve, then $\phi ^n(\gamma )=\gamma $, for some $n>0$%
. This means that exists $g\in Deck(\pi )$ such that $\widetilde \phi
^n(\widetilde \gamma )=g(\widetilde \gamma )$, where $\widetilde \gamma $ is
a lift of $\gamma $ $(\widetilde \gamma $ is a simple closed curve in $%
\mathbb{D})$ surrounding $\widetilde p_i^{+}$ and $\widetilde p_j^{+}$,
lifts of $p_i^{+}$ and $p_j^{+}$ respectively. By induction it follows that $%
\widetilde \phi ^{mn}(\widetilde \gamma )=g^m(\widetilde \gamma )$ encloses
both $\widetilde \phi ^{mn}(\widetilde p_i^{+})$ and $\widetilde \phi
^{mn}(\widetilde p_j^{+})$ for all $m\in \mathbb{Z}$. But this is a
contradiction because as $p_i^{+}\neq p_j^{+},$ $\lim _{l\to \infty
}\widetilde \phi ^l(\widetilde p_i^{+})$ and $\lim _{l\to \infty }\widetilde
\phi ^l(\widetilde p_j^{+})$ are different points of $\partial \mathbb{D}$.

Now, if $\gamma $ is a rotational reducing curve, let $\widetilde \gamma
\subset \mathbb{D}$ be an extended lift of $\gamma $. The curve $\widetilde
\gamma $ has two distinct endpoints at the ''boundary at infinity'' $%
\partial \mathbb{D}$ and $\mathbb{D}\setminus \widetilde \gamma $ has
exactly two connected components. Since $\widetilde \phi |_{\partial 
\mathbb{D}}=Id$, $\widetilde \phi (\widetilde \gamma ) $ has the same
endpoints on $\partial \mathbb{D}$ as $\widetilde \gamma $. Since $%
S\setminus \mathscr{C}$ is a union of topological disks, there exists $g\in
Deck(\pi )$ associated with some geodesic $\gamma _i$ in $\mathscr{C}$ such
that the fixed points of $g$ in $\partial \mathbb{D}$ separate the endpoints
of $\widetilde \gamma $.

\begin{figure}[!h]
	\centering
	\includegraphics[scale=0.21]{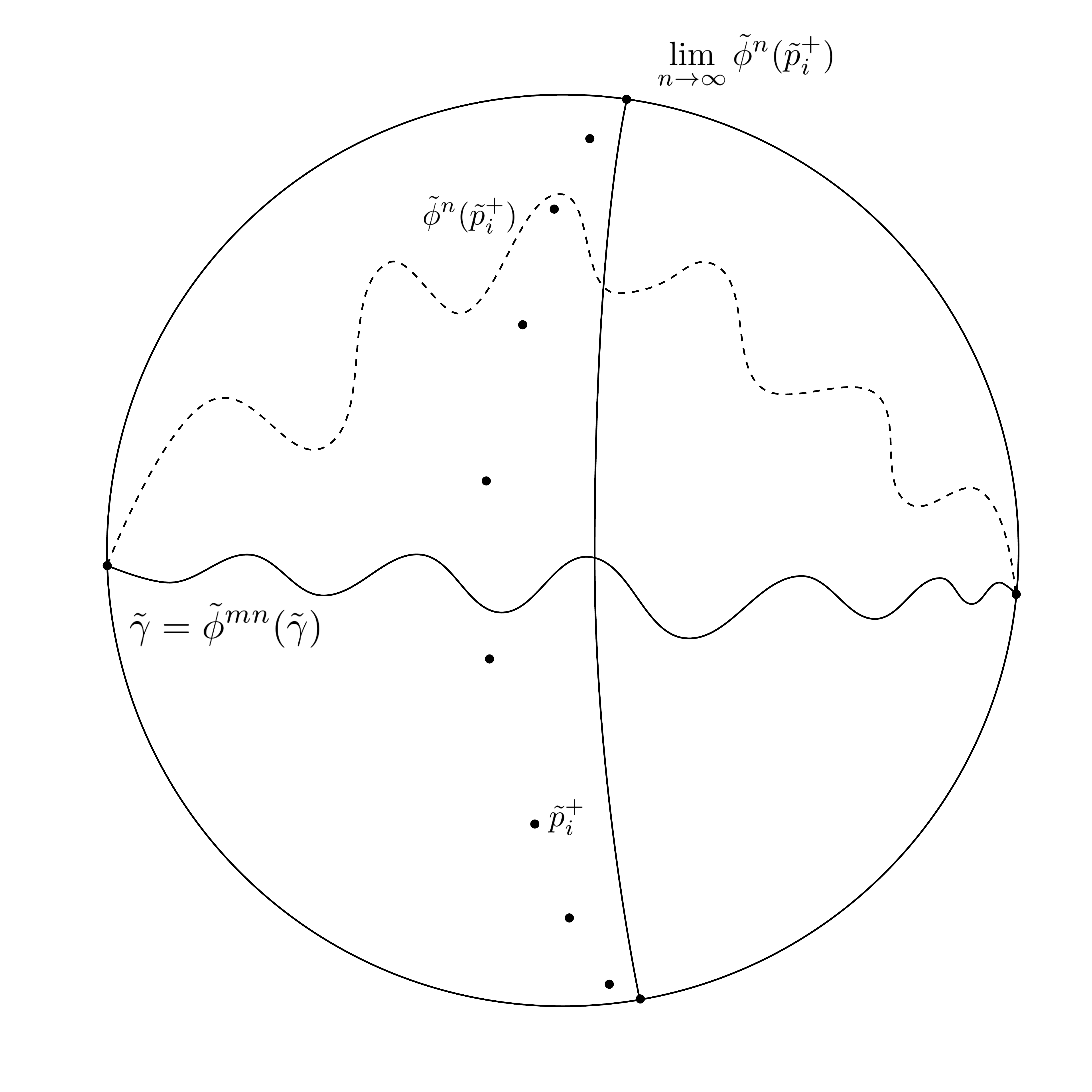}
	\caption{The final contradiction.}
	\label{contradiction}
\end{figure}

Finally, choose $\widetilde{p}_i^{+}\in \pi ^{-1}(p_i^{+})$ such that it
belongs to one connected component of $\mathbb{D}\setminus \widetilde{\gamma 
}$ and $\lim _{n\to \infty }\widetilde{\phi }^n(\widetilde{p}_i^{+})$ is in
the ''boundary at infinity'' of other connected component. Since $\widetilde{%
\phi }(\widetilde{\gamma })$ and $\widetilde{\gamma }$ have the same
endpoints in $\partial \mathbb{D}$ and $\phi ^n(\gamma )=\gamma ,$ then $
\widetilde{\phi }^{mn}(\widetilde{\gamma })=\widetilde{\gamma }$, for all $%
m>0$. As $\widetilde{\phi }$ preserves orientation, this clearly implies a
contradiction. See figure \ref{contradiction}.

This shows that $\phi $ cannot be of finite order or reducible by a system
of curves. So $\phi $ is pseudo-Anosov relative to $P.$

\qed

\subsection{On Handel's fixed point theorem}

\subsubsection{Preliminaries and a statement of Handel's theorem}

In \cite{han99} Michael Handel proved the existence of a fixed point for an
orientation preserving homeomorphism of the open unit disk that can be
extended to the closed disk as the identity on the boundary, provided that
for certain points in the open disk, their $\alpha $ and $\omega $-limit
sets are single points in the boundary of the disk, distributed with a
certain cyclic order. Later, in \cite{pat06}, Patrice Le Calvez gave a
different proof of this theorem based only on Brouwer theory and plane
topology arguments. In Le Calvez's proof, the existence of the fixed point
follows from the existence a simple closed curve contained in the open disk,
whose topological index can be calculated and is equal to $1.$

\begin{theorem}[Handel's fixed point theorem, \cite{pat06}]
\label{handelfixed}\index{Handel's fixed point theorem} Consider a
homeomorphism $\widetilde{h}:\overline{\mathbb{D}}\to \overline{\mathbb{D}}$
of the closed unit disk satisfying the following hypotheses:

\begin{enumerate}
\item  There exists $r\geq 3$ points $\widetilde{p}_1,\ldots ,\widetilde{p}_r
$ in $\mathbb{D}$ and $2r$ pairwise distinct points $\alpha _1,\omega
_1,\ldots ,\alpha _r,\omega _r$ on the boundary $\partial \mathbb{D}$ such
that, for every $1\leq i\leq r$, 
$$
\lim _{n\to \infty }\widetilde{h}^{-n}(\widetilde{p}_i)=\alpha _i,%
\hspace{4mm}\lim _{n\to \infty }\widetilde{h}^n(\widetilde{p}_i)=\omega _i. 
$$

\item  The cyclic order on $\partial \mathbb{D}$ is as represented on figure 
\ref{cyclic} below: 
$$
\alpha _1,\omega _r,\alpha _2,\omega _1,\alpha _3,\omega _2,\ldots ,\alpha
_r,\omega _{r-1},\alpha _1. 
$$
\end{enumerate}

Then there exists a fixed point free simple closed curve $\gamma \subset 
\mathbb{D}$ such that $ind(\widetilde{h},\gamma )=1$.
\end{theorem}

\begin{figure}[!h]
	\centering
	\includegraphics[scale=0.3]{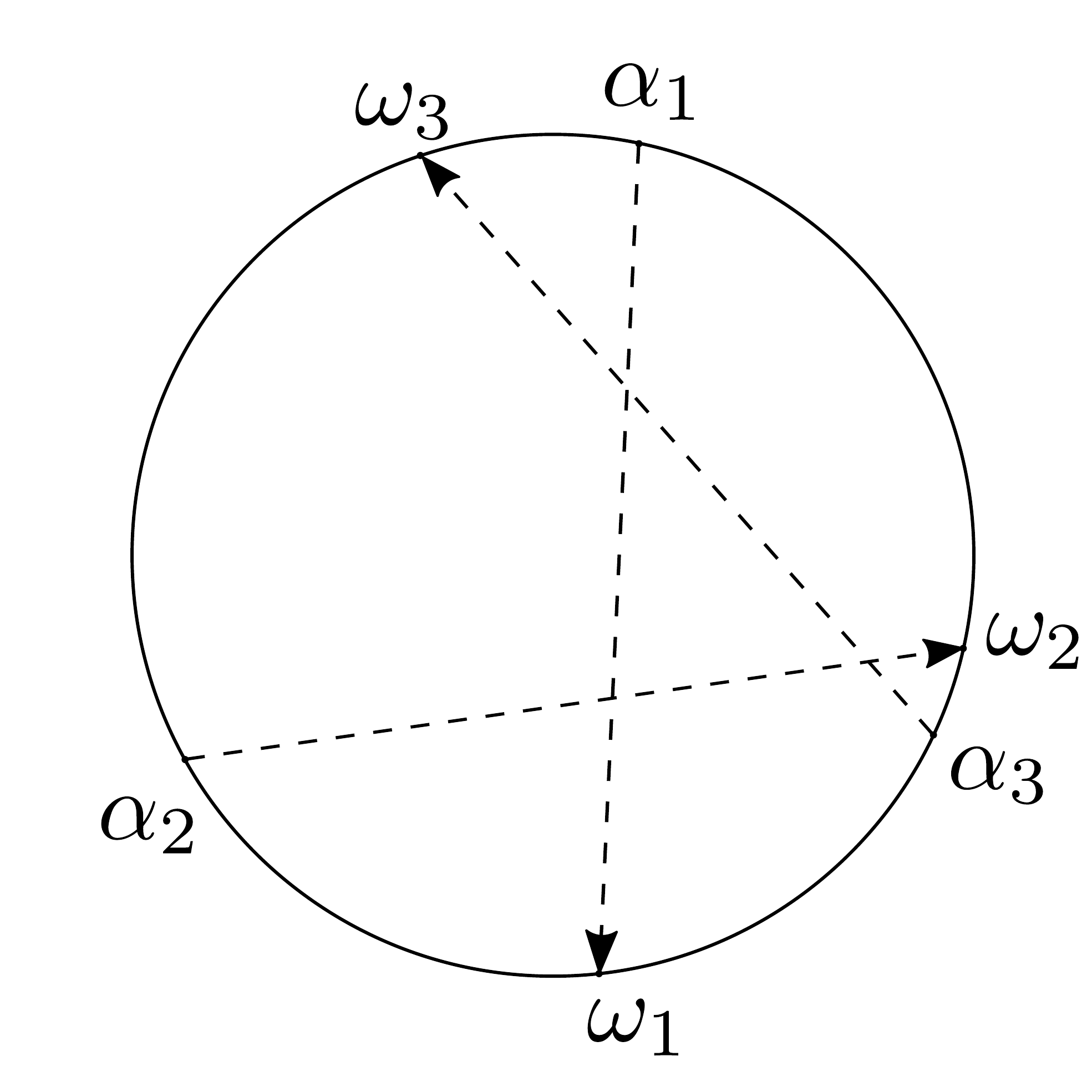}
	\hspace{10mm}
	\includegraphics[scale=0.3]{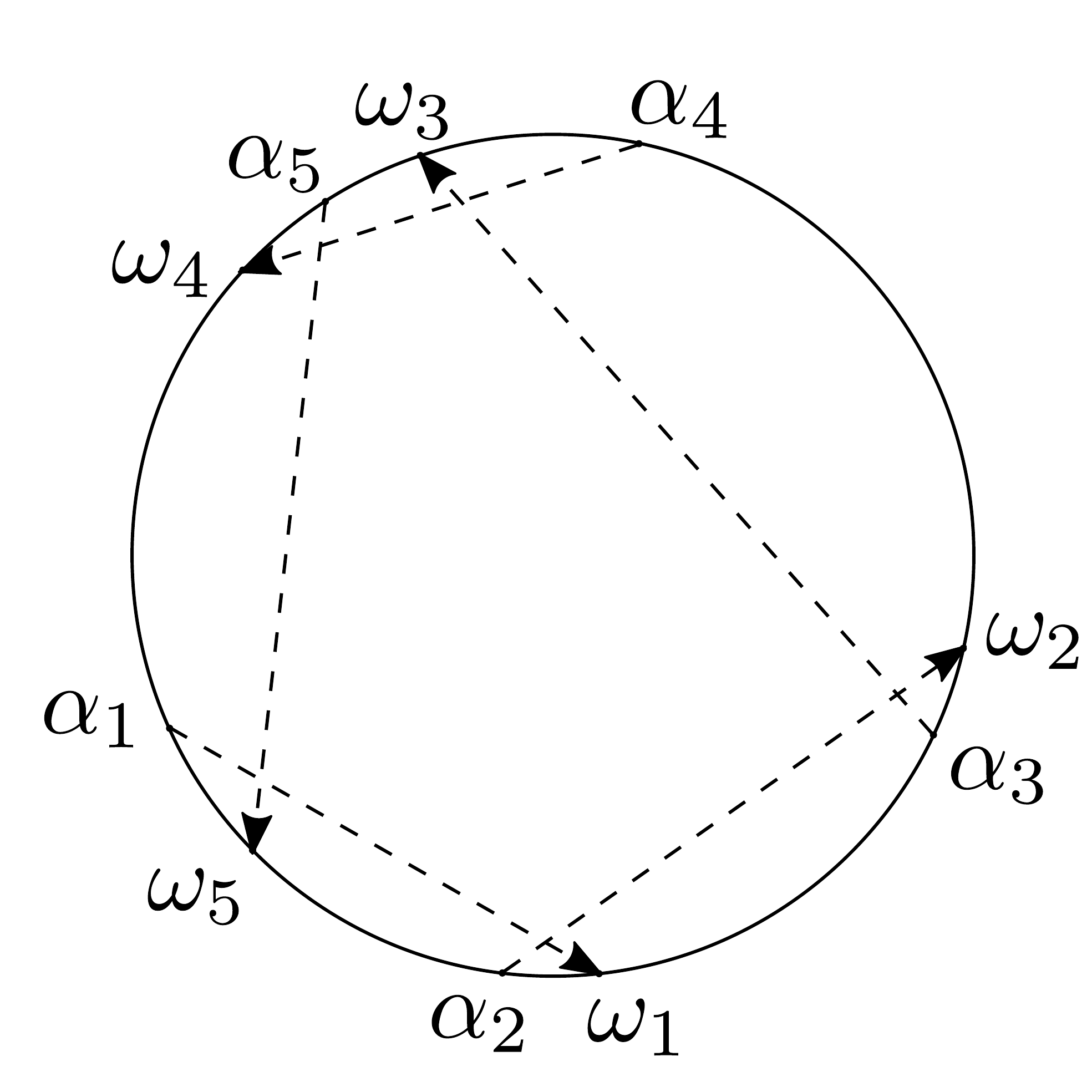}
	\caption{Cyclic order for Handel's fixed point theorem when $r = 3$ and $r = 5$.}
	\label{cyclic}
\end{figure}

Remember that, if $\widetilde p$ is an isolated fixed point of $\widetilde h$%
, the Poincar\'e-Lefschetz index of $\widetilde h$ at $\widetilde p$ is
defined as 
$$
ind(\widetilde h,\widetilde p)=ind(\widetilde h,\gamma ), 
$$
where $\gamma $ is a (small) simple closed curve surrounding $\widetilde p$
and no other fixed point. The index of $\widetilde h$ at $\widetilde p$ does
not depend of the choice of $\gamma .$

In case $\widetilde h$ has only isolated fixed points, if $int(\gamma )$ is
the bounded connected component of $\gamma ^c$ and $Fix(int(\gamma
))=\{\widetilde p\in int(\gamma ):\widetilde h(\widetilde p)=\widetilde p\}$
then, by properties of the Poincar\'e-Lefschetz index, 
$$
ind(\widetilde h,\gamma )=\sum_{\widetilde p\in Fix(int(\gamma
))}ind(\widetilde h,\widetilde p). 
$$

So if $\widetilde{h}:\mathbb{D}\rightarrow \mathbb{D}$ is a homeomorphism
with only isolated fixed points, satisfying the hypothesis of Handel's
theorem, as $ind(\widetilde{h},\gamma )=1,$ there exists a fixed point $
\widetilde{p}^{\prime }\in int(\gamma )$ with $ind(\widetilde{h},\widetilde{p%
}^{\prime })>0$.

\subsubsection{Existence of a hyperbolic $\widetilde{\phi }$-periodic point}

Remember that $\phi :S\rightarrow S$ is the relative to $P$ pseudo-Anosov
map given in lemma \ref{parel}, which is isotopic to the identity as a
homeomorphism of $S.$ $\widetilde{\phi }:\mathbb{D}\rightarrow \mathbb{D}$
is the natural lift of $\phi ,$ the one which commutes with all deck
transformations and extends as a homeomorphism of $\overline{\mathbb{D}},$
which is the identity on the ''boundary at infinity''. In the next
proposition we prove that $\widetilde \phi $ has a hyperbolic periodic
saddle point. When we say hyperbolic in this context, we mean that the local
dynamics at the point is obtained by gluing exactly four hyperbolic sectors,
or equivalently, the point is a regular point of the foliations $\mathscr{F}%
^u$ and $\mathscr{F}^s$.

\begin{proposition}
\label{periodic} The natural lift $\widetilde{\phi }:\mathbb{D}\to \mathbb{D}
$ of the map $\phi $ from lemma \ref{parel} has a hyperbolic periodic
(saddle) point $\widetilde{p}$.
\end{proposition}

{\it Proof: }We know that $\phi $ is a pseudo-Anosov map relative to a
finite set $P.$ And this set $P=\{p_1^{+},p_1^{-},\ldots ,p_k^{+},p_k^{-}\}$
is associated with a fully essential system of curves $\mathscr{C}.$ As a
consequence of that, for each geodesic $\gamma _j$ in $\mathscr{C}$ and
appropriate lifts $\widetilde{p}_j^{+}\in \pi ^{-1}(p_j^{+}),$ $\widetilde{p}%
_j^{-}\in \pi ^{-1}(p_j^{-})$, there are deck transformations $g_j^{+}$, $%
g_j^{-}$, with the following properties:

\begin{enumerate}
\item  $\widetilde{\phi }(\widetilde{p}_j^{+})=g_j^{+}(\widetilde{p}_j^{+})$
and $\widetilde{\phi }(\widetilde{p}_j^{-})=g_j^{-}(\widetilde{p}_j^{-});$

\item  $g_j^{+}\circ g_j^{-}=g_j^{-}\circ g_j^{+};$

\item  the invariant geodesic in $\mathbb{D}$ of the deck transformations $%
g_j^{-},g_j^{+}$ projects onto $\gamma _j,$ i. e. $\pi (\delta
_{g_j^{+}})=\gamma _j=\pi (\delta _{g_j^{-}});$
\end{enumerate}

So, by the above properties, 
$$
\lim _{n\to \infty }\widetilde \phi ^n(\widetilde p_j^{+})=\lim _{n\to
\infty }\widetilde \phi ^{-n}(\widetilde p_j^{-})\text{ and }\lim _{n\to
\infty }\widetilde \phi ^n(\widetilde p_j^{-})=\lim _{n\to \infty
}\widetilde \phi ^{-n}(\widetilde p_j^{+}), 
$$
where these limits are the fixed points of $g_j^{+}$ and $g_j^{-}$ in $%
\partial \mathbb{D}.$

Let $U$ be a connected component of $S\setminus \mathscr{C}$. By our
assumptions, $U$ is a topological disk and $\partial U$ is made by subarcs
of geodesics in $\mathscr{C}.$

If we consider $\widetilde{U}$ a connected component of the lift of $U,$ $%
\partial \widetilde U$ is made by lifts of subarcs of geodesics in $\partial
U$. Since $\partial U$ is homotopically trivial, $\partial \widetilde U$ is
a simple closed curve in $\mathbb{D}$. Fix $\widetilde \beta $ a subarc in $%
\partial \widetilde U$. Associated with this subarc there is an extended
lift of one geodesic in $\mathscr{C}$ and by the previous observations we
can find appropriate lifts $\widetilde p_{i_\beta }^{+}$, $\widetilde
p_{i_\beta }^{-},$ each one following the extended lift of the geodesic with
both possible orientations.

If $\widetilde{\beta }$ and $\widetilde{\beta }^{\prime }$ are two
consecutive subarcs of $\partial \widetilde{U}$, then the endpoints of the
extended lift of the geodesic associated to $\beta $ separate the endpoints
of the extended lift of the geodesic associated to $\beta ^{\prime }$.
Putting all these observations together we see that if we choose an
orientation for $\partial \widetilde{U}$, then it is possible to choose
lifts of points in $P$ such that $\widetilde{\phi }$ satisfies the
hypotheses of Handel's theorem.

\begin{figure}[!h]
	\centering
	\includegraphics[scale=0.21]{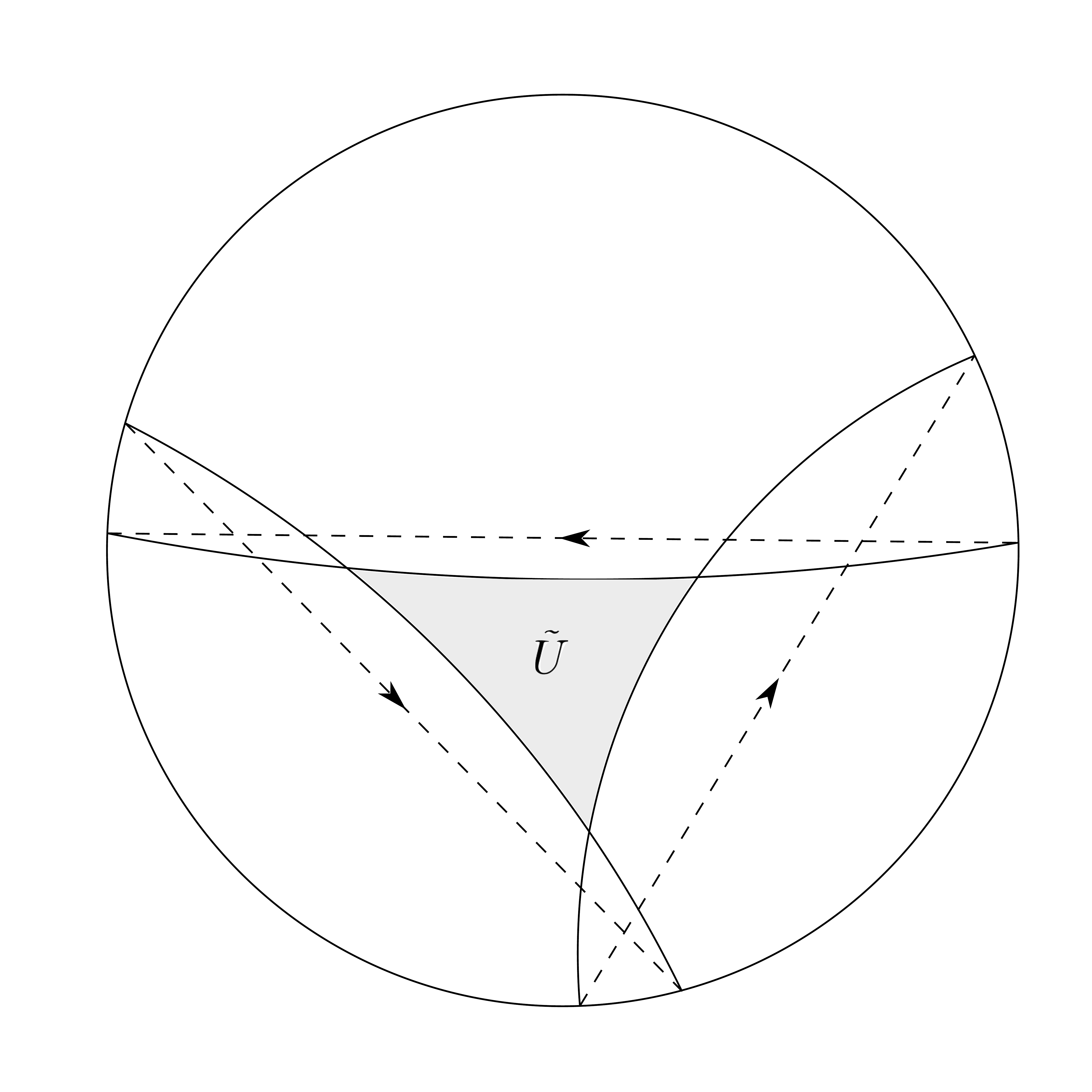}
	\caption{ $\widetilde U$ and how some points move with respect to its boundary.}
\end{figure}

Since $\phi $ is pseudo-Anosov relative to a finite set $P,$ for each
period, it has only isolated periodic points, and the same holds for $%
\widetilde \phi $. This means, by Handel's theorem, that there exists a
fixed point $\widetilde p_1$ of $\widetilde \phi $ such that%
$$
ind(\widetilde \phi ,\widetilde p_1)=ind(\phi ,\pi (\widetilde p_1))>0. 
$$
Observe that the same conclusion holds for $\widetilde \phi ^m$, for any $%
m>0 $.

But for some appropriate large $m_1>0,$ the local dynamics at points in $%
Fix(\phi )$ implies that 
$$
ind(\phi ^{m_1},p)\leq 0,\forall p\in Fix(\phi ). 
$$

This happens because all points in $Fix(\phi )$ with non-positive indexes
are saddle-like (maybe with more than four sectors), with $\phi $-invariant
separatrices and points with positive indexes are rotating saddles. So for
some $m_1>0$ sufficiently large, $\phi ^{m_1}$ fixes the separatrices of all
points in $Fix(\phi ),$ and thus they all have non-positive indexes with
respect to $\phi ^{m_1}.$ In particular $ind(\phi ^{m_1},\pi (\widetilde
p_1))<0. $

Now let us look at $\phi ^{m_1}.$ Again, as a consequence of Handel's
theorem there is a fixed point $\widetilde p_2$ of $\widetilde \phi ^{m_1}$
with $ind(\widetilde \phi ^{m_1},\widetilde p_2)=ind(\phi ^{m_1},\pi
(\widetilde p_2))>0.$ In the same way as above for some sufficiently large $%
m_2>0,$ the local dynamics at points in $Fix(\phi ^{m_1})$ implies that 
$$
ind(\phi ^{m_1m_2},p)\leq 0,\forall p\in Fix(\phi ^{m_1})), 
$$
and in particular $ind(\phi ^{m_1m_2},\pi (\widetilde p_2))<0$.

If we continue this process we get a sequence of pairwise different points $%
\widetilde p_1,\widetilde p_2,\widetilde p_3,\ldots $. In $S,$ the points $%
\pi (\widetilde p_1),\pi (\widetilde p_2),\pi (\widetilde p_3),\ldots $ are
also pairwise different.

So at some $j,$ the cardinality of $\{\widetilde p_1,\widetilde p_2,\ldots
,\widetilde p_j\}$ is larger than the number of singularities of the
foliations $\mathscr{F}^u$, $\mathscr{F}^s.$ This implies that for some $%
\widetilde \phi $-periodic point $\widetilde p$, $\pi (\widetilde p)$ does
not coincide with a singularity of the foliations $\mathscr{F}^u,$ $%
\mathscr{F}^s.$ Hence $\widetilde p$ is a hyperbolic periodic saddle point
for $\widetilde \phi.$ \qed

\subsection{A first result towards the proof of theorem 2 in the rel p.A.
case}

The stable and unstable foliations for $\phi $ lift to stable and unstable
foliations for $\widetilde \phi .$ If $F_p^s$ is the stable leaf of $%
\mathscr{F}^s$ that contains a point $p\in S,$ we will denote by $\widetilde{%
F}_{ \widetilde{p}}^s$ the lift of $F_p^s$ that contains a point $\widetilde
p\in \pi ^{-1}(p)$. The same for unstable leaves of $\mathscr{F}^u$.

Now we will state some definitions and properties of pseudo-Anosov maps
relative to finite invariant sets, which will be useful in the proof of the
next lemma.

Let $p\in S$ be a fixed point of $\phi .$ As we already said, the dynamics
of a sufficiently large iterate of $\phi $ in a neighborhood of $p$ can be
obtained gluing finitely many invariant hyperbolic sectors together. In each
sector the dynamics is locally like the dynamics in the first quadrant of
the map $(x,y)\mapsto (\lambda _1x,\lambda _2y)$, for some real numbers $%
0<\lambda _2<1<\lambda _1$.

We define the stable set of $p$ as the set $W^s(p)$ of points $z$ in $S$
such that $\phi ^n(z)\to p$ when $n\to \infty $ and the unstable set of $p$
as the set $W^u(p)$ of points $z$ in $S$ such that $\phi ^{-n}(z)\to p$ when 
$n\to \infty .$ If $p$ is a regular point of the foliations $\mathscr{F}^s,%
\mathscr{F}^u$ then $W^u(p)$ is the union of 2 branches, the same for $%
W^s(p).$ This is the situation we called the point a hyperbolic saddle point
in the previous proposition. In case $p$ is a singular point of the
foliations, $p$ is a $k$-prong singularity (for $k=1$ or some $k\geq 3),$
which implies that $W^u(p)$ is the union of $k$ branches, the same for $%
W^s(p).$ In this singular case, each branch is actually a leaf of the proper
foliation, which emanates from the singularity, while in the regular case
each leaf gives 2 branches. In both the regular and the singular cases, the
branches are either invariant, or rotated around $p$ under iterates of $\phi 
$ (thus, $\phi ^n$-invariant for some $n>0).$

In case $p^{\prime }\in S$ is a $\phi $-periodic point, if $n_{p^{\prime }}$
is the least period of $p^{\prime }$, then it is a fixed point of $\phi
^{n_{p^{\prime }}}$, so we define the stable and unstable sets of $p^{\prime
}$ accordingly, using $\phi ^{n_{p^{\prime }}}$ instead of $\phi $.

\begin{lemma}
\label{lemg1g2} Let $\widetilde{\phi }$ be the natural lift of $\phi $. Then
there exists $\widetilde{p}\in \mathbb{D}$ a $\widetilde{\phi }$-hyperbolic
periodic saddle point and deck transformations $g_1,$ $g_2$ such that $%
g_1\circ g_2\neq g_2\circ g_1$ and 
$$
\widetilde{F}_{\widetilde{p}}^{u+}\pitchfork \widetilde{F}_{g_i(\widetilde{p}%
)}^{s+},i\in \{1,2\}, 
$$
where $W^u(\widetilde{p})=\widetilde{F}_{\widetilde{p}}^{u+}\cup \widetilde{F%
}_{\widetilde{p}}^{u-},$ $W^s(\widetilde{p})=\widetilde{F}_{\widetilde{p}%
}^{s+}\cup \widetilde{F}_{\widetilde{p}}^{s-}$ and $\widetilde{F}_{
\widetilde{p}}^{u+},\widetilde{F}_{\widetilde{p}}^{u-},\widetilde{F}_{
\widetilde{p}}^{s+},\widetilde{F}_{\widetilde{p}}^{s-}$ are the four
branches at $\widetilde{p}.$
\end{lemma}

{\it Proof: Let }$\widetilde p\in \mathbb{D}$ be the $\widetilde \phi $%
-periodic point given in proposition \ref{periodic}. So, $p=\pi (\widetilde
p)$ is a hyperbolic $\phi $-periodic saddle point. Without loss of
generality, considering an iterate of $\phi $ if necessary, we will assume
that each point in $K=\{p\}\cup P$ is fixed and moreover, each stable or
unstable branch at a point in $K$ is also invariant under $\phi .$

The map $\phi $ is pseudo-Anosov relative to the finite $\phi $-invariant
set $P$. In particular, any stable leaf $F^s\in \mathscr{F}^s$ intersects
all unstable leaves $F^u\in \mathscr{F}^u$ $C^1$-transversely and
vice-versa. Let $F_p^u$ be the unstable leaf at the point $p$ (as $p$ is
regular, $F_p^u=W^u(p))$ and $F_{*p^{\prime }}^s$ be a stable leaf at some
point $p^{\prime }\in P=\{p_1^{+},p_1^{-},\ldots ,p_k^{+},p_k^{-}\}$. The
point $p^{\prime }$ may be singular or regular. From what we said above, $%
F_p^u\pitchfork F_{*p^{\prime }}^s.$ So, there exists an unstable branch at $%
p,$ denoted $F_p^{u+},$ and an unstable branch at $p^{\prime },$ denoted $%
F_{*p^{\prime }}^{u^{\prime }},$ such that $F_p^{u+}$ accumulates on $%
F_{*p^{\prime }}^{u^{\prime }}$ and $F_{*p^{\prime }}^{u^{\prime }}%
\pitchfork W^s(p).$ Let $F_p^{s+}$ be a stable branch at $p$ such that $%
F_{*p^{\prime }}^{u^{\prime }}\pitchfork F_p^{s+}.$ Lifting everything to
the universal cover, fixed some $\widetilde{p}\in \pi ^{-1}(p),$ there exist
deck transformations $g^{\prime }\neq Id$ and $h$ such that 
\begin{equation}
\label{firstinterac}\widetilde{F}_{\widetilde{p}}^{u+}\pitchfork \widetilde{F%
}_{(g^{\prime })^nh(\widetilde{p})}^{s+},
\end{equation}
for all sufficiently large $n>0.$ This follows from the fact that fixed some 
$\widetilde{p}\in \pi ^{-1}(p),$ there exist a $\widetilde{p}^{\prime }\in
\pi ^{-1}(p^{\prime })$ and deck transformations $g^{\prime }$ and $h$ such
that $\widetilde{\phi }(\widetilde{p}^{\prime })=g^{\prime }(\widetilde{p}%
^{\prime }),$ $\widetilde{F}_{\widetilde{p}}^{u+}\pitchfork \widetilde{F}_{*
\widetilde{p}^{\prime }}^s$ and $\widetilde{F}_{*\widetilde{p}^{\prime
}}^{u^{\prime }}\pitchfork \widetilde{F}_{h(\widetilde{p})}^{s+}.$

Let $g_1=(g^{\prime }){}^nh$ for some $n>0$ such that (\ref{firstinterac})
holds. Now consider $\widetilde \theta $ a path in $\mathbb{D}$ constructed
as follows: $\widetilde \theta =\widetilde \theta ^{\prime }*\widetilde
\theta ^{\prime \prime }$, where $\widetilde \theta ^{\prime }$ is a compact
subarc of $\widetilde{F}_{\widetilde{p}}^{u+}$ starting at $\widetilde p$
and ending at a point in $\widetilde{F}_{\widetilde{p}}^{u+}\cap \widetilde{F%
}_{g_1( \widetilde{p})}^{s+},$ and $\widetilde \theta ^{\prime \prime }$ is
a compact subarc of $\widetilde{F}_{g_1(\widetilde{p})}^{s+}$ starting at
the endpoint of $\widetilde \theta ^{\prime }$ and ending at $g_1(\widetilde
p).$

Let $\omega _1$ be the fixed point in $\partial \mathbb{D}$ of $g_1$ such
that $lim_{n\to \infty }g_1^n(\widetilde q)=\omega _1$, for all $\widetilde
q\in \overline{\mathbb{D}}$ and let $\alpha _1$ be the other fixed point.

Define 
$$
\Theta =\bigcup_{i\in \mathbb{Z}}g_1^i(\widetilde \theta ). 
$$

By construction, $\Theta $ is a path connected subset of $\mathbb{D},$
joining $\alpha _1$ to $\omega _1$. Since $S\setminus \mathscr{C}$ is a
union of open topological disks, there exists a geodesic in $\mathscr{C}$
and $m\in Deck(\pi )$ such that the projection of the axis of $m$ in $S$ is
the previous geodesic and the fixed points of $m$ in $\partial \mathbb{D}$
separate the endpoints $\omega _1$ and $\alpha _1$ of $\Theta $. This
follows from proposition 7 and the fact that the endpoints of the extended
lifts of the geodesics in $\pi ^{-1}(\mathscr{C})$ are dense in $\partial 
\mathbb{D}.$

Now consider the fixed points $\omega _m$ and $\alpha _m$ of $m$ in $%
\partial \mathbb{D}$ such that $\lim _{n\to \infty }m^n(\widetilde q)=\omega
_m$ and $\lim _{n\to -\infty }m^n(\widetilde q)=\alpha _m$, for all $%
\widetilde q\in \mathbb{D}$. Let $n_0>0$ be a sufficiently large integer,
such that $m^{n_0}(\omega _1)$ and $m^{n_0}(\alpha _1)$ are close to $\omega
_m$ and $\Theta \cap m^{n_0}(\Theta )=\emptyset .$ This is possible because $%
\Theta $ accumulates on $\omega _m$ under positive iterates of $m.$

The following holds: 
$$
\Theta =\bigcup_{i\in \mathbb{Z}}g_1^i(\widetilde \theta )\Rightarrow
m^{n_0}(\Theta )=\bigcup_{i\in \mathbb{Z}}m^{n_0}g_1^i(\widetilde \theta ). 
$$

As $\widetilde \phi $ commutes with all deck transformations, $\widetilde
\theta ^{\prime }\subset \widetilde{F}_{\widetilde p}^{u+}$ and $\widetilde{%
\phi }( \widetilde{F}_{\widetilde p}^{u+})=\widetilde{F}_{\widetilde
p}^{u+}, $ we get that for all $n>0$ and $t\in Deck(\pi ),$ $t(\widetilde
\theta ^{\prime })\subset \widetilde \phi ^n(t(\widetilde \theta ^{\prime
})).$ Similarly, since $\widetilde \theta ^{\prime \prime }\subset 
\widetilde{F}_{g_1(\widetilde p)}^{s+},$ $\widetilde \phi ^n(t(\widetilde
\theta ^{\prime \prime }))\subset t(\widetilde \theta ^{\prime \prime })$,
for all $n>0$ and $t\in Deck(\pi ).$

\begin{figure}[!h]
	\centering
	\includegraphics[scale=0.21]{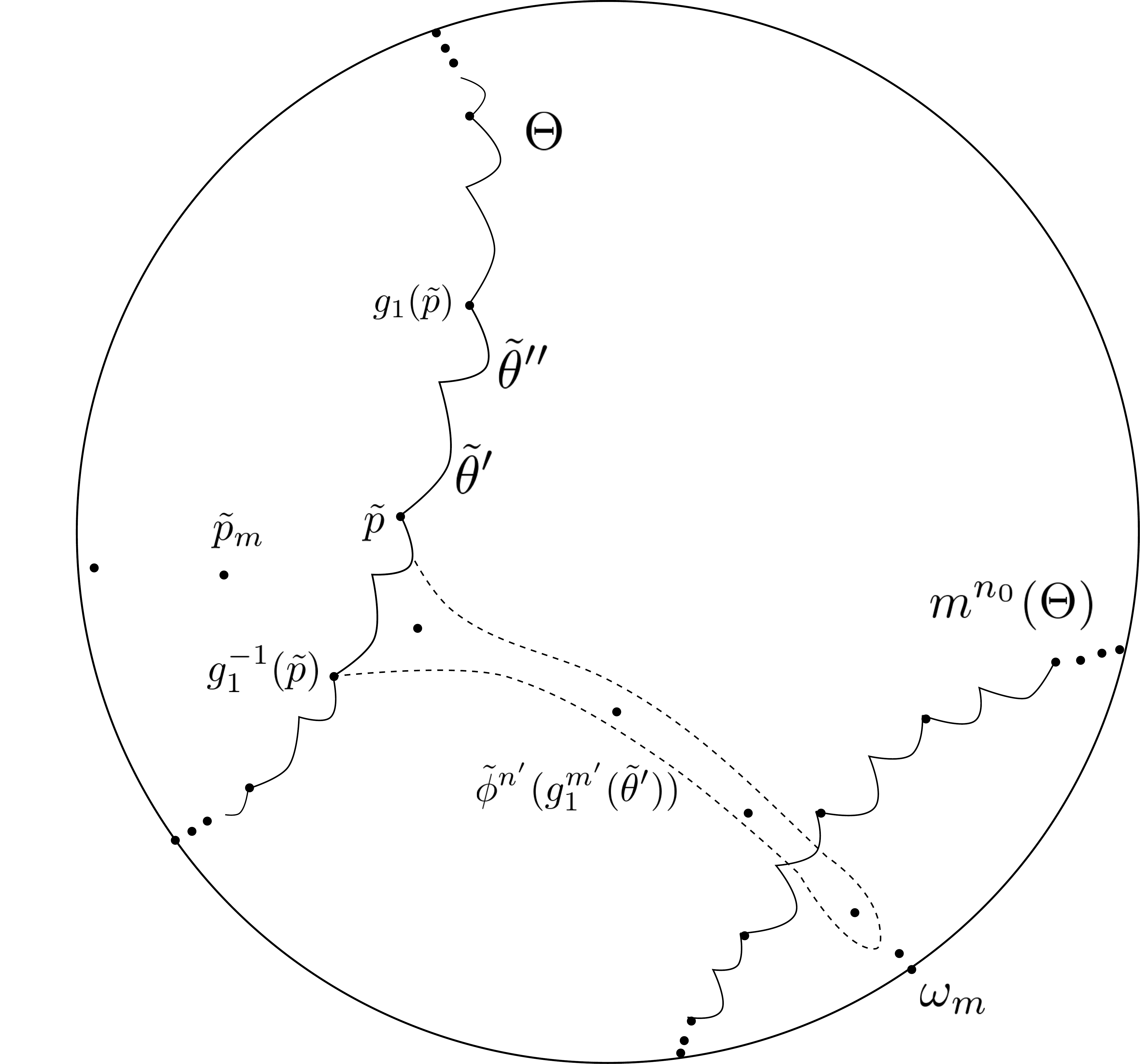}
	\caption{How to obtain $g_2$.}
	\label{g1g2}
\end{figure}

The hypotheses on $\mathscr{C}$ implies that there is a point $\widetilde
p_m\in \pi ^{-1}(P)$, such that $\widetilde{\phi }(\widetilde
p_m)=m(\widetilde p_m)$ and $\widetilde p_m$ is in the connected component
of $\mathbb{D}\setminus \Theta $ which contains $\alpha _m$ in its boundary.
As $m^{n_0}(\Theta )$ is in the other connected component of $\mathbb{D}%
\setminus \Theta ,$ $\lim _{n\to \infty }\widetilde \phi ^n(\widetilde
p_m)=\omega _m$ and $\widetilde \phi |_{\partial \mathbb{D}}=Id,$ we get
that for a sufficiently large $n^{\prime }>0$ there must exists two integers 
$i^{\prime },$ $i^{\prime \prime }$ such that 
$$
\widetilde \phi ^{n^{\prime }}(g_1^{i^{\prime }}(\widetilde \theta ^{\prime
}))\pitchfork m^{n_0}g_1^{i^{\prime \prime }}(\widetilde \theta ^{\prime
\prime }). 
$$

In particular, $\widetilde F_{g_1^{i^{\prime }}(\widetilde p)}^{u+}%
\pitchfork \widetilde F_{m^{n_0}g_1^{i^{\prime \prime }}(\widetilde p)}^{s+}$%
, and so, $\widetilde F_{\widetilde p}^{u+}\pitchfork \widetilde
F_{(g_1^{-i^{\prime }}m^{n_0}g_1^{i^{\prime \prime }})(\widetilde p)}^{s+}$
(see figure \ref{g1g2}).

Finally, let $g_2=g_1^{-i^{\prime }}m^{n_0}g_1^{i^{\prime \prime }}$. We
will show that $g_1$ and $g_2$ do not commute. If $g_1\circ g_2=g_2\circ
g_1, $ then there exists $l\in Deck(\pi )$ and integers $k_1 $, $k_2$ such
that $g_1=l^{k_1}$ and $g_2=l^{k_2}$. Thus 
$$
g_1^{-i^{\prime }}m^{n_0}g_1^{i^{\prime \prime }}=l^{k_2}\Rightarrow
l^{-i^{\prime }k_1}m^{n_0}l^{i^{\prime \prime }k_1}=l^{k_2}\Rightarrow
m^{n_0}=l^{k_2+k_1(i^{\prime }-i^{\prime \prime })}. 
$$

Since $m^{n_0}$ and $g_1$ are iterates of the same deck transformation, the
geodesics associated to the axes of $m$ and $g_1$ are equal. But this is a
contradiction with our choice of $m$. So, $g_1 $ and $g_2$ do not commute. $%
\qed$

\subsection{Proof of theorem 2 in a special case}

In this subsection we prove theorem 2 in case of relative pseudo-Anosov maps.

\begin{description}
\item[Remark]  \label{impremark}{\bf \ \ref{impremark}:\ }As $\phi
:S\rightarrow S$ is pseudo-Anosov relative to a finite invariant set, if for
some leaves $F^u$ of $\mathscr{F}^u$ and $F^s$ of $\mathscr{F}^s,$ there are
connected components $\widetilde{F}^u$ of $\pi ^{-1}(F^u)$ and $\widetilde{F}%
^s$ of $\pi ^{-1}(F^s)$ which intersect (not at a lift of a singularity of
the foliations), then they intersect in a $C^1$-transverse way. In the proof
of the next lemma we will not make use of this fact because when proving
theorem 2, at some point we say that the proof continues as the proof of the
next lemma. So, in the proof of lemma \ref{lema2}, although intersections
between stable and unstable leaves, either in $S$ or in $\mathbb{D},$ are
always $C^1$-transverse, we will not use this fact.

Moreover, as we said in the introduction, the main feature of topologically
transverse intersections is the fact that a $C^0$-version of the so called $%
\lambda $-lemma (see \cite{palis1}) holds: If $M$ is a surface, $%
f:M\rightarrow M$ is a $C^1$ diffemorphism, $p,q\in M$ are $f$-periodic
saddle points and $W^u(p)$ has a topologically transverse intersection with $%
W^s(q),$ then $W^u(p)$ $C^0$-accumulates on $W^u(q),$ in particular $
\overline{W^u(p)}\supset \overline{W^u(q)}.$ So if $p_1,p_2,p_3\in M$ are
hyperbolic $f$-periodic saddle points, $W^u(p_1)$ has a topologically
transverse intersection with $W^s(p_2)$ and $W^u(p_2)$ has a topologically
transverse intersection with $W^s(p_3),$ then $W^u(p_1)$ has a topologically
transverse intersection with $W^s(p_3)$. 
\end{description}



\begin{lemma}
\label{lema2} (Theorem 2 in case of relative p.A. maps) Let $\widetilde{\phi 
}$ be the natural lift of the map $\phi $. Then there exists a contractible
hyperbolic $\phi $-periodic point $p\in S,$ such that for any $\widetilde{p}%
\in \pi ^{-1}(p)$ and any given $g\in Deck(\pi )$, 
$$
W^u(\widetilde{p})\pitchfork W^s(g(\widetilde{p})). 
$$
\end{lemma}

{\it Proof:\ }Let $\widetilde p$ be the hyperbolic $\widetilde \phi $%
-periodic point from lemma \ref{lemg1g2} and $p=\pi (\widetilde{p}).$ That
lemma implies the existence of $g_1$ and $g_2$ in $Deck(\pi )$ and also the
existence of an unstable branch $\lambda _u$ of $W^u(p)$ and a stable branch 
$\beta _s$ of $W^s(p),$ such that if $\widetilde{\lambda }_u$ is the
connected component of $\pi ^{-1}(\lambda _u)$ contained in $W^u(\widetilde{p%
})$ and $\widetilde{\beta }_s$ is the connected component of $\pi
^{-1}(\beta _s)$ contained in $W^s(\widetilde{p}),$ then%
$$
\widetilde{\lambda }_u\pitchfork g_i(\widetilde{\beta }_s),i\in \{1,2\}. 
$$

Without loss of generality, as we did in lemma \ref{lemg1g2}, considering an
iterate of $\widetilde \phi $ if necessary, we will assume that $\widetilde
\phi (\widetilde p)=\widetilde p,$ $\widetilde \phi (\widetilde{\lambda }%
_u)= \widetilde{\lambda }_u$ and $\widetilde \phi (\widetilde{\beta }_s)= 
\widetilde{\beta }_s.$

Since $\widetilde \phi $ is the natural lift of $\phi ,$ every point of the
form $h(\widetilde p)$ with $h\in Deck(\pi ),$ is fixed under $\widetilde
\phi .$ Moreover, if we consider the stable set of the point $h(\widetilde
p) $ with respect to $\widetilde \phi ,$ the following equality holds: 
$$
W^s(h(\widetilde p))=h(W^s(\widetilde p)), 
$$
and the same is true for the unstable set of $\widetilde p$.

Consider the point $p\in S.$ Choose $\epsilon >0$ small enough, so that $%
B_\epsilon (\widetilde p)\cap \pi ^{-1}(p)=\widetilde p$, where $B_\epsilon
(\widetilde p)=\{\widetilde q\in \mathbb{D}|d_{\mathbb{D}}(\widetilde
p,\widetilde q)<\epsilon )\}$. Observe that, since every point on the fiber
of $p$ is of the form $h(\widetilde p)$ for some $h\in Deck(\pi ),$ and $h$
is an isometry, $B_\epsilon (h(\widetilde p))\cap \pi ^{-1}(p)=h(\widetilde
p)$.

\begin{figure}[!h]
	\centering
	\includegraphics[scale=0.18]{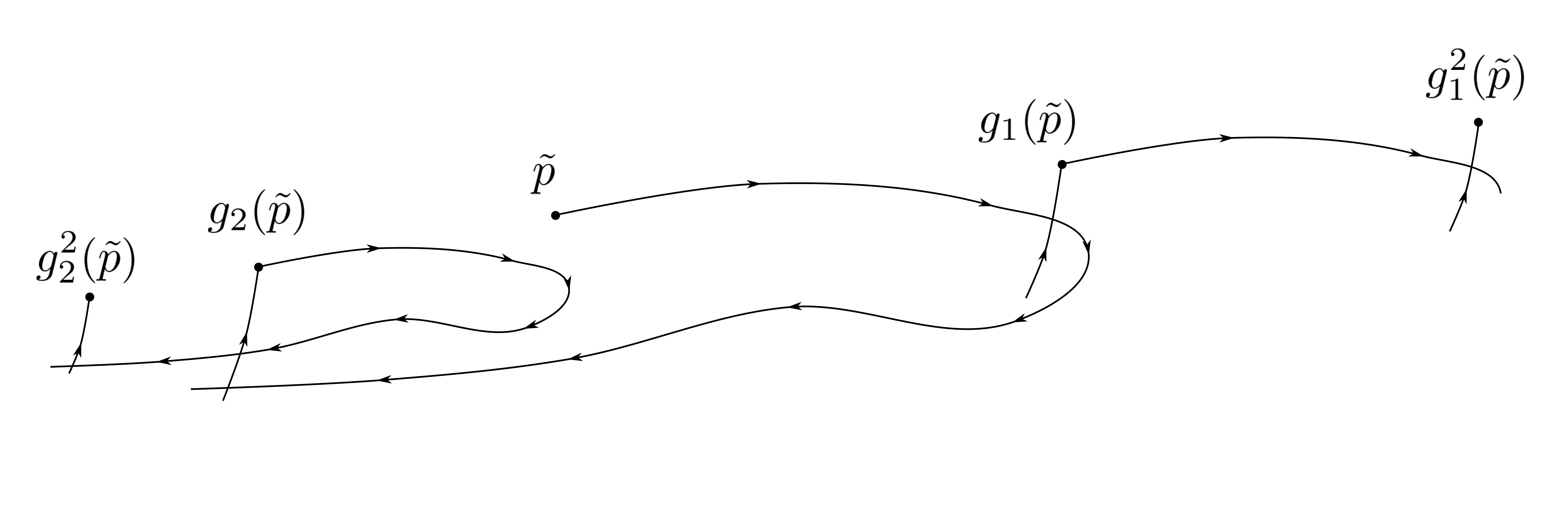}
	\caption{The construction of the path connected set $\theta$.}
	\label{curve}
\end{figure}

From the fact that $\widetilde{\lambda }_u\pitchfork g_1(\widetilde{\beta }%
_s),$ we can construct a path $\eta _1$ in $\mathbb{D}$ joining $\widetilde
p $ to $g_1(\widetilde p)$ exactly as in the previous lemma: $\eta _1$
starts at $\widetilde p,$ consists of a compact connected piece of $
\widetilde{\lambda }_u$ until it reaches $g_1(\widetilde{\beta }_s)$ 
and then it continues as a compact connected piece of $g_1(\widetilde{\beta }%
_s)$ until it reaches $g_1(\widetilde p)$. It is clear that we can choose
the piece that belongs to $g_1(\widetilde{\beta }_s)$ totally contained in $%
B_\epsilon (g_1(\widetilde p))$. Analogously, we construct a path $\eta _2$
in $\mathbb{D}$ joining $\widetilde p$ to $g_2(\widetilde p)$. Let $\theta
\subset \mathbb{D}$ be the path connected set obtained as follows (see
figure \ref{curve}):

\begin{equation}
\label{construct} \theta = \bigg(\bigcup_{i \geq 0} g_1^i(\eta_1)\bigg) \cup 
\bigg(\bigcup_{j \geq 0}g_2^j(\eta_2)\bigg). 
\end{equation}

One can think of $\theta $ geometrically as the concatenation of two curves
in $\overline{\mathbb{D}}$, one joining $\widetilde p$ to $\omega _{g_1}$
and another joining $\widetilde p$ to $\omega $$_{g_2}$, where $\omega $$%
_{g_1}$ and $\omega $$_{g_2}$ are the attractive fixed points at infinity of 
$g_1$ and $g_2$ respectively. The fact that $g_1$ and $g_2$ do not commute
implies that the fixed points at infinity of these deck transformations are
all different, so in particular $\omega _{g_1}\neq \omega _{g_2}$.

We want to show that for every $g\in Deck(\pi )$ 
\begin{equation}
\label{whatwewant}\widetilde{\lambda }_u\pitchfork g(\widetilde{\beta }%
_s)\Rightarrow W^u(\widetilde p)\pitchfork W^s(g(\widetilde p)). 
\end{equation}

Fix $g\in Deck(\pi )$ with $g\neq Id$. The case $g=Id$ will be considered in
the end.

As $\widetilde{\lambda }_u\pitchfork g_1(\widetilde{\beta }_s),$ we get that 
$$
\begin{array}{c}
g_1^{-1}( 
\widetilde{\lambda }_u)\pitchfork \widetilde{\beta }_s\text{ and so }%
gg_1^{-1}(\widetilde{\lambda }_u)\pitchfork g(\widetilde{\beta }_s) \\ \text{%
that can be rewritten as } \\ gg_1^{-1}g^{-1}(g(\widetilde{\lambda }_u))%
\pitchfork g(\widetilde{\beta }_s). 
\end{array}
$$

Notice that an analogous statement holds for $g_2.$ Using this, let us
construct a path connected set $\theta ^{\prime }$ containing $g(\widetilde
p)$ in a similar way as $\theta .$ An important simple observation here is
the fact that for a fixed point of $\widetilde \phi ,$ its stable set with
respect to $\widetilde \phi $ coincides with its unstable set with respect
to $\widetilde \phi ^{-1}$. This duality allows us to construct the set $%
\theta ^{\prime }$ in the same way as $\theta ,$ but using the point $%
g(\widetilde p),$ the deck transformations $gg_1^{-1}g^{-1}$, $%
gg_2^{-1}g^{-1}$ and the map $\widetilde \phi ^{-1}$. Hence we construct a
path $\eta _1^{\prime }$ joining $g(\widetilde p)$ to $gg_1^{-1}g(g(%
\widetilde p))$ such that $\eta _1^{\prime }$ starts at $g(\widetilde p)$,
consists of a compact connected piece of $g(\widetilde{\beta }_s)$ until it
reaches $gg_1^{-1}g^{-1}(g(\widetilde{\lambda }_u))$ 
and then it continues as compact connected piece of $gg_1^{-1}g^{-1}(g( 
\widetilde{\lambda }_u))\cap B_\epsilon (gg_1^{-1}g^{-1}(g(\widetilde p)))$
until it reaches $gg_1^{-1}g^{-1}(g(\widetilde p))$. Constructing $\eta
_2^{\prime }$ analogously, we define 
$$
\theta ^{\prime }=\bigg(\bigcup_{i\geq 0}gg_1^{-i}g^{-1}(\eta _1^{\prime })%
\bigg) \cup \bigg(\bigcup_{j\geq 0}gg_2^{-j}g^{-1}(\eta _2^{\prime })\bigg). 
$$

Similar to what was said about $\theta ,$ one can think of $\theta ^{\prime
} $ as the concatenation of two curves in $\overline{\mathbb{D}}$, one
joining $g(\widetilde p)$ to $g(\alpha _{g_1})$ and another joining $%
g(\widetilde p)$ to $g(\alpha _{g_2}),$ where $\alpha _{g_1}$ and $\alpha
_{g_2}$ are the repulsive fixed points at infinity of $g_1$ and $g_2$
respectively.

The sets $\theta $ and $\theta ^{\prime }$ have analogous properties:

\begin{description}
\item[Properties of]  {\bf $\theta $ and $\theta ^{\prime }\label
{proptetatetal}$} {\bf \ref{proptetatetal}:}

\begin{itemize}
\item  for $i\in \{1,2\},$ all points of the form $g_i^m(\widetilde{p})\in
\theta $ and $\forall $$m>0,$ $\widetilde{\lambda }_u\pitchfork g_i^m(
\widetilde{\beta }_s).$ Remember that $\widetilde{\lambda }_u$ is a branch
of $W^u(\widetilde{p})$ and $g_i^m(\widetilde{\beta }_s)$ is a branch of $%
W^s(g_i^m(\widetilde{p}));$

\item  for $i\in \{1,2\},$ all points of the form $gg_i^{-m}g^{-1}(g(
\widetilde{p}))\in \theta ^{\prime }$ and

\noindent $\forall m>0,$ $gg_i^{-m}g^{-1}(g(\widetilde{\lambda }_u))%
\pitchfork g(\widetilde{\beta }_s).$ And in this case, $gg_i^{-m}g^{-1}(g(
\widetilde{\lambda }_u))$ is a branch of $W^u(gg_i^{-m}g^{-1}(g(\widetilde{p}%
)))$ and $g(\widetilde{\beta }_s)$ is a branch of $W^s(g(\widetilde{p}));$
\end{itemize}
\end{description}

If $\theta $ and $\theta ^{\prime }$ have a topologically transverse
intersection, then the lemma is proved.

Indeed, if there is such an intersection, then at least one of the following
four possibilities holds:

\begin{itemize}
\item  there exists $j^{\prime }\in \{1,2\}$ and $m^{\prime }>0$ such that $%
g_{j^{\prime }}^{m^{\prime }}(\widetilde{p})=g(\widetilde{p}).$

\item  there exists $j^{\prime }\in \{1,2\}$ and $m^{\prime }>0$ such that $%
gg_{j^{\prime }}^{-m^{\prime }}g^{-1}(g(\widetilde{p}))=\widetilde{p}.$

\item  there exists $j^{\prime },j^{\prime \prime }\in \{1,2\}$ and $%
m^{\prime },m^{\prime \prime }>0$ such that $g_{j^{\prime }}^{m^{\prime }}(
\widetilde{p})=gg_{j^{\prime \prime }}^{-m^{\prime \prime }}g^{-1}(g(
\widetilde{p})).$
\end{itemize}

In the three possibilities above, using properties \ref{proptetatetal}, we
get that (\ref{whatwewant}) holds. The last possibility is the following:

\begin{itemize}
\item  there exists $j^{\prime },j^{\prime \prime }\in \{1,2\}$ and $%
m^{\prime },m^{\prime \prime }\geq 0$ such that some compact piece of $%
\theta \cap W^u(g_{j^{\prime }}^{m^{\prime }}(\widetilde{p}))$ has a
topologically transverse intersection with some compact piece of $\theta
^{\prime }\cap W^s(gg_{j^{\prime \prime }}^{-m^{\prime \prime }}g^{-1}(g(
\widetilde{p}))).$ This happens because for $i\in \{1,2\}$ and all $n\geq 0,$
$$
\begin{array}{c}
\theta \cap W^s(g_i^n(
\widetilde{p}))\subset B_\epsilon (g_i^n(\widetilde{p})) \\ \text{and} \\ 
\theta ^{\prime }\cap W^u(gg_i^{-n}g^{-1}(g(\widetilde{p})))\subset
B_\epsilon (gg_i^{-n}g^{-1}(g(\widetilde{p})))
\end{array}
$$
and all of these balls are disjoint. So, by the $C^{0\text{ }}\lambda $%
-lemma mentioned in remark \ref{impremark}, we are done.
\end{itemize}

Hence let us suppose that $\theta $ and $\theta ^{\prime }$ do not have
topologically transverse intersections. Our goal is to show that, in this
case, using the fully essential system of curves $\mathscr{C}$ and the
periodic points associated with the geodesics, we can force a topologically
transverse intersection between $\theta ^{\prime }$ and a path connected set 
$\theta _0\in \mathbb{D}$ that has the same properties, and is obtained from 
$\theta .$

\begin{figure}[!h]
	\centering
	\includegraphics[scale=0.22]{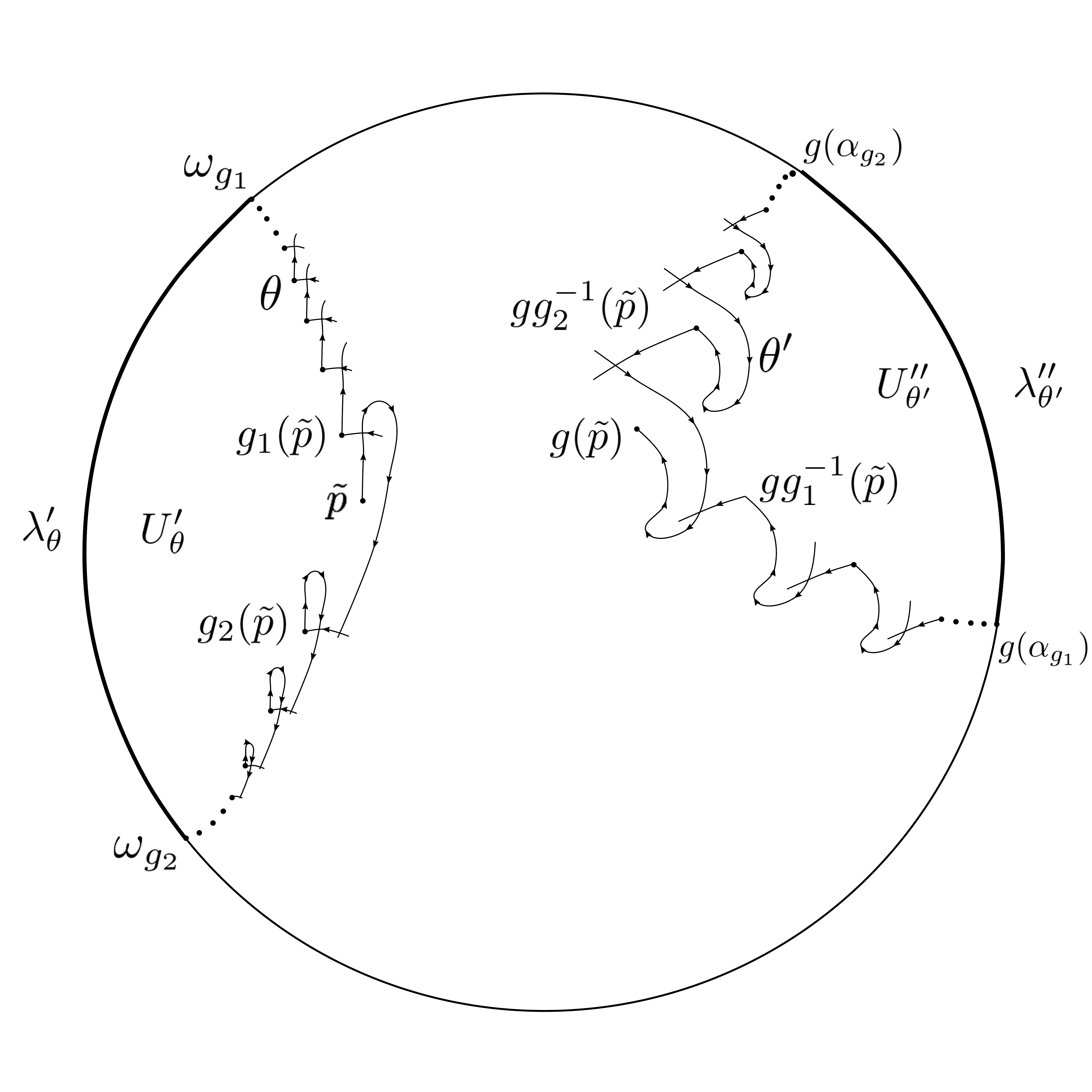}
	\caption{The case when $\theta$ and $\theta'$ do not have topologically transverse intersections.}
	\label{tetatetal}
\end{figure}

So suppose that $\theta $ and $\theta ^{\prime }$ do not have topologically
transverse intersections. The set $\mathbb{D}\setminus \theta $ has two
unbounded connected components $U_\theta ^{\prime }$ and $U_\theta ^{\prime
\prime },$ the closure of one of them containing $\theta ^{\prime }$. We
will assume that $\theta ^{\prime }\subset closure(U_\theta ^{\prime \prime
})$. The boundary at infinity of $U_\theta ^{\prime }$ is equal to a segment
of $\partial \mathbb{D}$ delimited by $\omega _{g_1}$ and $\omega _{g_2}$
that will be denoted $\lambda _\theta ^{\prime }.$ Similarly, the boundary
at infinity of $U_\theta ^{\prime \prime }$ is equal to a segment of $%
\partial \mathbb{D}$ delimited by $\omega _{g_1}$ and $\omega _{g_2}$ that
will be denoted $\lambda _\theta ^{\prime \prime }.$ In the same way, $%
\mathbb{D}\setminus \theta ^{\prime }$ has two unbounded connected
components $U_{\theta ^{\prime }}^{\prime }$ and $U_{\theta ^{\prime
}}^{\prime \prime }$. We will assume that $\theta \subset closure(U_{\theta
^{\prime }}^{\prime })$ and call $\lambda _{\theta ^{\prime }}^{\prime }$, $%
\lambda _{\theta ^{\prime }}^{\prime \prime }$ the segments of $\partial 
\mathbb{D}$ delimited by $g(\alpha _{g_1})$ and $g(\alpha _{g_2})$ that are
equal to the boundary at infinity of $U_{\theta ^{\prime }}^{\prime }$ and $%
U_{\theta ^{\prime }}^{\prime \prime }$ respectively. Then $\lambda _\theta
^{\prime }\subseteq \lambda _{\theta ^{\prime }}^{\prime }$ and $\lambda
_{\theta ^{\prime }}^{\prime \prime }\subseteq \lambda _\theta ^{\prime
\prime }$. See figure \ref{tetatetal}.

By proposition \ref{connectedC}, $\pi ^{-1}(\mathscr{C})$ is a closed
connected subset of $\mathbb{D}$. Moreover $\pi ^{-1}(\mathscr{C})\cap
U_\theta ^{\prime }\neq \emptyset $ and $\pi ^{-1}(\mathscr{C})\cap
U_{\theta ^{\prime }}^{\prime \prime }\neq \emptyset .$ This happens because
the set $A=\{\omega _g\in \partial \mathbb{D}|g\in Deck(\pi )\}$ is dense in 
$\partial \mathbb{D},$ so given an extended lift $\widetilde \gamma $ of one
geodesic $\gamma \in \mathscr{C},$ there exists $h\in Deck(\pi )$ such that $%
\omega _h$ is in the interior of the segment $\lambda _\theta ^{\prime }.$
Then for a sufficiently large $n>0,$ $h^n(\widetilde \gamma )$ is contained
in an arbitrarily small Euclidean neighborhood of $\omega _h,$ so small that
it is contained in $U_\theta ^{\prime }$. Since $\pi ^{-1}(\mathscr{C})$ is
invariant under $h,$ this implies that $\pi ^{-1}(\mathscr{C})\cap U_\theta
^{\prime }\neq \emptyset .$ Analogously, $\pi ^{-1}(\mathscr{C})\cap
U_{\theta ^{\prime }}^{\prime \prime }\neq \emptyset .$

Let $C(U_\theta ^{\prime },U_{\theta ^{\prime }}^{\prime \prime })$ be the
set of paths in $\pi ^{-1}(\mathscr{C})$ joining a point $\widetilde
r^{\prime }\in \pi ^{-1}(\mathscr{C})\cap U_\theta ^{\prime }$ to a point $%
\widetilde r^{\prime \prime }\in \pi ^{-1}(\mathscr{C})\cap U_{\theta
^{\prime }}^{\prime \prime }$ and formed by a finite number of subarcs of
extended lifts of geodesics in $\mathscr{C}.$ Proposition \ref{finitepath}
implies that $C(U_\theta ^{\prime },U_{\theta ^{\prime }}^{\prime \prime
})\neq \emptyset . $

For every $\beta $ in $C(U_\theta ^{\prime },U_{\theta ^{\prime }}^{\prime
\prime })$ we can write 
$$
\beta =\beta _1*\beta _2*\ldots *\beta _l, 
$$
where each $\beta _i,$ $i\in \{1,2,\ldots ,l\},$ is a subarc of an extended
lift of one geodesic in $\mathscr{C}.$ We will consider $C^{\prime
}(U_\theta ^{\prime },U_{\theta ^{\prime }}^{\prime \prime })\subseteq
C(U_\theta ^{\prime },U_{\theta ^{\prime }}^{\prime \prime })$ the set of
all paths $\beta \in C(U_\theta ^{\prime },U_{\theta ^{\prime }}^{\prime
\prime })$ satisfying the following property: if $\widetilde \gamma _1$ is
the extended lift of the geodesic in $\mathscr{C}$ such that $\beta
_1\subset \widetilde \gamma _1$ and $\widetilde \gamma _l$ is the extended
lift of the geodesic in $\mathscr{C}$ such that $\beta _l\subset \widetilde
\gamma _l,$ then at least one of the points at infinity of $\widetilde
\gamma _1$ is in the interior of $\lambda _\theta ^{\prime }\subset \partial 
\mathbb{D}$ and at least one of the points at infinity of $\widetilde \gamma
_l$ is in the interior of $\lambda _{\theta ^{\prime }}^{\prime \prime
}\subset \partial \mathbb{D}$. Clearly $C^{\prime }(U_\theta ^{\prime
},U_{\theta ^{\prime }}^{\prime \prime })\neq \emptyset $.

Observe that when we consider some $\beta \in C^{\prime }(U_\theta ^{\prime
},U_{\theta ^{\prime }}^{\prime \prime })$, since $\beta $ connects a point $%
\widetilde r^{\prime }\in \pi ^{-1}(\mathscr{C})\cap U_\theta ^{\prime }$ to
a point $\widetilde r^{\prime \prime }\in \pi ^{-1}(\mathscr{C})\cap
U_{\theta ^{\prime }}^{\prime \prime }$, there is a natural orientation on $%
\beta ,$ from $\widetilde r^{\prime }$ to $\widetilde r^{\prime \prime }$.
If $\beta =\beta _1*\beta _2*\ldots *\beta _l,$ then this orientation
induces an orientation on each $\beta _j,$ $1\leq j\leq l.$

Let $\widetilde \gamma _j,$ $1\leq j\leq l$ be the extended lifts of
geodesics in $\mathscr{C}$ such that $\beta _j\subset \widetilde \gamma _j$.
We will say that $\beta \in C^{\prime }(U_\theta ^{\prime },U_{\theta
^{\prime }}^{\prime \prime })$ has a cycle if $\mathbb{D}\setminus \cup
_{j=1}^l\widetilde \gamma _j$ contains a bounded connected component.
Otherwise we will say that $\beta $ has no cycles.

\begin{figure}[!h]
	\centering
	\includegraphics[scale=0.2]{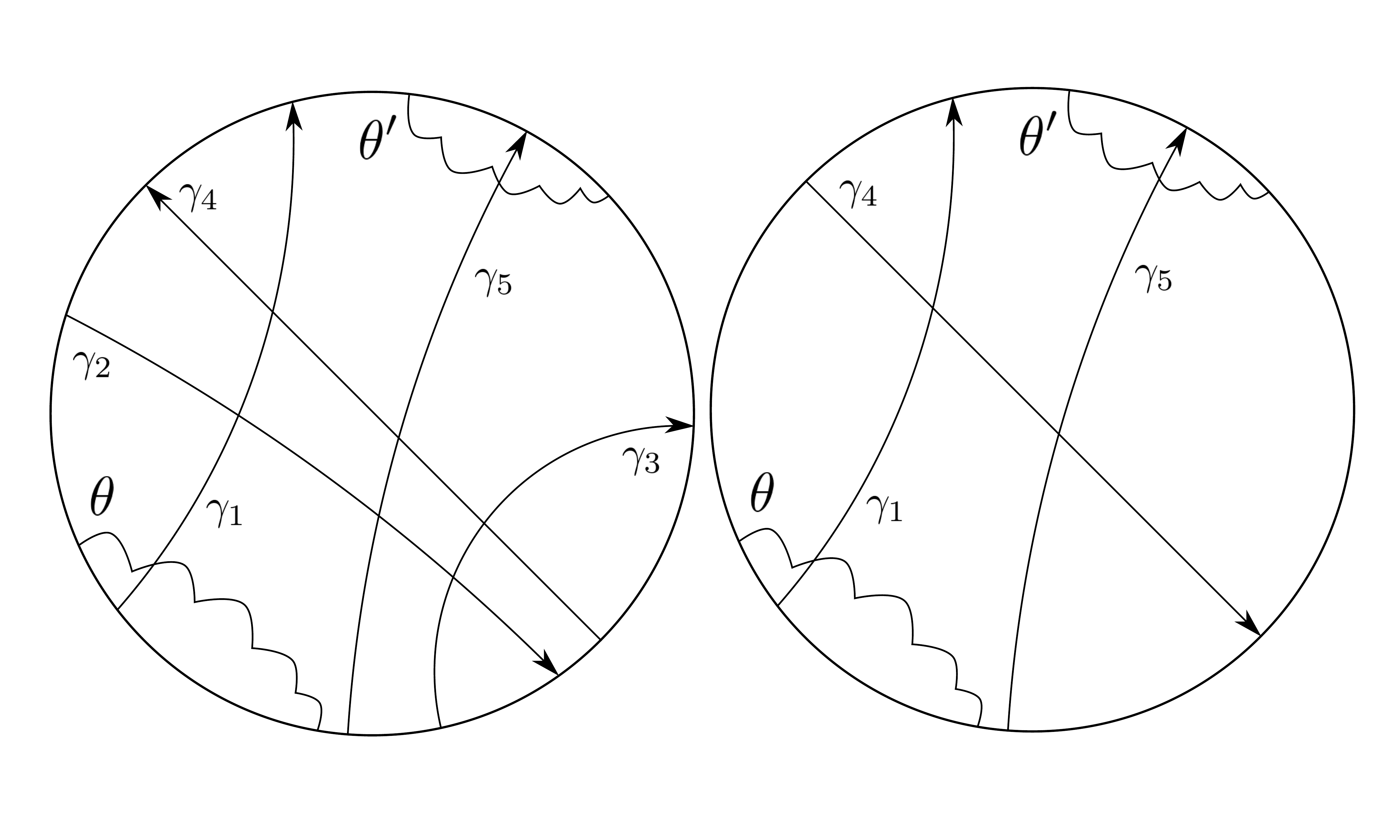}
	\caption{Example of how to eliminate cycles in $\beta$.}
	\label{cicle}
\end{figure}

Let $k=\min \{l\in \mathbb{N}|\beta \in C^{\prime }(U_\theta ^{\prime
},U_{\theta ^{\prime }}^{\prime \prime }),\beta =\beta _1*\beta _2*\ldots
*\beta _l\}$ and take $\beta _{\min }\in C^{\prime }(U_\theta ^{\prime
},U_{\theta ^{\prime }}^{\prime \prime })$ such that $\beta _{\min }=\beta
_{\min \_1}*\ldots *\beta _{\min \_k}$. We claim that the path $\beta _{\min
}$ has no cycles. Indeed, if it had, we could modify it and arrive at a path 
$\alpha \in C^{\prime }(U_\theta ^{\prime },U_{\theta ^{\prime }}^{\prime
\prime }),$ $\alpha $$=\alpha _1*\ldots *\alpha _{k^{\prime }},$ for some $%
k^{\prime }<k,$ a contradiction. See figure \ref{cicle}.

We will prove by induction on $k=\min \{l\in \mathbb{N}|\beta \in C^{\prime
}(U_\theta ^{\prime },U_{\theta ^{\prime }}^{\prime \prime }),\beta =\beta
_1*\beta _2*\ldots *\beta _l\}$ that $W^u(\widetilde p)\pitchfork %
W^s(g(\widetilde p)).$

Recall that, by definition of a fully essential system of curves $\mathscr{C}
$ (see definition \ref{thesystem}), for any geodesic $\gamma \in \mathscr{C}$
and any $\widetilde \gamma ,$ extended lift of that geodesic, there are
points $\widetilde p^{-}$ and $\widetilde p^{+}$ in $\mathbb{D}$
''following'' $\widetilde \gamma $ with both possible orientations, i.e., if 
$\alpha $ and $\omega $ are the endpoints of $\widetilde \gamma $ at $%
\partial \mathbb{D},$ then 
$$
\lim _{n\to \infty }\widetilde \phi ^n(\widetilde p^{+})=\omega =\lim _{n\to
\infty }\widetilde \phi ^{-n}(\widetilde p^{-})\text{ and }\lim _{n\to
\infty }\widetilde \phi ^n(\widetilde p^{-})=\alpha =\lim _{n\to \infty
}\widetilde \phi ^{-n}(\widetilde p^{+}). \text{ } 
$$

So, if $k=1,$ there exists $\beta _1\in C^{\prime }(U_\theta ^{\prime
},U_{\theta ^{\prime }}^{\prime \prime })$ and $\widetilde{\gamma }_1$ an
extended lift of a geodesic in $\mathscr{C}$ with $\beta _1\subset 
\widetilde{\gamma }_1$. It is clear that the orientation on $\widetilde{%
\gamma }_1$ induced by $\beta _1$ is from $\lambda _\theta ^{\prime }$ to $%
\lambda _{\theta ^{\prime }}^{\prime \prime }$. Associated to the extended
lift $\widetilde{\gamma }_1,$ there is a point $\widetilde{p}_1\in \pi
^{-1}(P)$ such that for some $n_1>0$ and some $h_1\in Deck(\pi )$ with $h_1(
\widetilde{\gamma }_1)=\widetilde{\gamma }_1,$ 
$$
\widetilde{\phi }^{n_1}(\widetilde{p}_1)=h_1(\widetilde{p}_1), 
$$
where 
$$
\lim _{m\to \infty }\widetilde{\phi }^{-mn_1}(\widetilde{p}_1)=\lim _{m\to
\infty }h_1^{-m}(\widetilde{p}_1)=\alpha _{h_1}\text{ and }\lim _{m\to
\infty }\widetilde{\phi }^{mn_1}(\widetilde{p}_1)=\lim _{m\to \infty }h_1^m(
\widetilde{p}_1)=\omega _{h_1}. 
$$
Note that $\alpha _{h_1}$ is the point at infinity of $\widetilde{\gamma }_1$
in $interior(\lambda _\theta ^{\prime })$ and $\omega _{h_1}$ is the point
at infinity of $\widetilde{\gamma }_1$ in $interior(\lambda _{\theta
^{\prime }}^{\prime \prime }).$ The point $\widetilde{p}_1$ can be chosen as
close as we want (in the Euclidean distance) to the point $\alpha _{h_1},$
something that forces $\widetilde{p}_1$ to belong to $U_\theta ^{\prime }$.
Since $\widetilde{\phi }|_{\partial \mathbb{D}}\equiv Id,$ for all $m>0$, $
\widetilde{\phi }^{mn_1}(\theta )$ is a path connected set in $\overline{%
\mathbb{D}}$ joining the points $\omega _{g_1},$ $\omega _{g_2}\in \partial 
\mathbb{D}$. As $\widetilde{p}_1\in U_\theta ^{\prime }$ and $\widetilde{%
\phi }$ preserves orientation, we get that for sufficiently large $m>0$, 
$$
\widetilde{\phi }^{mn_1}(U_\theta ^{\prime })\cap U_{\theta ^{\prime
}}^{\prime \prime }\neq \emptyset , 
$$
something that implies that 
$$
\widetilde{\phi }^{mn_1}(\theta )\pitchfork \theta ^{\prime }. 
$$

\begin{figure}[!h]
	\centering
	\includegraphics[scale=0.21]{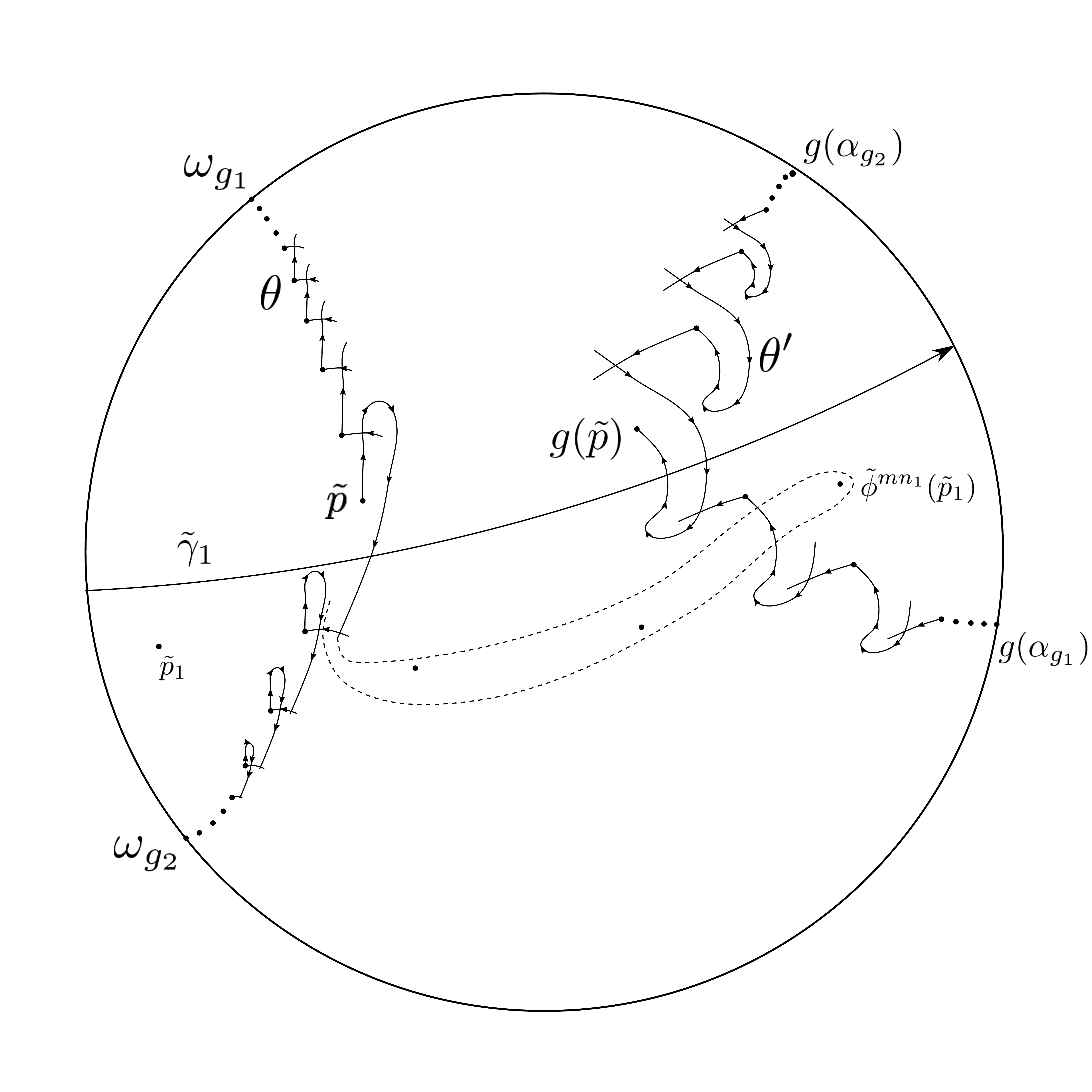}
	\caption{Intersection between $\widetilde{\phi}^{mn_1}(\theta)$ and $\theta'$.}
	\label{k1}
\end{figure}

Since the parts of $\theta \cap W^s(g_i^j(\widetilde p))$, $i=\{1,2\},$ $j>0$
are invariant under $\widetilde \phi $ and shrinking, we conclude that there
exists $i^{\prime },i^{\prime \prime }\in \{1,2\}$ and $j^{\prime
},j^{\prime \prime }>0$ such that $W^u(g_{i^{\prime }}^{j^{\prime
}}(\widetilde p))\pitchfork W^s(gg_{i^{\prime \prime }}^{-j^{\prime \prime
}}g^{-1}(g(\widetilde p)))$ (see figure \ref{k1}). So by the $C^0$ $\lambda $%
-lemma we get that 
$$
W^u(\widetilde p)\pitchfork W^s(g(\widetilde p)). 
$$

If $k=2$ there exist $\beta \in C^{\prime }(U_\theta ^{\prime },U_{\theta
^{\prime }}^{\prime \prime }),$ $\beta =\beta _1*\beta _2$ and $\widetilde
\gamma _1,\widetilde \gamma _2$ extended lifts of geodesics in $\mathscr{C}$
such that $\beta _1\subset \widetilde \gamma _1$ and $\beta _2\subset
\widetilde \gamma _2$. Moreover, in the same way as in the previous case,
for $i\in \{1,2\}$ there exists a point $\widetilde p_i\in \pi ^{-1}(P)$ and 
$h_i\in Deck(\pi )$ that leaves $\widetilde \gamma _i$ invariant, such that
for some $n_i>0,$ 
$$
\widetilde \phi ^{n_i}(\widetilde p_i)=h_i(\widetilde p_i). 
$$

\begin{figure}[!h]
	\centering
	\includegraphics[scale=0.21]{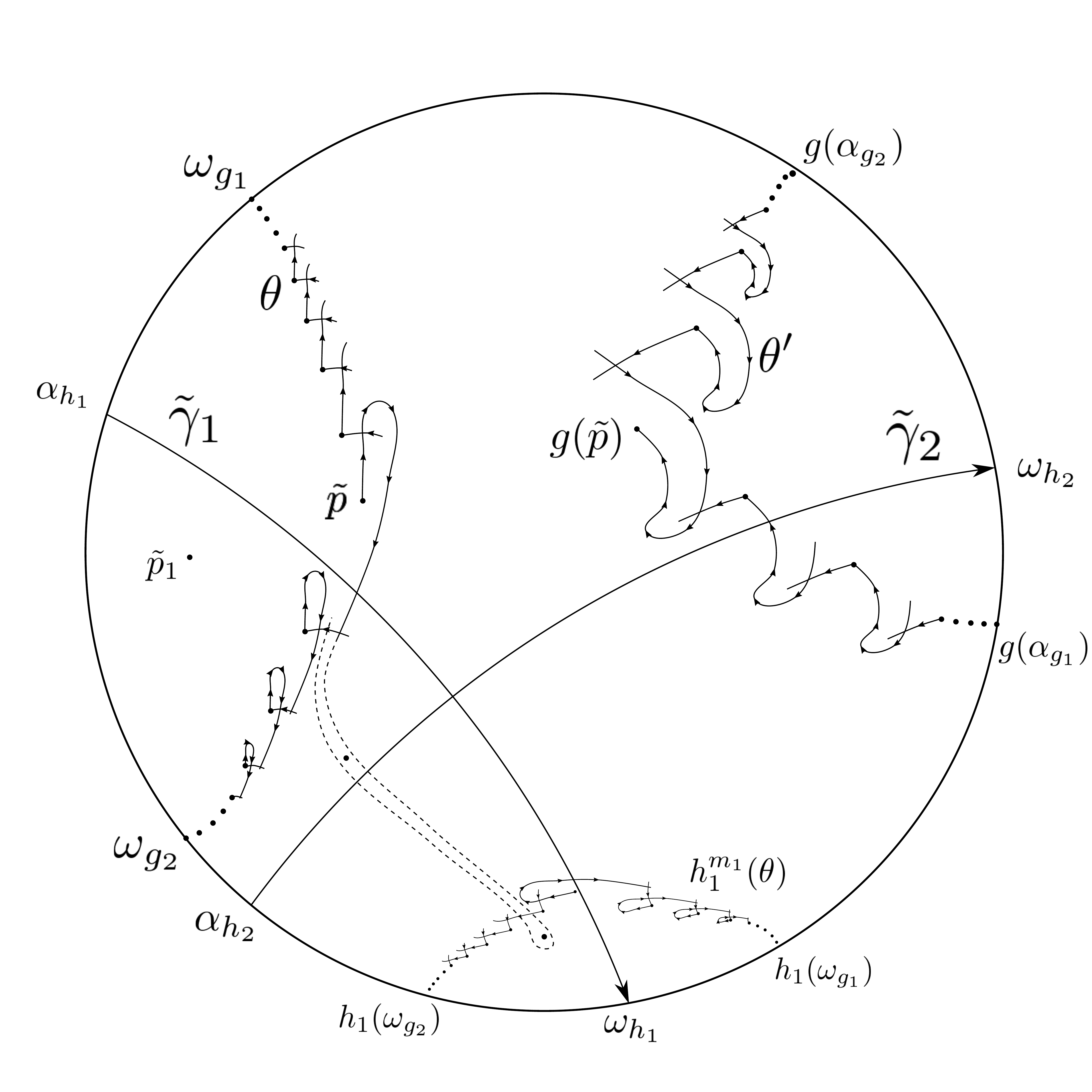}
	\caption{Intersection between $\widetilde{\phi}^{mn_1}(\theta)$ and $h_1^{m_1}(\theta)$.}
	\label{k2}
\end{figure}

The points $\alpha _{h_1},\omega _{h_1}$ separate the points $\alpha
_{h_2},\omega _{h_2}$ at $\partial \mathbb{D}$ and $\alpha _{h_1}$ is in the
interior of $\lambda _\theta ^{\prime }$ and $\omega _{h_2}$ is in the
interior of $\lambda _{\theta ^{\prime }}^{\prime \prime },$ see figure \ref
{k2}. Let us consider a sufficiently large $m_1>0$ in a way that $%
h_1^{m_1}(\theta )$ is close (in the Euclidean distance) to the point $%
\omega _{h_1}$ and $\theta \cap h_1^{m_1}(\theta )=\emptyset .$ In
particular, the points $h_1^{m_1}(\omega _{g_1})$ and $h_1^{m_1}(\omega
_{g_2})$ are very close to $\omega _{h_1}$.

Exactly as in the case $k=1,$ $\widetilde{p}_1$ can be chosen sufficiently
close to $\alpha _{h_1}$ in the Euclidean distance, in a way that $
\widetilde{p}_1\in U_\theta ^{\prime },$ and then for a sufficiently large $%
m>0,$ $\widetilde{\phi }^{mn_1}(\widetilde{p}_1)$ is so close to $\omega
_{h_1},$ something that forces $\widetilde{\phi }^{mn_1}(\theta )\pitchfork %
h_1^{m_1}(\theta ).$

Exactly as in the previous case, there exists $i^{\prime },i^{\prime \prime
}\in \{1,2\}$ and $j^{\prime },j^{\prime \prime }>0$ such that 
\begin{equation}
\label{its} 
\begin{array}{c}
W^u(g_{i^{\prime }}^{j^{\prime }}(\widetilde p)) 
\pitchfork W^s(h_1^{m_1}g_{i^{\prime \prime }}^{j^{\prime \prime
}}(\widetilde p)), \\ \text{something that implies, again by the }C^0\text{ }%
\lambda \text{-lemma, that} \\ W^u(\widetilde p)\pitchfork %
W^s(h_1^{m_1}g_{i^{\prime \prime }}^{j^{\prime \prime }}(\widetilde p)). 
\end{array}
\end{equation}

Note that since $k=2,$ $\alpha _{h_2}$ is not in the interior of $\lambda
_\theta ^{\prime },$ otherwise there would be a path in $C^{\prime
}(U_\theta ^{\prime },U_{\theta ^{\prime }}^{\prime \prime })$ with size $1,$
a contradiction with the fact that $k=2.$ So, we can always choose $i_0\in
\{1,2\}$ and construct a new path connected set $\theta _{h_1}$ using $\eta
_{i_0}$ of expression (\ref{construct}) and an analogous construction
obtained from expression (\ref{its}). Before getting into details, we
emphasize that the choice of $i_0$ is very important, because $\alpha _{h_2}$
is not in the interior of $\lambda _\theta ^{\prime },$ but it could be one
of its endpoints. So, if $\alpha _{h_2}$ is one of the endpoints of $\lambda
_\theta ^{\prime },$ $g_{i_0}$ must be chosen associated with the other
endpoint of $\lambda _\theta ^{\prime }.$

In order to construct $\theta _{h_1},$ first consider the set $\eta _{i_0}$
of expression (\ref{construct}) associated to $g_{i_0}$ chosen as before.
Since expression (\ref{its}) holds, there exists a path $\eta $ joining $%
\widetilde p$ to $h_1^{m_1}g_{i^{\prime \prime }}^{j^{\prime \prime
}}h_1^{-m_1}(h_1^{m_1}(\widetilde p))$ as follows: $\eta $ starts at $%
\widetilde p,$ consists of a compact connected piece of $\widetilde{\lambda }%
_u$ until it reaches $h_1^{m_1}g_{i^{\prime \prime }}^{j^{\prime \prime
}}h_1^{-m_1}(h_1^{m_1}(\widetilde{\beta }_s))$ 
and then it continues as a compact connected piece of $h_1^{m_1}g_{i^{\prime
\prime }}^{j^{\prime \prime }}h_1^{-m_1}(h_1^{m_1}(\widetilde{\beta }_s))$
until it reaches $h_1^{m_1}g_{i^{\prime \prime }}^{j^{\prime \prime
}}h_1^{-m_1}(h_1^{m_1}(\widetilde p)).$ Note that $i^{\prime \prime }\in
\{1,2\} $ was defined before expression (\ref{its}). Finally, pick the $\eta
_{i^{\prime \prime }}$ as in expression (\ref{construct}).

Then, define

\begin{equation}
\label{h1}\theta _{h_1}=\bigg(\bigcup_{i\geq 0}g_{i_0}^i(\eta _{i_0})\bigg) %
\cup \eta \cup \bigg(\bigcup_{j\geq j^{\prime \prime }}h_1^{m_1}g_{i^{\prime
\prime }}^jh_1^{-m_1}(h_1^{m_1}(\eta _{i^{\prime \prime }}))\bigg). 
\end{equation}

The new path connected set $\theta _{h_1}$ has similar properties to $\theta
,$ and as was explained for $\theta $ and $\theta ^{\prime },$ it can be
understood as the concatenation of two curves in $\mathbb{D},$ one joining $
\widetilde{p}$ to $\omega _{g_{i_0}}$ and another joining $\widetilde{p}$ to 
$h_1^{m_1}(\omega _{g_{i^{\prime \prime }}}).$ If $\theta _{h_1}\pitchfork %
\theta ^{\prime },$ then $W^u(\widetilde{p})\pitchfork W^s(g(\widetilde{p})).
$ And in case there is no topologically transverse intersection between $%
\theta _{h_1}$ and $\theta ^{\prime },$ as $\min \{l\in \mathbb{N}|\beta \in
C^{\prime }(U_{\theta _{h_1}}^{\prime },U_{\theta ^{\prime }}^{\prime \prime
}),\beta =\beta _1*\beta _2*\ldots *\beta _l\}=1,$ the situation is reduced
to the previous case. This happens because $\alpha _{h_2}$ is in the
interior of $\lambda _{\theta _{h_1}}^{\prime }$ if $i_0$ is chosen as
explained above. Hence arguing exactly as in case $k=1,$ we conclude that $%
W^u(\widetilde{p})\pitchfork W^s(g(\widetilde{p})).$

Now by induction suppose the result is true for 
$$
\min \{l\in \mathbb{N}|\beta \in C^{\prime }(U_\theta ^{\prime },U_{\theta
^{\prime }}^{\prime \prime }),\beta =\beta _1*\beta _2*\ldots *\beta
_l\}=1,2,...,k-1 
$$
and let us prove that it holds for $k.$ We can assume that $k\geq 3.$ Fix $%
\beta \in C^{\prime }(U_\theta ^{\prime },U_{\theta ^{\prime }}^{\prime
\prime })$ with $\beta =\beta _1*\beta _2*\beta _3*\ldots *\beta _k.$
Remember that we are assuming that there exists no path in $C^{\prime
}(U_\theta ^{\prime },U_{\theta ^{\prime }}^{\prime \prime })$ using less
than $k$ geodesics. Let $\widetilde{\gamma }_i,$ $0\leq i\leq k,$ be the
extended lifts of the geodesics in $\mathscr{C}$ such that $\beta _i\subset 
\widetilde{\gamma }_i$. For each $i\in \{1,...,k\},$ as the $\beta
_{i^{\prime }s}$ are oriented, there exists a point $\widetilde{p}_i\in \pi
^{-1}(P)$ and $h_i\in Deck(\pi )$ that leaves $\widetilde{\gamma }_i$
invariant and moves points according to the orientation of $\widetilde{%
\gamma }_i$, such that for some $n_i>0$ 
$$
\widetilde{\phi }^{n_i}(\widetilde{p}_i)=h_i(\widetilde{p}_i). 
$$

We claim that the following properties are true:

\begin{itemize}
\item  for every $2\leq i\leq k,$ $\alpha _{h_i}$ and $\omega _{h_i}$ are
not in the interior of $\lambda _\theta ^{\prime }$. If this happened, it
would be possible to create a path $\beta ^{\prime }\in C^{\prime }(U_\theta
^{\prime },U_{\theta ^{\prime }}^{\prime \prime })$ using less than $k$
geodesics.

\item  for every $3\leq i\leq k,$ $\alpha _{h_i}$ and $\omega _{h_i}$ are
not inside the segment $[\alpha _{h_1},\omega _{h_1}]$ in $\partial 
\mathbb{D}$ delimited by $\alpha _{h_1}$ and $\omega _{h_1}$ and containing
the point $\alpha _{h_2}.$ This follows from the fact that $\beta $ has no
cycles.
\end{itemize}

Now we can proceed as in case $k=2$ and construct the path connected set $%
\theta _{h_1}$ in the same way as in (\ref{h1}). Here we just have to be
careful and choose an integer $m_1>0$ sufficiently large so that in the
segment $[h_1^{m_1}(\omega _{g_1}),h_1^{m_1}(\omega _{g_2})]$ of $\partial 
\mathbb{D}$ delimited by $h_1^{m_1}(\omega _{g_1})$ and $h_1^{m_1}(\omega
_{g_2}),$ and containing the point $\omega _{h_1},$ there are no others $%
\alpha ^{\prime }s$ and $\omega ^{\prime }s.$ But since there is a finite
number of $\alpha ^{\prime }s$ and $\omega ^{\prime }s,$ this is always
possible. Note that $\lambda _{\theta _{h_1}}^{\prime }\subset \lambda
_\theta ^{\prime }\cup [\alpha _{h_1},\omega _{h_1}]\cup [h_1^{m_1}(\omega
_{g_1}),h_1^{m_1}(\omega _{g_2})],$ and using the properties proved before,
we conclude that $\alpha _{h_2}$ is in the interior of $\lambda _{\theta
_{h_1}}^{\prime }$. Moreover, $\lambda _{\theta _{h_1}}^{\prime }$ contains
no other $\alpha _{h_i}$ or $\omega _{h_i},$ $2\leq i\leq k.$ If we pick a
point in $\widetilde \gamma _2$ close to $\alpha _{h_2}$ and make $\beta
_2^{\prime }$ the subarc of $\widetilde \gamma _2$ joining this point to the
intersection point of $\widetilde \gamma _2$ and $\widetilde \gamma _3,$ we
get that $\beta _{mod}=\beta _2^{\prime }*\beta _3*\ldots *\beta _k$ is a
path in $C^{\prime }(U_{\theta _{h_1}}^{\prime },U_{\theta ^{\prime
}}^{\prime \prime })$ formed by $k-1$ subarcs of extended lifts of geodesics
in $\mathscr{C}$ that has no cycles. So using the induction hypothesis we
conclude that $W^u(\widetilde p)\pitchfork W^s(g(\widetilde p)).$ %

This proves that for all $g\in Deck(\pi )$, $g\neq Id$, $W^u(\widetilde p)%
\pitchfork W^s(g(\widetilde p))$. In order to deal with $g=Id,$ consider
some $h,h^{-1}\in Deck(\pi ),$ $h\neq Id.$ Then $W^u(\widetilde p)\pitchfork %
W^s(h(\widetilde p))$ and $W^u(\widetilde p)\pitchfork W^s(h^{-1}(\widetilde
p))$. As $\widetilde \phi $ commutes with all deck transformations, $%
W^u(h^{-1}(\widetilde p))\pitchfork W^s(\widetilde p)$ and so, by the $C^0$ $%
\lambda $-lemma, $W^u(\widetilde p)\pitchfork W^s(\widetilde p)$.

Actually, as we said in the beginning, we proved something a little bit
stronger: for all $g\in Deck(\pi )$, $\widetilde{\lambda }_u\pitchfork g( 
\widetilde{\beta }_s).$

%

\qed

\subsection{Handel's global shadowing}

The next result tells us that, as our map $\phi $ is pseudo-Anosov relative
to some finite invariant set, the complicated dynamics of $\phi $ is in a
certain sense inherited by $f,$ see \cite{mh85}.

\begin{theorem}[Handel's global shadowing]
If $f:S\to S$ is a homeomorphism of a closed orientable surface $S,$
homotopic to the identity, $P$ is a finite $f$-invariant set and $f$ is
isotopic relative to $P$ to some map $\phi :S\to S$ which is pseudo-Anosov
relative to $P,$ then there exists a compact $f$-invariant set $W\subset S$
and a continuous surjection $s:W\to S$ that is homotopic to the inclusion
map $i:W\to S,$ such that $s$ semi-conjugates $f|_W$ to $\phi ,$ that is, $%
s\circ f|_W=\phi \circ s.$
\end{theorem}

Observe that, as $s:W\to S$ is homotopic to the inclusion map $i:W\to S,$ $s$
has a lift $\widetilde s:\pi ^{-1}(W)\to \mathbb{D},$ such that 
$$
\widetilde s\circ \widetilde f|_{\pi ^{-1}(W)}=\widetilde \phi \circ
\widetilde s, 
$$
where $\widetilde{\phi }$ and $\widetilde{f}$ are the natural lifts of $\phi 
$ and $f,$ and sup$\{d_{\mathbb{D}}(\widetilde s(\widetilde q),\widetilde
q)|\widetilde q\in \pi ^{-1}(W)\}<C_f,$ for some constant $C_f>0.$

\subsection{Special horseshoes for the pseudo-Anosov map $\phi $}

In this subsection we prove a simple lemma used in the proofs of theorems 1
and 2. The setting is the following: let $f:S\to S$ be a homeomorphism
isotopic to the identity with a fully essential system of curves $\mathscr{C}
$ and let $P$ be the set of periodic points associated with the geodesics in 
$\mathscr{C}.$ From lemma \ref{parel} we know that there exists an integer $%
m_0>0$ such that $f^{m_0}$ is isotopic relative to $P$ to $\phi :S\to S,$ a
homeomorphism which is pseudo-Anosov relative to $P$ and isotopic to the
identity as a homeomorphism of $S.$ From lemma \ref{lema2} there exists a
contractible hyperbolic $\phi $-periodic point $p\in S,$ such that for any $%
\widetilde p\in \pi ^{-1}(p)$ and any given $g\in Deck(\pi ),$ $%
W^u(\widetilde p)\pitchfork W^s(g(\widetilde p)).$ In $\mathbb{D}$ we are
considering the natural lift of $\phi ,$ denoted $\widetilde{\phi }.$ As we
already done before, without loss of generality, assume that 
$p$ is fixed under $\phi $ and also all four branches at $p$ are invariant
under $\phi .$

\begin{lemma}
\label{basic} For any $g\in Deck(\pi )$ and any fundamental domain $
\widetilde{Q}\subset \mathbb{D}$ of $S$ such that $\widetilde{p}=\pi
^{-1}(p)\cap \widetilde{Q}$ is in the interior of $\widetilde{Q}$, there
exists arbitrarily small rectangles $\widetilde{R}\subset \widetilde{Q}$
such that:
\end{lemma}

\begin{enumerate}
\item  $interior(\widetilde{R})$ contains $\widetilde{p},$ and two sides of $
\widetilde{R}$ are very close to an arc $\widetilde{\beta }$, \ $\widetilde{p%
}\in \widetilde{\beta }\subset W^s(\widetilde{p})\cap R$ and the two other
sides of $\widetilde{R}$ have very small length;

\item  for some $N>0,$ $\widetilde{\phi }^N(\widetilde{R})\cap \widetilde{R}%
\supset \widetilde{R}_0$ and $\widetilde{\phi }^N(\widetilde{R})\cap g(
\widetilde{R})\supset \widetilde{R}_1,$ where $R_0=\pi (\widetilde{R}_0),$ $%
R_1=\pi (\widetilde{R}_1)$ are rectangles contained in $R=\pi (\widetilde{R})
$ which have two sides contained in the sides of $R$ which are very close to 
$\beta =\pi (\widetilde{\beta })$ and two sides contained in the interior of 
$R.$
\end{enumerate}

\begin{figure}[!h]
	\centering
	\includegraphics[scale=0.2]{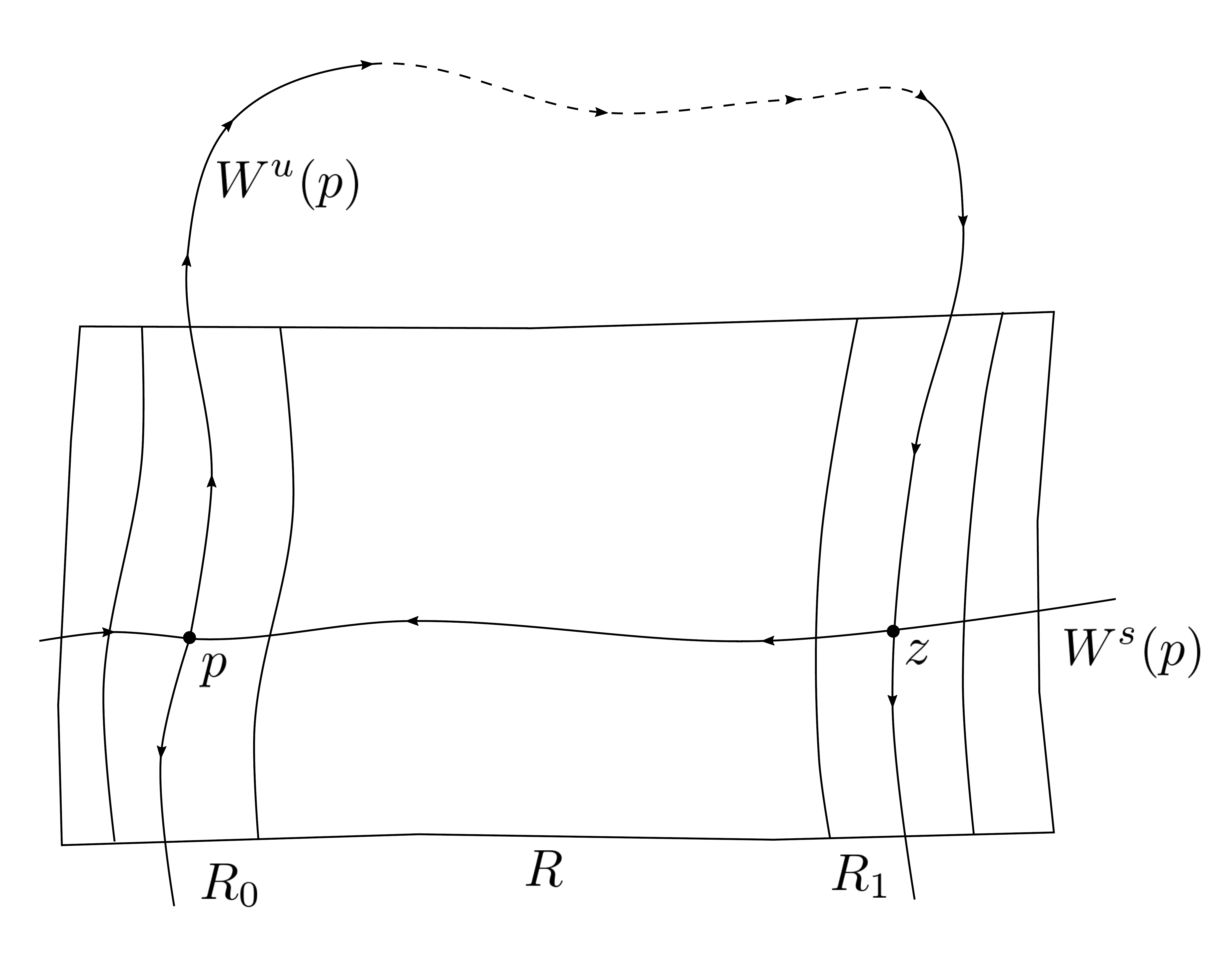}
	\caption{Horseshoe associated to $W^u(\widetilde{p}) \pitchfork W^s(g(\widetilde{p}))$.}
	\label{horseshoe1}
\end{figure}

{\it Proof: }Fix $g\in Deck(\pi )$ and a fundamental domain $\widetilde{Q}%
\subset \mathbb{D}$ of $S.$ Consider $\widetilde{p}\in interior(\widetilde{Q}%
)$ the $\widetilde{\phi }$-hyperbolic periodic point given by lemma \ref
{lema2}. As $W^u(\widetilde{p})\pitchfork W^s(g(\widetilde{p})),$ if we
consider the projection of this intersection to the surface $S,$ we obtain a
horseshoe associated to $g$ as follows:

Let $\widetilde{z}$ be a point in $W^u(\widetilde{p})\cap W^s(g(\widetilde{p}%
))$ as close as we want to $g(\widetilde{p})$ and let $R\subset S$ be a
rectangle, whose interior contains the connected arc $\beta \subset W^s(p)$
whose endpoints are $p=\pi (\widetilde{p})$ and $z=\pi (\widetilde{z}).$ As $%
p$ is in the interior of $R,$ the rectangle $R$ also contains an arc in $%
W^u(p)$ with one endpoint at $p.$ As usual, $R$ is very thin in the unstable
direction, being very close to $\beta .$ This implies that for some large $%
N>0,$ $\phi ^N(R)\cap R\supseteq R_0\cup R_1,$ where $p\in R_0$ and $z\in R_1
$ are rectangles, each of them exactly as in the statement of the lemma. If $
\widetilde{R}$ is the connected component of $\pi ^{-1}(R)$ containing $
\widetilde{p},$ then clearly $\widetilde{\phi }^N(\widetilde{R})\cap 
\widetilde{R}\supset \widetilde{R}_0$ and $\widetilde{\phi }^N(\widetilde{R}%
)\cap g(\widetilde{R})\supset \widetilde{R}_1,$ for some connected
components $\widetilde{R}_0$ of $\pi ^{-1}(R_0)$ and $\widetilde{R}_1$ of $%
\pi ^{-1}(R_1).$ See figure \ref{horseshoe1}. Clearly, there may be other
rectangles in $\phi ^N(R)\cap R.$

\qed

\subsection{The $C^{1+\epsilon }$ case: Some background in Pesin theory}

\label{pesin}

In this subsection, assume that $f:S\to S$ is a $C^{1+\epsilon }$
diffeomorphism, for some $\epsilon >0.$ Recall that a $f$-invariant Borel
probability measure $\mu $ is hyperbolic if all the Lyapunov exponents of $f$
are non-zero at $\mu $-almost every point (for instance, see the supplement
of \cite{akbh}). The next paragraphs were taken from \cite{asmp03}. They
consist of an informal description of the theory of non-uniformly hyperbolic
systems, together with some definitions and lemmas from \cite{asmp03}.

Let $\mu $ be a non-atomic hyperbolic ergodic $f$-invariant Borel
probability measure. Given $0<\delta <1,$ there exists a compact Pesin set $%
\Lambda _\delta $ with $\mu (\Lambda _\delta )>1-\delta ,$ having the
following properties: for every $p\in \Lambda _\delta ,$ there exists an
open neighborhood $U_p,$ a compact neighborhood $V_p\subset U_p$ and a
diffeomorphism $F:(-1,1)^2\to U_p,$ with $F(0,0)=p$ and $%
F([-1/10,1/10]^2)=V_p,$ such that the local unstable manifolds $W_{loc}^u(q)$
of all points $q$ in $\Lambda _\delta \cap V_p$ are the images under $F$ of
graphs of the form $\{(x,F_2(x))|x\in (-1,1)\},$ $F_2$ a function with small
Lipschitz constant. Any two such local unstable manifolds are either
disjoint or equal and they depend continuously on the point $q\in \Lambda
_\delta \cap V_p.$ Similarly, the local stable manifolds $W_{loc}^s(q)$ of
points $q\in \Lambda _\delta \cap V_p$ are the images under $F$ of graphs of
the form $\{(F_1(y),y)|y\in (-1,1)\},$ $F_1$ a function with small Lipschitz
constant. Any two such local stable manifolds are either disjoint or equal
and they depend continuously on the point $q\in \Lambda _\delta \cap V_p.$

It follows that there exists a continuous product structure in $\Lambda
_\delta \cap V_p:$ given any $r,r^{\prime }\in \Lambda _\delta \cap V_p,$
the intersection $W_{loc}^u(r)\cap W_{loc}^s(r^{\prime })$ is transversal
and consists of exactly one point, which will be denoted $[r,r^{\prime }]$.
This intersection varies continuously with the two points and may not be in $%
\Lambda _\delta .$ Hence we can define maps $P_p^s:\Lambda _\delta \cap
V_p\to W_{loc}^s(p)$ and $P_p^u:\Lambda _\delta \cap V_p\to W_{loc}^u(p)$ as 
$P_p^s(q)=[q,p]$ and $P_p^u(q)=[p,q].$

Let $R^{\pm }$ denote the set of all points in $S$ which are both forward
and backward recurrent. By the Poincar\'e recurrence theorem, $\mu (R^{\pm
})=1.$

\begin{description}
\item[Definition]  \label{inacess}{\bf (Accessible and inaccessible points) 
\ref{inacess}.} A point $p\in \Lambda _\delta \cap V_p\cap R^{\pm }$ is
inaccessible if it is accumulated on both sides of $W_{loc}^s(p)$ by points
in $P_p^s(\Lambda _\delta \cap V_p\cap R^{\pm })$ and also accumulated on
both sides of $W_{loc}^u(p)$ by points in $P_p^u(\Lambda _\delta \cap
V_p\cap R^{\pm })$. Otherwise, $p$ is accessible.
\end{description}

After this definition, we can state two lemmas from \cite{asmp03} about
accessible and inaccessible points and the relation between these points and
hyperbolic periodic points close to them.

\begin{lemma}
\label{andre1} Let $q\in \Lambda _\delta \cap V_p\cap R^{\pm }$ be an
inaccessible point. Then there exist rectangles enclosing $q,$ having sides
along the invariant manifolds of hyperbolic periodic saddles in $V_p$ and
having arbitrarily small diameter.
\end{lemma}

A rectangle is a Jordan curve made up of alternating segments of stable and
unstable manifolds, two of each. The segments forming the boundary are its
sides and the intersection points of the sides are the corners. A rectangle
is said to enclose $p$ if it is the boundary of an open topological disk
containing $p.$

\begin{lemma}
\label{andre2} The subset of accessible points in $\Lambda _\delta \cap
V_p\cap R^{\pm }$ has zero $\mu $ measure.
\end{lemma}

Another concept that will be a crucial hypothesis for us is positive
topological entropy. In the following we describe why.

When the topological entropy $h(f\mid _K)$ is positive, for some compact $f$%
-invariant set $K,$ by the variational principle, there exists a $f$%
-invariant Borel probability measure $\mu _0$ with $supp(\mu _0)\subset K$
and $h_{\mu _0}(f)>0.$ Using the ergodic decomposition of $\mu _0$ we find
an extremal point $\mu $ of the set of Borel probability $f$-invariant
measures, such that $supp(\mu )$ is also contained in $K$ and $h_\mu (f)>0.$
Since the extremal points of this set are ergodic measures, $\mu $ is
ergodic. The ergodicity and the positiveness of the entropy imply that $\mu $
has no atoms and applying the Ruelle inequality to $f$ we see that $\mu $
has a positive Lyapunov exponent, see \cite{akbh}. Working with $f^{-1}$ and
using the fact that $h_\mu (f^{-1})=h_\mu (f)>0,$ we see that $f^{-1}$ must
also have a positive Lyapunov exponent with respect to $\mu ,$ which is the
negative of the negative Lyapunov exponent for $f$.

Hence when $K$ is a compact $f$-invariant set and the topological entropy of 
$f\mid _K$ is positive, there always exists a ergodic, non-atomic, invariant
measure supported on $K,$ with non-zero Lyapunov exponents, one positive and
one negative, the measure having positive entropy: A hyperbolic measure.

The existence of this kind of measure will be important for us because of
the following theorem, that can be proved combining the main lemma and
theorem $4.2$ of \cite{ak80}.

\begin{theorem}
\label{katok} Let $f$ be a $C^{1+\epsilon }$ (for some $\epsilon >0$)
diffeomorphism of a surface $M$ and suppose $\mu $ is an ergodic hyperbolic
Borel probability $f$-invariant measure with $h_\mu (f)>0$ and compact
support. Then, for any $\alpha >0$ and any $p\in supp(\mu )$, there exists a
hyperbolic periodic point $q\in B_\alpha (p)$ which has a transversal
homoclinic intersection, and the whole orbit of $q$ is contained in the $%
\alpha $-neighborhood of $supp(\mu )$.
\end{theorem}

\section{Proof of theorem 1}

Let $f:S\to S$ be a homeomorphism satisfying the theorem hypotheses. If we
remember subsection 2.8 and lemma \ref{basic}, for any fixed $g\in Deck(\pi)$
and any fundamental domain $\widetilde{Q}\subset \mathbb{D}$ of $S,$ there
exist arbitrarily small rectangles $R\subset S$ such that a connected
component $\widetilde{R}$ of $\pi^{-1}(R)$ is contained in $interior( 
\widetilde{Q})$ (we may have to perturb $\widetilde{Q}$ a little bit) and
for some $N>0,$ $\phi ^N(R)\cap R\supseteq R_0\cup R_1.$ Associated with
this horseshoe, if we consider the $\phi ^N$-fixed point $q\in R_1,$ then
for $\widetilde q=\pi ^{-1}(q)\cap \widetilde R$ the following holds: 
$$
\widetilde \phi ^N(\widetilde q)=g(\widetilde q)\Rightarrow \text{ for all }%
j>0,\widetilde \phi ^{jN}(\widetilde q)=g^j(\widetilde q). 
$$

Let $s:W\to S$ be the semi-conjugacy given by Handel's global shadowing and $%
\widetilde s:\pi ^{-1}(W)\to \mathbb{D}$ its lift which relates the natural
lifts $\widetilde f$ and $\widetilde \phi $. Fix $\widetilde z\in \widetilde
s^{-1}(\widetilde q) $. Since $\widetilde s\circ \widetilde f(\widetilde
z)=\widetilde \phi \circ \widetilde s(\widetilde z),$ we get that 
$$
\widetilde s(\widetilde f^{jN}(\widetilde z))=\widetilde \phi
^{jN}(\widetilde s(\widetilde z))=\widetilde \phi ^{jN}(\widetilde
q)=g^j(\widetilde q). 
$$

As we explained in subsection 2.7, the fact that $s$ is isotopic to the
inclusion implies the existence of $C_f>0$ such that $d_{\mathbb{D}%
}(\widetilde s(\widetilde w),\widetilde w)<C_f,$ for all $\widetilde w\in
\pi ^{-1}(W).$ In particular 
$$
d_{\mathbb{D}}(\widetilde f^{jN}(\widetilde z),\widetilde s(\widetilde
f^{jN}(\widetilde z)))=d_{\mathbb{D}}(\widetilde f^{jN}(\widetilde
z),g^j(\widetilde q))<C_f,\text{ for all }j>0. 
$$

As $g^{-1}\in Deck(\pi )$ is an isometry of $d_{\mathbb{D}},$ 
$$
d_{\mathbb{D}}(\widetilde f^{jN}(\widetilde z),g^j(\widetilde q))=d_{%
\mathbb{D}}(g^{-j}(\widetilde f^{jN}(\widetilde z)),\widetilde q)<C_f. 
$$

This means that for any $\widetilde z\in \widetilde s^{-1}(\widetilde q)$
and for all $j>0,$ $(g^{-1}\widetilde f^N)^j(\widetilde z)\in
B_{C_f}(\widetilde q)$. So the positive orbit of $\widetilde z$ with respect
to $g^{-1}\widetilde f^N$ is bounded. Thus, defining $\widetilde{K}_g$ as
the $\omega $-limit set of the point $\widetilde z$ under $g^{-1}\widetilde{f%
}^N,$ $\widetilde{K}_g$ is a compact $g^{-1}\widetilde f^N$-invariant set
contained in $V_{C_f}(Q),$ and hence $\widetilde f^N( \widetilde{K}_g)=g( 
\widetilde{K}_g)$. By Brouwer's lemma on translation arcs \cite{franks}, $%
g^{-1}\widetilde f^N$ has a fixed point, that is, there exists $\widetilde
r\in \mathbb{D}$ with $g^{-1}\widetilde f^N(\widetilde r)=\widetilde r$, and
so 
$$
\widetilde f^N(\widetilde r)=g(\widetilde r). 
$$

\qed

\section{Proof of theorem 2}

Let $f:S\to S$ be a $C^{1+\epsilon }$ diffeomorphism isotopic to the
identity with a fully essential system of curves $\mathscr{C}.$ As in
theorem 1, let $\phi :S\to S$ be the pseudo-Anosov map relative to $P,$
which is isotopic to $f^{m_0}$  relative to $P$ (for some $m_0>0,$ which as
before, to simplify the notation, we assume to be $1).$ The finite set $P$
is the set of periodic points associated with the geodesics in $\mathscr{C}.$
By lemma \ref{lema2}, for any given $g\in Deck(\pi )$ and any fundamental
domain $\widetilde{Q}$ of $S,$ if $\widetilde{\phi }:\mathbb{D}\to \mathbb{D}
$ is the natural lift of $\phi ,$ there exists a hyperbolic $\widetilde{\phi 
}$-periodic point $\widetilde{p}\in \widetilde{Q}\subset \mathbb{D}$ such
that 
$$
W^u(\widetilde{p})\pitchfork W^s(g(\widetilde{p})). 
$$

Again, as we did in previous results, without loss of generality,
considering an iterate of $\widetilde \phi $ if necessary, assume that $%
\widetilde p$ is fixed under $\widetilde \phi $ and each branch at $
\widetilde{p}$ is also $\widetilde{\phi}$-invariant.

Using lemma \ref{basic}, if the transversal intersection $W^u(\widetilde p)%
\pitchfork W^s(g(\widetilde p))$ at some $\widetilde z\in \mathbb{D}$ is
projected to the surface $S,$ it corresponds to a transversal homoclinic
point $z=\pi (\widetilde z)\in W^u(p)\cap W^s(p).$ Associated with this
intersection, a horseshoe in $S$ can be obtained, i.e. on the surface there
exist a small rectangle $R,$ containing the arc $\beta $ in $W^s(p)$ from $p$
to $z$ (as always $R$ is very close to $\beta ),$ and a positive integer $%
N>0 $ such that $\phi ^N(R)\cap R\supseteq R_0\cup R_1$, where $R_0$ is a
rectangle inside $R$ containing $p$ and $R_1$ is another rectangle inside $R$
containing $z,$ see figure \ref{horseshoe1}.

As $\widetilde z\in W^u(\widetilde p)\cap W^s(g(\widetilde p))$ can be
chosen as close as we want to $g(\widetilde p),$ the rectangle $R\subset S$
can be chosen small enough so that all the singular points of the stable and
unstable foliations of $\phi $ do not belong to $R.$ Moreover, considering
the compact set $\Omega =\cap _{k\in \mathbb{Z}}\phi ^{kN}(R_0\cup R_1),$ we
know that if $R$ is sufficiently close to $\beta ,$ to every bi-infinite
sequence in $\{0,1\}^{\mathbb{Z}}$ denoted $(a_n)_{n\in \mathbb{Z}}$ there
is a single point $z_{*}\in \Omega $ which realizes it, that is $\phi
^{kN}(z_{*})$ belongs to $R_{a_k}$ for all integers $k.$

Let $q_1=p,$ $q_2$ and $q_3$ be the $\phi ^N$-periodic points in $\Omega $
satisfying: 
$$
\begin{array}{c}
sequence(q_1)=\ldots 000000000000\ldots  \\ 
sequence(q_2)=\ldots 001001001001\ldots  \\ 
sequence(q_3)=\ldots 011011011011\ldots 
\end{array}
$$

Since there are no singular points of the stable and unstable foliations
inside $R,$ the points $q_2,$ $q_3$ are regular points of the stable and
unstable foliations and are $\phi ^{3N}$-periodic. Moreover, there is a
local product structure inside $R:$ given $i,j\in \{1,2,3\},$ $i\neq j,$ the
intersection $W_{loc}^u(q_i)\cap W_{loc}^s(q_j)$ is transversal and consists
of exactly one point. By $W_{loc}^{s,u}(q_i),$ we mean the connected
components of  $W^{s,u}(q_i)\cap R$ containing $q_i.$

\begin{figure}[!h]
	\centering
	\includegraphics[scale=0.6]{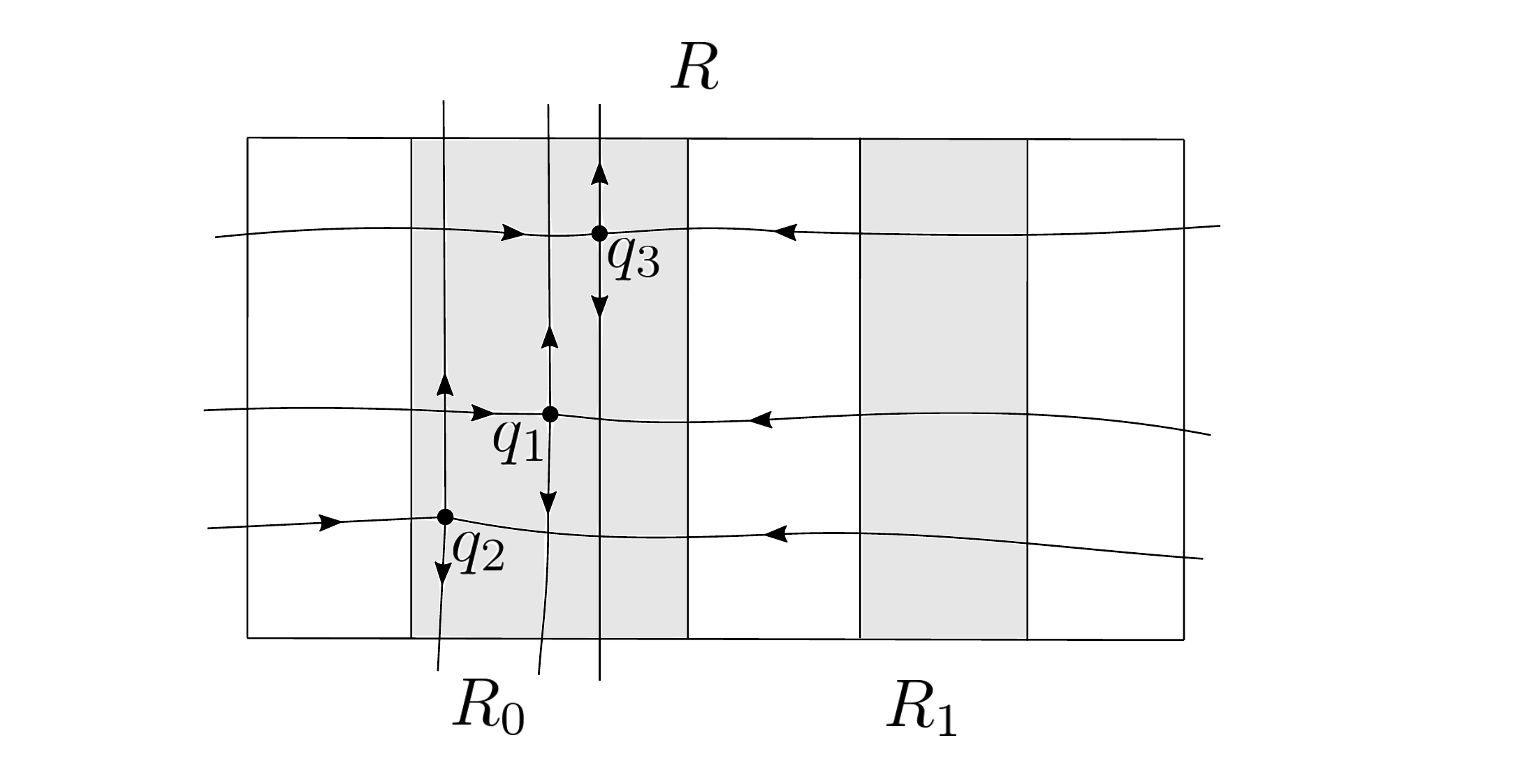}
	\caption{The product structure inside $R$.}
\end{figure}

Consider now $\widetilde R$ the connected component of $\pi ^{-1}(R)$ such
that $\widetilde p=\widetilde q_1\in \widetilde R$ and fix $\widetilde
q_2=\pi ^{-1}(q_2)\cap \widetilde R$ and $\widetilde q_3=\pi ^{-1}(q_3)\cap
\widetilde R$. By construction

$$
\begin{array}{cl}
\widetilde \phi ^{3N}(\widetilde q_1)= & \widetilde q_1, \\ 
\widetilde \phi ^{3N}(\widetilde q_2)= & g(\widetilde q_2), \\ 
\widetilde \phi ^{3N}(\widetilde q_3)= & g^2(\widetilde q_3). 
\end{array}
$$

If we set $\widetilde \psi =g^{-1}\widetilde \phi ^{3N}$, then

$$
\begin{array}{cl}
\widetilde \psi (\widetilde q_1)= & g^{-1}(\widetilde q_1), \\ 
\widetilde \psi (\widetilde q_2)= & \widetilde q_2, \\ 
\widetilde \psi (\widetilde q_3)= & g(\widetilde q_3). 
\end{array}
$$

In particular, $\widetilde{q}_2$ is a hyperbolic fixed saddle point for $
\widetilde{\psi }$. As $W_{loc}^u(q_2)\pitchfork W_{loc}^s(q_1),$ we get
that $W^u(\widetilde{q}_2)\pitchfork W^s(\widetilde{q}_1)$ (note that for
all $m>0,$ $W^{s,u}(g^{-m}(\widetilde{q}_1))$ is the lift of $W^{s,u}(q_1)$
to $g^{-m}(\widetilde{R})).$ Since $\widetilde{\psi }(\widetilde{q}%
_1)=g^{-1}(\widetilde{q}_1),$ using that $W^u(\widetilde{q}_2)$ is invariant
under $\widetilde{\psi }$ and the fact that $g^{-1}$ commutes with $
\widetilde{\psi },$ we conclude that for all $m>0,$

\begin{equation}
\label{int1}W^u(\widetilde q_2)\pitchfork W^s(g^{-m}(\widetilde q_1)). 
\end{equation}

Note that, as $W_{loc}^u(q_1)$ intersects $W_{loc}^s(q_2)$ transversely,
there exists $m^{\prime }>0$ such that

\begin{equation}
\label{int2}W^u(\widetilde q_2)\pitchfork W^s(g^{-m^{\prime }}(\widetilde
q_2)). 
\end{equation}

The same argument considering the point $q_3$ instead of $q_1$, gives an
integer $m^{\prime \prime }>0$ such that

\begin{equation}
\label{int3}W^u(\widetilde q_2)\pitchfork W^s(g^{m^{\prime \prime
}}(\widetilde q_2)). 
\end{equation}

So, by the $\lambda $-lemma, 
$$
W^u(\widetilde{q}_2)\pitchfork W^s(g^{m^{\prime}m^{\prime \prime }}(\widetilde{q}_2))\text{
and }W^u(\widetilde{q}_2)\pitchfork W^s(g^{-m^{\prime}m^{\prime \prime }}(\widetilde{q}_2)), 
$$
which finally imply that 
$$
W^u(\widetilde{q}_2)\pitchfork W^s(\widetilde{q}_2). 
$$

Associated with the transversal intersection $W^u(\widetilde q_2)\pitchfork %
W^s(\widetilde q_2),$ there is a compact $\widetilde \psi ^{N^{\prime }}$%
-invariant set $\Omega _g$, for some $N^{\prime }>0$, such that $%
h(\widetilde \psi ^{N^{\prime }}|_{\Omega _g})>0$. Defining $\Omega
_g^{*}=\cup _{i=0}^{N^{\prime }-1}\widetilde \psi ^i(\Omega _g)$, it is a $%
\widetilde \psi $-invariant compact set with $h(\widetilde \psi |_{\Omega
_g^{*}})>0.$ We are looking for a similar statement for the map $%
g^{-1}\widetilde f^{3N}$.

As in the proof of theorem 1, Handel's global shadowing implies that there
exists a compact $f$-invariant set $W$ and a continuous surjective map $%
s:W\to S$ homotopic to the inclusion such that $s\circ f|_W=\phi \circ s$.
Instead of $W,$ we will consider a compact $f$-invariant subset $W^{\prime
\prime }\subseteq W$ constructed in the following way: Since $\phi $ is
pseudo-Anosov relative to a finite set, there exists a point $z_0$ in $S$
such that $Orb_\phi ^{+}(z_0)=\{\phi ^n(z_0)|n\geq 0\}$ is dense in $S$.
Choose some point $w_0\in s^{-1}(z_0)$ and let $W^{\prime }=\overline{%
Orb_f^{+}(w_0)}$. Clearly, $f(W^{\prime })\subseteq W^{\prime }$ and
defining $W^{\prime \prime }=\cap _{n\geq 0}f^n(W^{\prime }),$ we get that $%
f(W^{\prime \prime })=W^{\prime \prime },$ it is compact and $s(W^{\prime
\prime })=S.$ In particular, $Orb_f^{+}(w_0)$ is dense in $W^{\prime \prime
}.$

Lifting $s:W^{\prime \prime }\rightarrow S$ to $\widetilde{s}:\pi
^{-1}(W^{\prime \prime })\rightarrow \mathbb{D},$ we obtain a compact $%
g^{-1}\widetilde f^{3N}$-invariant set 
$$
K_g\subseteq \cup _{n\geq 0}(g^{-1}\widetilde f^{3N})^n\left( \widetilde
s^{-1}(\Omega _g^{*})\right) \subset \pi ^{-1}(W^{\prime \prime })\subset 
\mathbb{D},\text{ with }h(g^{-1}\widetilde f^{3N}|_{K_g})>0. 
$$
As we explained in subsection 2.9, the fact that $h(g^{-1}\widetilde
f^{3N}|_{K_g})>0,$ implies the existence of a non-atomic, hyperbolic,
ergodic $g^{-1}\widetilde f^{3N}$-invariant Borel probability measure $\mu
_g $ with positive entropy, whose support is contained in $K_g$.

As $\mu _g(R^{\pm }\cap supp(\mu _g$$))=1,$ for $0<\delta <1,$ choosing any
point $\widetilde{r}\in \Lambda _\delta \cap supp(\mu _g)\cap R^{\pm },$ we
get that $\mu _g(V_{\widetilde{r}}\cap \Lambda _\delta \cap R^{\pm })>0.$ So
lemmas \ref{andre1}, \ref{andre2} and theorem \ref{katok} assure that, given 
$\eta >0,$ there exists an inaccessible point $\widetilde{z}_g\in supp(\mu
_g)\subset K_g$ (see definition (\ref{inacess})) such that arbitrarily small
rectangles enclosing $\widetilde{z}_g$ can be obtained, the sides of these
rectangles contained in the invariant manifolds of two hyperbolic $g^{-1}
\widetilde{f}^{3N}$-periodic saddle points, $\widetilde{r}_g^{\prime },$ $
\widetilde{r}_g^{\prime \prime },$ whose orbits are contained in the $\eta $%
-neighborhood of $supp(\mu _g)$. Moreover $W^u(\widetilde{r}_g^{\prime })%
\pitchfork W^s(\widetilde{r}_g^{\prime \prime })$ and $W^u(\widetilde{r}%
_g^{\prime \prime })\pitchfork W^s(\widetilde{r}_g^{\prime })$ in a $C^1$%
-transverse way. So, $W^u(\widetilde{r}_g^{\prime })\pitchfork W^s(
\widetilde{r}_g^{\prime })$ and $W^u(\widetilde{r}_g^{\prime \prime })%
\pitchfork W^s(\widetilde{r}_g^{\prime \prime })$ also in a $C^1$-transverse
way. An important observation that will be used later is that each of these
rectangles contains infinitely many points belonging to $supp(\mu _g)\subset
\pi ^{-1}(W^{\prime \prime }),$ because $\mu _g$ is non-atomic. For these
two periodic points, there exists $k^{\prime },k^{\prime \prime }>0$ such
that 
$$
\widetilde{f}^{3Nk^{\prime }}(\widetilde{r}_g^{\prime })=g^{k^{\prime }}(
\widetilde{r}_g^{\prime })\text{ }\ \text{and }\widetilde{f}^{3Nk^{\prime
\prime }}(\widetilde{r}_g^{\prime \prime })=g^{k^{\prime \prime }}(
\widetilde{r}_g^{\prime \prime }). 
$$

Going back to the surface $S,$ defining $z_g=\pi (\widetilde z_g),$ $%
r_g^{\prime }=\pi (\widetilde r_g^{\prime })$ and $r_g^{\prime \prime }=\pi
(\widetilde r_g^{\prime \prime }),$ then $r_g^{\prime }$ and $r_g^{\prime
\prime }$ are hyperbolic $f$-periodic saddles for which $W^u(r_g^{\prime })%
\pitchfork W^s(r_g^{\prime \prime })$ and $W^u(r_g^{\prime \prime })%
\pitchfork W^s(r_g^{\prime })$ in a $C^1$-transverse way and so, $\overline{%
W^s(r_g^{\prime })}=\overline{W^s(r_g^{\prime \prime })}$ and $\overline{%
W^u(r_g^{\prime })}=\overline{W^u(r_g^{\prime \prime })}.$ Associated with
these points there are small rectangles in $S$ whose sides are contained in
their invariant manifolds and enclosing the point $z_g;$ they are the
projection under $\pi $ of the rectangles in $\mathbb{D}.$

Choose now $g_1,g_2\in Deck(\pi )$ such that they correspond to different
geodesics in $S$. In this way, $g_1g_2\neq g_2g_1$, and their powers are
never conjugated, i.e. for all $h\in Deck(\pi),$ and $n,m$ integers, $%
hg^n_1h^{-1}\neq g^m_2.$

For the maps $Id,g_1,g_2$ we consider the compact sets $K_{Id},$ $K_{g_1}$
and $K_{g_2}$ contained in $\mathbb{D}$ and inaccessible points $\widetilde
z_{Id}\in K_{Id},$ $\widetilde z_{g_1}\in K_{g_1}$ and $\widetilde
z_{g_2}\in K_{g_2}.$ From what we just did, there are hyperbolic $f$%
-periodic saddles $r_{Id}^{\prime },$ $r_{Id}^{\prime \prime },$ $%
r_{g_1}^{\prime },$ $r_{g_1}^{\prime \prime },$ $r_{g_2}^{\prime }$ and $%
r_{g_2}^{\prime \prime }$ in $S,$ with $R_0$ a small rectangle in $S$ whose
sides are contained in the invariant manifolds of $r_{Id}^{\prime }$ and $%
r_{Id}^{\prime \prime }$ and enclosing the point $\pi (\widetilde
z_{Id})=z_{Id}\in W^{\prime \prime }.$ Similarly, for $i\in \{1,2\},$ $R_i$
is a small rectangle in $S$ whose sides are contained in the invariant
manifolds of $r_{g_i}^{\prime }$ and $r_{g_i}^{\prime \prime }$ and
enclosing the point $\pi (\widetilde z_{g_i})=z_{g_i}\in W^{\prime \prime }$.

\begin{figure}[!h]
	\centering
	\includegraphics[scale=0.5]{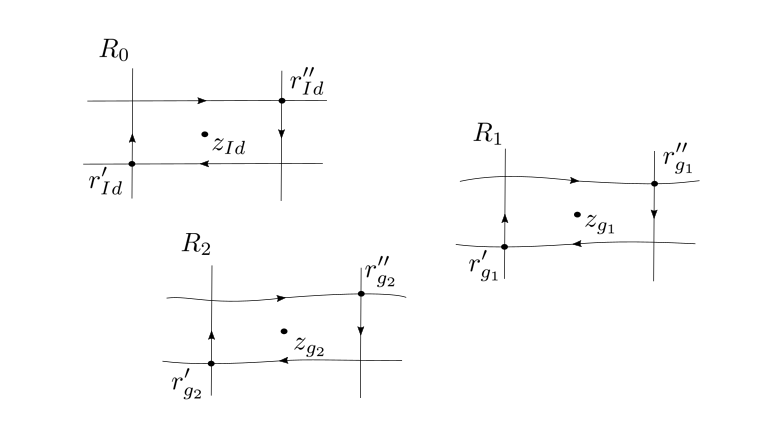}
	\caption{The rectangles $R_0$, $R_1$ and $R_2$.}
	\label{horseshoe}
\end{figure}

Let $n^{comm}>0$ be a natural number that is a common period of all the
points $r_{Id}^{\prime },$ $r_{Id}^{\prime \prime },$ $r_{g_1}^{\prime },$ $%
r_{g_1}^{\prime \prime },$ $r_{g_2}^{\prime }$ and $r_{g_2}^{\prime \prime
}, $ which also leaves invariant all stable and unstable branches of these
points. Clearly, the orbits of all the previous points can be assumed to be
disjoint.

As we said, $R_0$ is the small rectangle enclosing the point $z_{Id}.$ For
all $0\leq i\leq n^{comm}-1$, $f^i(R_0)$ is a rectangle in $S.$ If we denote
the segments in the boundary of $R_0$ as $\alpha _0^{\prime }\in
W^s(r_{Id}^{\prime })$, $\omega _0^{\prime }\in W^u(r_{Id}^{\prime })$, $%
\alpha _0^{\prime \prime }\in W^s(r_{Id}^{\prime \prime })$ and $\omega
_0^{\prime \prime }\in W^u(r_{Id}^{\prime \prime }),$ then for large $m>0$
and for all $0\leq i\leq n^{comm}-1$, 
$$
\partial (f^{n^{comm}m}(f^i(R_0)))\subset f^{n^{comm}m}(f^i(\alpha
_0^{\prime }))\cup W^u(r_{Id}^{\prime })\cup f^{n^{comm}m}(f^i(\alpha
_0^{\prime \prime }))\cup W^u(r_{Id}^{\prime }), 
$$
and the sets $f^{n^{comm}m}(f^i(\alpha _0^{\prime })),$ $f^{n^{comm}m}(f^i(%
\alpha _0^{\prime \prime }))$ are as close as we want to the points $%
f^i(r_{Id}^{\prime })$ and $f^i(r_{Id}^{\prime \prime })$ respectively.
Using an analogous notation with the rectangles $R_1$ and $R_2,$ we can find
a natural number $m_0>0$ such that for $0\leq i,j\leq n^{comm}-1,$ $k,t\in
\{0,1,2\},$ $k\neq t$ and $m>m_0,$

\begin{equation}
\label{int} 
\begin{array}{c}
f^{n^{comm}m}(f^i(\alpha _k^{\prime }))\cap f^j(\omega _t^{\prime
})=\emptyset , \\  
\\ 
f^{n^{comm}m}(f^i(\alpha _k^{\prime }))\cap f^j(\omega _t^{\prime \prime
})=\emptyset , \\  
\\ 
f^{n^{comm}m}(f^i(\alpha _k^{\prime \prime }))\cap f^j(\omega _t^{\prime
})=\emptyset , \\  
\\ 
f^{n^{comm}m}(f^i(\alpha _k^{\prime \prime }))\cap f^j(\omega _t^{\prime
\prime })=\emptyset . 
\end{array}
\end{equation}

As $f^i(z_{Id})$ and $f^i(z_{g_1})$ are in the interior of $f^i(R_0)$ and $%
f^i(R_1),$ respectively, they are both accumulated by points in $W^{\prime
\prime }$ and as there exists a point whose orbit is dense in $W^{\prime
\prime },$ we get that for all $0\leq i\leq n^{comm}-1$, there exist
integers $l_0(i),l_1(i)>m_0n^{comm}$ such that 
$$
f^{l_1(i)}(f^i(R_1))\cap R_0\neq \emptyset \text{ and }f^{l_0(i)}(f^i(R_0))%
\cap R_1\neq \emptyset . 
$$

So, for any $0\leq i\leq n^{comm}-1$, there exists integers $%
m_0(i),m_1(i)\geq m_0$ and other integers $0\leq j_0(i),j_1(i)\leq
n^{comm}-1 $ such that 
$$
f^{n^{comm}m_1(i)}(f^i(R_1))\cap f^{j_0(i)}(R_0)\neq \emptyset \text{ and }%
f^{n^{comm}m_0(i)}(f^i(R_0))\cap f^{j_1(i)}R_1\neq \emptyset . 
$$

Since for the rectangle $R_0,$ its boundary is contained in the invariant
manifolds of $r_{Id}^{\prime },$ $r_{Id}^{\prime \prime },$ and 
$$
\begin{array}{c}
W^u(r_{Id}^{\prime }) 
\pitchfork W^s(r_{Id}^{\prime \prime }); \\ W^u(r_{Id}^{\prime \prime })%
\pitchfork W^s(r_{Id}^{\prime }); 
\end{array}
$$
the same being true for $R_1,$ $r_{g_1}^{\prime }$ and $r_{g_1}^{\prime
\prime },$ we conclude by expression (\ref{int}) that for all $0\leq i\leq
n^{comm}-1$, $W^u(f^i(r_{Id}^{\prime }))\pitchfork W^s(f^{j_1(i)}(r_{g_1}^{%
\prime }))$ and $W^u(f^i(r_{g_1}^{\prime }))\pitchfork %
W^s(f^{j_0(i)}(r_{Id}^{\prime }))$. Then a combinatorial argument implies
that there exist $0\leq i,j\leq n^{comm}-1$ such that $W^u(f^i(r_{g_1}^{%
\prime }))\pitchfork W^s(f^j(r_{Id}^{\prime }))$ and $W^u(f^j(r_{Id}^{\prime
}))\pitchfork W^s(f^i(r_{g_1}^{\prime }))$, see \cite{saz15}.

Doing the same for $R_2$, $r_{g_2}^{\prime }$ and $r_{g_2}^{\prime \prime }$%
, and using that topologically transverse intersections are mapped into
themselves under $f,$ we can find $0\leq k\leq n^{comm}-1$ such that for the
same $j$ above, $W^u(f^k(r_{g_2}^{\prime }))\pitchfork W^s(f^j(r_{Id}^{%
\prime }))$ and $W^u(f^j(r_{Id}^{\prime }))\pitchfork W^s(f^k(r_{g_2}^{%
\prime }))$.

Set $r_0=f^j(r_{Id}^{\prime })$, $r_1=f^i(r_{g_1}^{\prime })$ and $%
r_2=f^k(r_{g_2}^{\prime })$. Then, these are hyperbolic $f$-periodic saddle
points and the following intersections hold, for $i\in \{1,2\}:$ 
$$
W^u(r_0)\pitchfork W^s(r_i)\text{ and }W^u(r_i)\pitchfork W^s(r_0). 
$$

Fix any $\widetilde r_0$ in $\pi ^{-1}(r_0)$. By our construction, since
there is a point $\widetilde r_0^{\prime }\in \pi ^{-1}(r_0)$ whose orbit is
forever close to $K_{Id},$ we get that $\widetilde f^{n^{comm}}(\widetilde
r_0)=\widetilde r_0.$

Recall that $n^{comm}$ is a common period for $r_0$, $r_1$ and $r_2$. The
fact that $W^u(r_0)\pitchfork W^s(r_1)$ implies that there exists a point $%
\widetilde r_1\in \pi ^{-1}(r_1)$ for which $W^u(\widetilde r_0)\pitchfork %
W^s(\widetilde r_1)$. Moreover, arguing as above, there exists an integer $%
n_1>0$ such that $\widetilde f^{n^{comm}}(\widetilde r_1^{\prime
})=g_1^{n_1}(\widetilde r_1^{\prime })$ for some $\widetilde r_1^{\prime
}\in \pi ^{-1}(r_1),$ close to $K_{g_1}.$ As $\widetilde r_1,\widetilde
r_1^{\prime }\in \pi ^{-1}(\widetilde r_1),$ there exists $h_1\in Deck(\pi )$
with $\widetilde r_1=h_1(\widetilde r_1^{\prime }).$ Hence, $\widetilde
f^{n^{comm}}(\widetilde r_1)=h_1g_1^{n_1}h_1^{-1}(\widetilde r_1)$.

Set $g_1^{\prime }=h_1g_1^{n_1}h_1^{-1}$. As we did before, since $%
W^u(\widetilde r_0)\pitchfork W^s(\widetilde r_1)$ and $\widetilde
f^{n^{comm}}(\widetilde r_1)=g_1^{\prime }(\widetilde r_1),$ for all $m\geq
0 $ 
$$
W^u(\widetilde r_0)\pitchfork W^s((g_1^{\prime })^m(\widetilde r_1)). 
$$

As $W^u(r_1)$ intersects $W^s(r_0)$ in a topologically transverse way, there
is a compact connected piece of a branch of $W^u(r_1),$ denoted $\lambda _1,$
such that one of its endpoints is $r_1$ and it has a topologically
transversal intersection with $W^s(r_0).$ If $\widetilde{\lambda }_1$ is the
lift of $\lambda _1$ starting at the point $\widetilde{r}_1,$ then there
exists $h_1^{\prime }\in Deck(\pi )$ such that 
$$
\widetilde{\lambda }_1\pitchfork W^s(h_1^{\prime }(\widetilde{r}_0)). 
$$

This implies that if $m_1>0$ is sufficiently large, then a piece of $%
W^u(\widetilde r_0)$ is sufficiently close in the Hausdorff topology to $%
(g_1^{\prime })^{m_1}(\widetilde \lambda _1)$, something that forces $%
W^u(\widetilde r_0)$ to have a topological transverse intersection with $%
W^s((g_1^{\prime })^{m_1}h_1^{\prime }(\widetilde r_0))$. In other words,
for all $m_1>0$ sufficiently large,%
$$
W^u(\widetilde r_0)\pitchfork W^s((g_1^{\prime })^{m_1}h_1^{\prime
}(\widetilde r_0)). 
$$

Arguing in an analogous way with respect to the point $r_2$ we find $%
h_2,h_2^{\prime }\in Deck(\pi )$ and an integer $n_2>0,$ 
such that, if $g_2^{\prime }=h_2g_2^{n_2}h_2^{-1}$, then for all $m_2>0$
sufficiently large,%
$$
W^u(\widetilde r_0)\pitchfork W^s((g_2^{\prime })^{m_2}h_2^{\prime
}(\widetilde r_0)). 
$$

In order to conclude, let us show that $m_1,m_2>0$ can be chosen in a way
that $(g_1^{\prime })^{m_1}h_1^{\prime }$ and $(g_2^{\prime
})^{m_2}h_2^{\prime }$ do not commute. We started with deck transformations $%
g_1$ and $g_2$ for which $g_1g_2\neq g_2g_1$ and $g^n_1$ is not conjugated
to $g^m_2$, for all integers $n,m.$ As we already explained, the above
conditions follow from the fact that $g_1$ and $g_2$ correspond, in $S,$ to
different geodesics.

In particular this implies that the deck transformations $g_1^{\prime }$ and 
$g_2^{\prime }$ do not commute and the fixed points of $g_1^{\prime }$ and $%
g_2^{\prime }$ at the boundary at infinity $\partial \mathbb{D}$ are all
different, i.e $Fix(g_1^{\prime })\cap Fix(g_2^{\prime })=\emptyset $.

Fix two large integers $m_1,m_2>0$ and let us analyze $(g_1^{\prime
})^{m_1}h_1^{\prime }$ and $(g_2^{\prime })^{m_2}h_2^{\prime }$. If they do
not commute, there is nothing to do.

So assume that $(g_1^{\prime })^{m_1}h_1^{\prime }$ and $(g_2^{\prime
})^{m_2}h_2^{\prime }$ commute. Since they commute, $Fix((g_1^{\prime
})^{m_1}h_1^{\prime })=Fix((g_2^{\prime })^{m_2}h_2^{\prime })$. Observe
that either $g_1^{\prime }$ does not commute with $(g_1^{\prime
})^{m_1}h_1^{\prime }$ or $g_2^{\prime }$ does not commute with $%
(g_2^{\prime })^{m_2}h_2^{\prime }$.

In fact, if they both commute, then 
$$
Fix(g_1^{\prime })=Fix((g_1^{\prime })^{m_1}h_1^{\prime })=Fix((g_2^{\prime
})^{m_2}h_2^{\prime })=Fix(g_2^{\prime }), 
$$
and this contradicts the fact that $g_1^{\prime }$ and $g_2^{\prime }$ do
not commute. So, without loss of generality , assume that $g_1^{\prime }$
and $(g_1^{\prime })^{m_1}h_1^{\prime }$ do not commute. Hence $%
Fix(g_1^{\prime })\cap Fix((g_1^{\prime })^{m_1}h_1^{\prime })=\emptyset $.

We claim that $(g_1^{\prime })^{m_1+1}h_1^{\prime }=g_1^{\prime
}(g_1^{\prime })^{m_1}h_1^{\prime }$ and $(g_2^{\prime })^{m_2}h_2^{\prime }$
do not commute. Otherwise, 
$$
Fix((g_1^{\prime })^{m_1+1}h_1^{\prime })=Fix((g_2^{\prime
})^{m_2}h_2^{\prime })=Fix((g_1^{\prime })^{m_1}h_1^{\prime }). 
$$

So, for all $\widetilde q\in Fix((g_1^{\prime })^{m_1}h_1^{\prime }),$ 
$$
\widetilde q=g_1^{\prime }((g_1^{\prime })^{m_1}h_1^{\prime }(\widetilde
q))=g_1^{\prime }(\widetilde q), 
$$
which means that $Fix((g_1^{\prime })^{m_1}h_1^{\prime })=Fix(g_1^{\prime
}), $ that is a contradiction with our previous assumption that $g_1^{\prime
}$ and $(g_1^{\prime })^{m_1}h_1$ do not commute. Hence $(g_1^{\prime
})^{m_1+1}h_1^{\prime }$ and $(g_2^{\prime })^{m_2}h_2^{\prime }$ do not
commute.

So, we can always can find arbitrarily large integers $m_1,m_2>0$ such that
if $g_1^{\prime \prime }=(g_1^{\prime })^{m_1}h_1^{\prime }$ and $%
g_2^{\prime \prime }=(g_2^{\prime })^{m_2}h_2^{\prime }$, then $g_1^{\prime
\prime }g_2^{\prime \prime }\neq g_2^{\prime \prime }g_1^{\prime \prime }$
and 
$$
W^u(\widetilde r_0)\pitchfork W^s(g_1^{\prime \prime }(\widetilde r_0))\text{
and }W^u(\widetilde r_0)\pitchfork W^s(g_2^{\prime \prime }(\widetilde
r_0)). \text{ } 
$$

Now, as in the proof of lemma \ref{lema2}, we construct the path connected
sets $\theta $ and $\theta ^{\prime }$ using the point $\widetilde r_0$ and
the deck transformations $g_1^{\prime \prime }$ and $g_2^{\prime \prime }$.
Since $f$ has a fully essential system of curves $\mathscr{C}$ and the
periodic points $P$ associated to $\mathscr{C},$ the exact same proof of the
lemma \ref{lema2} without any modifications shows that for every $g\in
Deck(\pi )$, 
$$
W^u(\widetilde{r}_0)\pitchfork W^s(g(\widetilde{r}_0)). 
$$
As $\widetilde{r}_0\in \pi ^{-1}(r_0)$ was arbitrary, redefining $p=r_0,$
the proof is over.

\qed

\section{Proof of theorem 3}

Let $\widetilde{p}\in \mathbb{D}$ be a hyperbolic periodic saddle point for $
\widetilde{f}$ as in theorem 2 (as before, assume without loss of generality
that $\widetilde{p}$ is fixed and all four branches at $\widetilde{p}$ are $
\widetilde{f}$-invariant, otherwise consider some iterate of $\widetilde{f}).
$ For all $g\in Deck(\pi ),$ 
$$
W^u(\widetilde{p})\pitchfork W^s(g(\widetilde{p})). 
$$

In fact, a stronger statement holds: the proof of theorem 2 gives an
unstable branch $\widetilde \lambda _u$ of $W^u(\widetilde p)$ and a stable
branch $\widetilde{\beta }_s$ of $W^s(\widetilde p)$ such that for all $g\in
Deck(\pi ),$

\begin{equation}
\label{branch}\widetilde{\lambda }_u\pitchfork g(\widetilde{\beta }_s).
\end{equation}

Fix some $0<\epsilon <1/10$ small enough so that for any $z\in S,$ if $
\widetilde{z}_1,\widetilde{z}_2\in \pi ^{-1}(z),$ $\widetilde{z}_1\neq 
\widetilde{z}_2,$ then $B_{2\epsilon }(\widetilde{z}_1)\cap B_{2\epsilon }( 
\widetilde{z}_2)=\emptyset .$ Let $\widetilde \lambda $ be a compact subarc
of $\widetilde \lambda _u$ small enough so that one of its endpoints is $%
\widetilde p$ and $\widetilde \lambda \subset B_\epsilon (\widetilde q)$. In
a similar way, let $\widetilde \beta $ be a compact subarc of $\widetilde{%
\beta }_s,$ so that $\widetilde p$ is one of its endpoints and $\widetilde
\beta \subset B_\epsilon (\widetilde p).$ The arc $\widetilde{\beta }$
satisfies another property: Its endpoint which is not $\widetilde{p}$
belongs to $W^u(\widetilde p)$ and actually, this end point is a $C^1$%
-transversal homoclinic point. It is possible to choose $\widetilde{\beta }$
in this way because the proof of theorem 2 implies the existence of a $C^1$%
-transversal intersection between $W^s(Id(\widetilde p))$ and $%
W^u(\widetilde p).$ When instead of $Id,$ we consider any other deck
transformation, only topologically transverse intersections are assured, but
for the $Id,$ $C^1$-transversality was obtained.

Now choose $h_1,h_2,\ldots ,h_{2g}\in Deck(\pi ),$ where $g>0$ is the genus
of $S,$ such that the complement in $S$ of the union of the geodesics in $S$
associated to $\{h_1,h_2,\ldots ,h_{2g}\}$ is a union of open disks.
Expression (\ref{branch}) implies the existence of a compact arc $\widetilde{%
\Lambda }$ such that $\widetilde \lambda _u\supset \widetilde{\Lambda }%
\supset \widetilde \lambda $ and 
\begin{equation}
\label{defineM} 
\begin{array}{c}
\widetilde{\Lambda }\text{ contains both endpoints of }\widetilde \beta , 
\text{ } \widetilde{\Lambda }\pitchfork h_i(\widetilde \beta ),\forall 1\leq
i\leq 2g, \text{ the } \\ \text{endpoint of }\widetilde{\Lambda }\text{
which is not } \widetilde{p}\text{ is contained in the interior of }%
\widetilde \beta \text{ and } \\ \text{it is a }C^1\text{-transversal
homoclinic point. } 
\end{array}
\end{equation}
Clearly the above choice implies that every connected component of the
complement of $\pi (\widetilde{\Lambda }\cup \widetilde \beta )$ is an open
disk in $S.$

Let $R\subset B_{2\epsilon }(p)$ be a closed rectangle which has $p$ as a
vertex and $\beta =\pi ($$\widetilde \beta )$ as one side; $R$ is very thin,
close to $\beta $ in the Hausdorff topology. $\partial R$ is given by the
union of 4 arcs: $\alpha ,\alpha ^{\prime },\beta $ and $\beta ^{\prime }.$
The arcs $\alpha $ and $\alpha ^{\prime }$ are contained in $W^u(p):$ $%
\alpha \subset \pi (\widetilde{\lambda })$ and $\alpha ^{\prime }$ contains
the endpoint of $\beta $ which is not $p.$

From the choice of $\widetilde{\Lambda },$ 
\begin{equation}
\label{naoprecisatirar}\pi (\widetilde{\Lambda })\supset \alpha ^{\prime }. 
\end{equation}
Clearly, $\beta $ and $\beta ^{\prime }$ are contained in $W^s(p)$ and $%
\beta ^{\prime }$ is $C^1$-close to $\beta .$ As was explained when defining 
$\widetilde{\beta },$ the existence of such a rectangle $R$ follows from
theorem 2, which says that $W^u(\widetilde p)$ has $C^1$-transverse
intersections with $W^s(\widetilde p).$

At this point we need to determine the size of $\alpha $ and $\alpha
^{\prime }$ and a number $N>0,$ as follows: we know (from (\ref{defineM}))
that $\widetilde{\Lambda }\pitchfork h_i(\widetilde{\beta }),\forall 1\leq
i\leq 2g.$ Choose $\beta ^{\prime }$ sufficiently close to $\beta $ (so $%
\alpha $ and $\alpha ^{\prime }$ are very small), in a way that if $
\widetilde{R}$ is the connected component of $\pi ^{-1}(R)$ which contains $
\widetilde{\beta }$ (the sides of $\widetilde{R}$ are denoted as $\widetilde{%
\alpha }\subset \widetilde{\lambda },\widetilde{\alpha }^{\prime },
\widetilde{\beta }$ and $\widetilde{\beta }^{\prime };$ $\widetilde{\alpha },
\widetilde{\alpha }^{\prime }\subset W^u(\widetilde{p})$ and $\widetilde{%
\beta },\widetilde{\beta }^{\prime }\subset W^s(\widetilde{p})),$ then

$$
\widetilde{\Lambda }\pitchfork h_i(\widetilde \beta ^{\prime }),\forall
1\leq i\leq 2g. 
$$

Now, fix some $N>0$ such that 
\begin{equation}
\label{defineNunif} 
\begin{array}{c}
\widetilde f^N( 
\widetilde{\beta }^{\prime })\subset \widetilde{\beta },\text{ }\widetilde
f^N( \widetilde{\alpha })\supset \widetilde{\Lambda }\supset \widetilde{%
\alpha }^{\prime }\text{ and }\widetilde f^N(\widetilde{\alpha }^{\prime
})\supset \widetilde{\Lambda }^{\prime },\text{ an arc } \\ \text{%
sufficiently }C^1 \text{-close to }\widetilde{\Lambda },\text{ whose
endpoints} \\ \text{ are also in }\widetilde{\beta },\text{ in a way that } 
\widetilde{\Lambda }^{\prime }\pitchfork h_i(\widetilde{\beta })\text{ }\ 
\text{and} \\ \widetilde{\Lambda }^{\prime }\pitchfork h_i(\widetilde{\beta }%
^{\prime }),\forall 1\leq i\leq 2g.\text{ Moreover, the arcs in }\widetilde{%
\beta } \\ \text{ connecting the appropriate endpoints (the ones } \\ \text{%
which are closer) of }\widetilde{\Lambda }\text{ and }\widetilde{\Lambda }%
^{\prime }\text{ are disjoint from}\ \widetilde{\Lambda }\cup \ \widetilde{%
\Lambda }^{\prime }. 
\end{array}
\end{equation}

Define

\begin{equation}
\label{os4m}\widetilde{M}_{\widetilde{\Lambda }}=filled\left( \widetilde{%
\beta }\cup \widetilde{\Lambda }\right) ,\text{ }\widetilde{M}_{\widetilde{%
\Lambda }^{\prime }}=filled\left( \widetilde{\beta }\cup \widetilde{\Lambda }%
^{\prime }\right) \text{ and }\widetilde{M}_{\min }=\widetilde{M}_{ 
\widetilde{\Lambda }}\cap \widetilde{M}_{\widetilde{\Lambda }^{\prime }}, 
\end{equation}
where for any compact connected subset $\widetilde{K}$ of $\mathbb{D},$ 
$$
filled(\widetilde{K})=\widetilde{K}\cup \{\text{all bounded connected
components of }\widetilde{K}^c\}. 
$$
It is well known than $fill(\widetilde{K})^c$ is open, connected and
unbounded.

\begin{figure}[!h]
	\centering
	\includegraphics[scale=0.22]{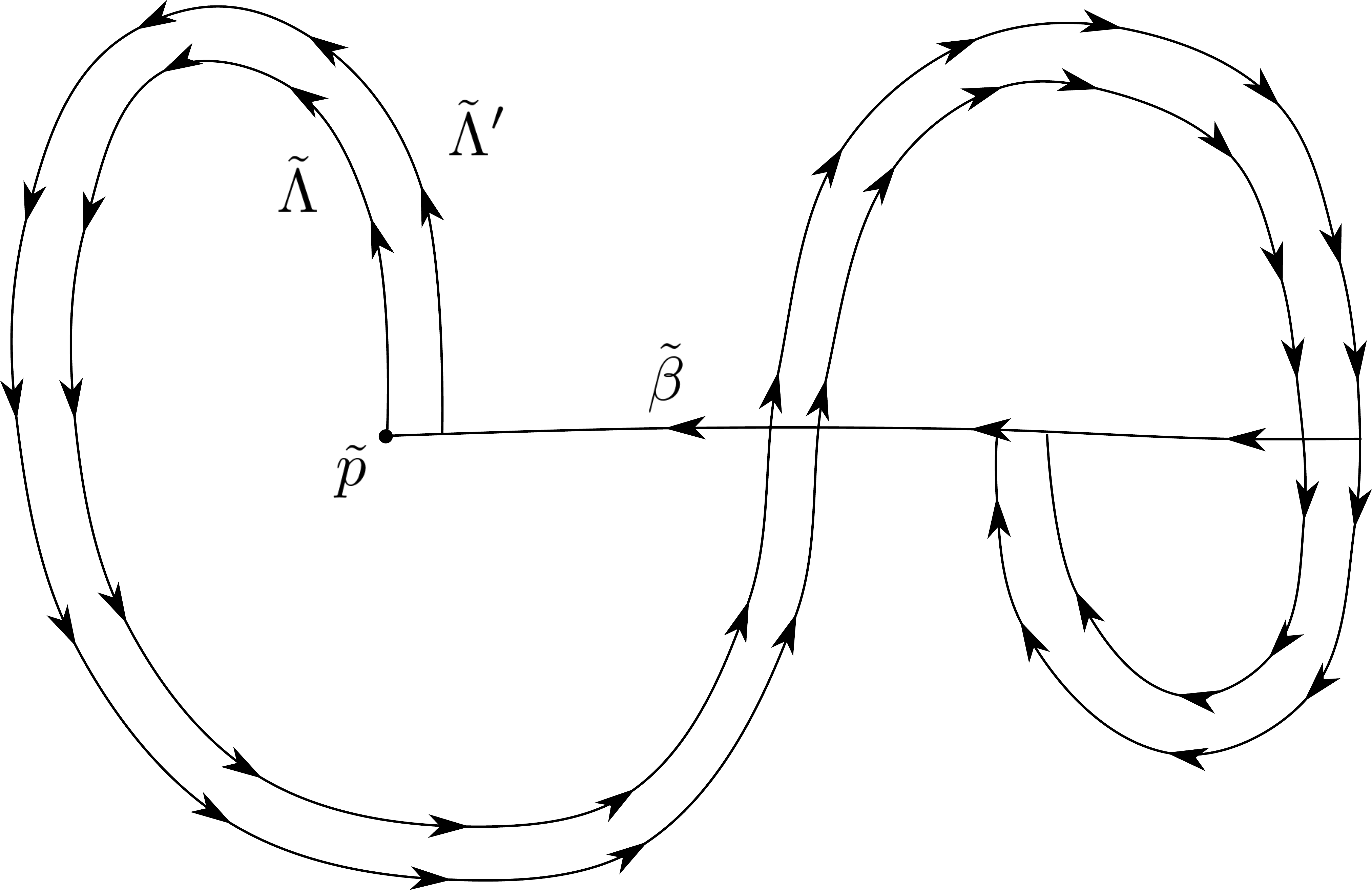}
	\caption{How to obtain the sets $\widetilde{M}_{\widetilde{\Lambda }}$, $\widetilde{M}_{\widetilde{\Lambda }^{\prime }}$ and $\widetilde{M}_{\min }$.}
\end{figure}

From the choice of $\widetilde{\Lambda },\widetilde{\Lambda }^{\prime }$ and 
$\widetilde{\beta },$ the sets $\widetilde{M}_{\widetilde{\Lambda }}, 
\widetilde{M}_{\widetilde{\Lambda }^{\prime }},\widetilde{M}_{\min }$ are
connected and the complement of any of the three sets $\pi (\widetilde{M}_{ 
\widetilde{\Lambda }}),\pi (\widetilde{M}_{\widetilde{\Lambda }^{\prime
}}),\pi (\widetilde{M}_{\min })$ is a union of open disks. So, given a
fundamental domain $\widetilde{Q}\subset \mathbb{D}$ of $S,$ there exist
deck transformations $\left\{ m_1,m_2,...,m_J\right\} ,$ for some $J>0,$
such that

$$
\cup _{i=1}^Jm_i(\widetilde{M}_{\min }) 
$$
is a (bounded) connected closed set, its complement has a bounded connected
component denoted $\widetilde{\theta },$ which contains $\widetilde{Q}.$
Moreover, 
$$
d_{\mathbb{D}}\left( \cup _{i=1}^Jm_i(\widetilde{M}_{\min }),\widetilde{Q}%
\right) >1. 
$$
In particular, this implies that $\widetilde{Q}\subset filled(\cup
_{i=1}^Jm_i(\widetilde{M}_{\min })).$

The reason for the above to be true is the following: $\pi ^{-1}\left( \pi (
\widetilde{M}_{\min })\right) $ is a closed connected equivariant subset of $%
\mathbb{D}$ and its complement has only open topological disks as connected
components, all with uniformly bounded diameters. Let $\widetilde{\Gamma }$
be a simple closed curve which surrounds $\widetilde{Q}$ and such that 
$$
d_{\mathbb{D}}(\widetilde{\Gamma },\widetilde{Q})>1. 
$$
As $\widetilde{\Gamma }$ is compact and $\widetilde{\Lambda }\pitchfork h_i(
\widetilde{\beta }),$ $\widetilde{\Lambda }^{\prime }\pitchfork h_i(
\widetilde{\beta }),$ $\forall 1\leq i\leq 2g,$ there exists deck
transformations $\left\{ m_1,m_2,...,m_J\right\} ,$ for some $J>0,$ such
that $\cup _{i=1}^Jm_i(\widetilde{M}_{\min })$ is connected and its
complement has a bounded connected component (the one we previously denoted $
\widetilde{\theta })$ which contains $\widetilde{\Gamma }.$ Moreover, if $%
\nu $ is a simple arc which avoids unstable manifolds of periodic saddle
points and $\nu $ connects a point in the unbounded connected component of $%
\left( Neighborhood_{1/5}\left( \cup _{i=1}^Jm_i(\widetilde{M}_{\min
})\right) \right) ^c$ to a point in $\widetilde{Q},$ then, for some $i\in
\{1,2,...,J\},$ $\nu $ must cross $m_i(\widetilde{R})$ from $m_i(\widetilde{%
\beta })$ to $m_i(\widetilde{\beta }^{\prime })$ or vice-versa. This happens
because as $diameter(R)<2.\epsilon <1/5,$ 
$$
Neighborhood_{1/5}\left( \cup _{i=1}^Jm_i(\widetilde{M}_{\min })\right)
\supset \cup _{i=1}^Jm_i(\widetilde{R}). 
$$

Remembering subsection 1.2, the following equalities hold: 
$$
\rho _m(f)=Conv(\rho _{erg}(f))=Conv(\rho _{mz}(f)) 
$$

We also know that every extremal point of the convex hull of $\rho _{mz}(f)$
is the rotation vector of some recurrent point.

Let $w$ be a extremal point of $Conv(\rho _{mz}(f))$, and let $q_w\in S$ be
a recurrent point with 
$$
w=\lim _{n\to \pm \infty }\frac{\Psi _f^n(q_w)}n. 
$$
From the existence of a fully essential system of curves $\mathscr{C},$ it
is easy to see that $\underline{0}=(0,...,0)$ belongs to the interior of the 
$Conv(\rho _{erg}(f)).$ So, $w\neq \underline{0}.$

Since $q_w$ is a recurrent point, if we fix $\widetilde q_w\in \pi
^{-1}(q_w)\cap \widetilde{Q},$ there exists a sequence $n_k\to \infty $ as $%
k\to \infty ,$ such that for some $g_k\in Deck(\pi )$ (depending on $k),$ 
$$
\widetilde f^{n_k}(\widetilde q_w)\in g_k(\widetilde{Q})\text{ and }d_{%
\mathbb{D}}(\widetilde f^{n_k}(\widetilde q_w),g_k(\widetilde q_w))<\frac
1k,\forall k>0. 
$$

For all $k>0,$ let $\beta _k$ be a path in $S$ joining $f^{n_k}(q_w)$ to $%
q_w $ with $l(\beta _k)<1/k.$ As, $||[I_{q_w}^{n_k}*\beta _k]-\Psi
_f^{n_k}(q_w)||\leq 2C_{\mathcal{A}}+1$ and $w=\lim _{k\to \infty }\Psi
_f^{n_k}(q_w)/n_k,$ we get that 
$$
w=\lim _{k\to \infty }\frac{[I_{q_w}^{n_k}*\beta _k]}{n_k}. 
$$

Let $\widetilde I_{\widetilde q_w}^{n_k}$ be the lift of $I_{q_w}^{n_k}$
with base point $\widetilde q_w$. Then $\widetilde I_{\widetilde q_w}^{n_k}$
is a path in $\mathbb{D}$ joining $\widetilde q_w$ to $\widetilde
f^{n_k}(\widetilde q_w)$. Since $d_{\mathbb{D}}(\widetilde
f^{n_k}(\widetilde q_w),g_k(\widetilde q_w))<1/k$, the loop $%
I_{q_w}^{n_k}*\beta _k$ lifts to a path $\widetilde I_{\widetilde
q_w}^{n_k}*\widetilde \beta _k$ joining $\widetilde q_w$ to $g_k(\widetilde
q_w)$.

For any $g\in Deck(\pi )$, a path $\widetilde \gamma _g$ joining any point $%
\widetilde q\in \mathbb{D}$ to $g(\widetilde q)$ projects into a loop $%
\gamma _g=\pi (\widetilde \gamma _g)$ whose free homotopy class (and in
particular its homology class) is determined only by $g$. We denote by $%
[g]=[\gamma _g]$ this homology class. Hence, we can write 
$$
w=\lim _{k\to \infty }\frac{[g_k]}{n_k}. 
$$

We want to prove the following fact:

\begin{description}
\item[Fact]  \label{factfact} {\bf \ref{factfact}.} {\it For all
sufficiently large $k,$ and $\widetilde{R}$ the connected component of $\pi
^{-1}(R)$ which contains $\widetilde{\beta },$ there exists $i_0=i_0(k)$ and 
$i_1=i_1(k)$ in $\{1,...,J\}$ such that } 
$$
\widetilde{f}^{N+n_k}(\widetilde{R})\cap \left( m{}_{i_0}^{-1}g_km_{i_1}(
\widetilde{R})\right) \supset \widetilde{R}_1, 
$$
{\it where $N>0$ is given in (\ref{defineNunif}) and $\widetilde{R}_1$ is a
''vertical rectangle'' in $m_{i_0}^{-1}g_km_{i_1}(\widetilde{R}):$ Two of
its sides are contained, one in $m{}_{i_0}^{-1}g_km_{i_1}(\widetilde{\beta })
$ and the other in $m{}_{i_0}^{-1}g_km_{i_1}(\widetilde{\beta }^{\prime })$
and the two other sides are contained in the interior of $%
m{}_{i_0}^{-1}g_kh_{i_1}(\widetilde{R}),$ each one connecting a point from
one of the previous sides to the other. Clearly }

$$
\widetilde{f}^{N+n_k}(\widetilde{R})\cap \widetilde{R}\supset \widetilde{R}%
_0, 
$$
{\it a rectangle similar to $\widetilde{R}_1,$ but contained in $\widetilde{R%
}.$}
\end{description}

\vskip 0.2truecm

{\it Proof: }

Assume $k>0$ is sufficiently large, so that $n_k>2.N$ and 
$$
filled(\cup _{i=1}^Jm_i(\widetilde{M}_{\widetilde{\Lambda }}\cup \widetilde{M%
}_{\widetilde{\Lambda }^{\prime }}))\cap g_k(filled(\cup _{i=1}^Jm_i( 
\widetilde{M}_{\widetilde{\Lambda }}\cup \widetilde{M}_{\widetilde{\Lambda }%
^{\prime }}\cup \widetilde{R})))=\emptyset . 
$$
This implies the following:

\begin{description}
\item[Remark]  \label{remarkremark} {\bf \ref{remarkremark}.} $\widetilde{f}%
^{n_k}(filled(\cup _{i=1}^Jm_i(\widetilde{M}_{\widetilde{\Lambda }}\cup 
\widetilde{M}_{\widetilde{\Lambda }^{\prime }})))$ {\it does not intersect} $%
g_k\left( \cup _{i=1}^Jm_i(\widetilde{\Lambda }\cup \widetilde{\Lambda }%
^{\prime })\right) .$
\end{description}

\vskip 0.2truecm

{\it Proof of the remark:} Otherwise, if some point 
$$
\widetilde{z}\in g_k\left( \cup _{i=1}^Jm_i(\widetilde{\Lambda }\cup 
\widetilde{\Lambda }^{\prime })\right) \cap \widetilde{f}^{n_k}(filled(\cup
_{i=1}^Jm_i(\widetilde{M}_{\widetilde{\Lambda }}\cup \widetilde{M}_{ 
\widetilde{\Lambda }^{\prime }}))), 
$$
then $\widetilde{f}^{-n_k}(\widetilde{z})\in g_k\left( \cup _{i=1}^Jm_i( 
\widetilde{\alpha })\right) \cap filled(\cup _{i=1}^Jm_i(\widetilde{M}_{ 
\widetilde{\Lambda }}\cup \widetilde{M}_{\widetilde{\Lambda }^{\prime }})),$
which is contained in

$$
g_k(filled(\cup _{i=1}^Jm_i(\widetilde{M}_{\widetilde{\Lambda }}\cup 
\widetilde{M}_{\widetilde{\Lambda }^{\prime }}\cup \widetilde{R})))\cap
filled(\cup _{i=1}^Jm_i(\widetilde{M}_{\widetilde{\Lambda }}\cup \widetilde{M%
}_{\widetilde{\Lambda }^{\prime }}))=\emptyset 
$$
and this is a contradiction. \qed

\vskip 0.2truecm

The previous remark, although simple, will be very important.

\noindent As $\widetilde q_w\in \widetilde{Q}\subset filled(\cup
_{i=1}^Jm_i( \widetilde{M}_{\min }))$ and $\widetilde f^{n_k}(\widetilde
q_w)\in g_k(\widetilde{Q}),$ we can argue as follows: Consider the connected
components of 
$$
interior\left\{ filled\left[ \widetilde f^{n_k}\left( \cup _{i=1}^Jm_i( 
\widetilde{M}_{\min })\right) \cup (\cup _{i=1}^Jm_i(\widetilde{M}_{\min
}))\right] \cap \left( filled(\cup _{i=1}^Jm_i(\widetilde{M}_{\min
}))\right) ^c\right\} 
$$
%
From the existence of $\widetilde q_w$ as above, there is one such connected
component, denoted $\widetilde{C}_k,$ which intersects $g_k(\widetilde{Q}).$
The boundary of $\widetilde{C}_k$ is a Jordan curve, made of two simple arcs
which only intersect at their endpoints: one arc is contained in $\partial
\left( filled(\cup _{i=1}^Jm_i(\widetilde{M}_{\min }))\right) $ and its
endpoints are in $\cup _{i=1}^Jm_i(\widetilde{\beta })$ and the other arc is
equal to $\widetilde f^{n_k}(\widetilde{\xi }),$ where $\widetilde{\xi }$ is
an arc either contained in $m_{i_0}(\widetilde{\Lambda })$ or $m_{i_0}( 
\widetilde{\Lambda }^{\prime })$ (assume it is $m_{i_0}(\widetilde{\Lambda }%
)),$ for some $i_0\in \{1,...,J\}.$

As both endpoints of $\widetilde f^{n_k}(\widetilde{\xi })$ are contained in 
$\cup _{i=1}^Jm_i(\widetilde{\beta })\subset filled(\cup _{i=1}^Jm_i( 
\widetilde{M}_{\min })),$ there exists some $i_1\in \{1,...,J^{\prime }\},$
such that $\widetilde{f}^{n_k}(\widetilde{\xi })$ crosses $g_k.m_{i_1}( 
\widetilde{R})$ from outside

\noindent $g_k(filled(\cup _{i=1}^Jm_i(\widetilde{M}_{\widetilde{\Lambda }%
}\cup \widetilde{M}_{\widetilde{\Lambda }^{\prime }}\cup \widetilde{R})))$
to inside $g_k(\widetilde{Q}),$ that is, it crosses $g_k.m_{i_1}(\widetilde{R%
})$ from $g_k.m_{i_1}(\widetilde{\beta })$ to $g_k.m_{i_1}(\widetilde{\beta }%
^{\prime })$ or vice-versa, in order to intersect $g_k(\widetilde{Q})$. 

From the definition of $\widetilde{M}_{\min }$ (see (\ref{os4m})), and our
assumption that $\widetilde{\xi }$ is contained in $m_{i_0}(\widetilde{%
\Lambda }),$ there exists an arc $\widetilde{\xi }^{\prime }\subset m_{i_0}(
\widetilde{\Lambda }^{\prime }),$ whose endpoints are also contained in $%
\cup _{i=1}^Jm_i(\widetilde{\beta })),$ such that 
\begin{equation}
\label{inclusion}\widetilde{\xi }\subset interior\left( filled(\cup
_{i=1}^Jm_i(\widetilde{M}_{\min })\cup \widetilde{\xi }^{\prime })\right) .
\end{equation}

This implies that 
$$
Strip_{[\widetilde{\xi },\widetilde{\xi }^{\prime }]}=closure\left(
filled(\cup _{i=1}^Jm_i(\widetilde{M}_{\min })\cup \widetilde{\xi }^{\prime
})\backslash filled(\cup _{i=1}^Jm_i(\widetilde{M}_{\min }))\right) 
$$
has two types of boundary points:

\begin{itemize}
\item  an inner boundary, contained in $\partial \left( filled(\cup
_{i=1}^Jm_i(\widetilde{M}_{\min }))\right) $ and containing $\widetilde{\xi }%
;$

\item  an outer boundary, equal to $\widetilde{\xi }^{\prime };$
\end{itemize}

The inclusion in (\ref{inclusion}), together with the facts that $\widetilde{%
f}^{n_k}(\widetilde{\xi })$ is the part of the boundary of $\widetilde{C}_k$
which crosses $g_k.m_{i_1}(\widetilde{R})$ from outside $g_k(filled(\cup
_{i=1}^Jm_i(\widetilde{M}_{\widetilde{\Lambda }}\cup \widetilde{M}_{ 
\widetilde{\Lambda }^{\prime }}\cup \widetilde{R})))$ to inside and $
\widetilde{f}^{n_k}(Strip_{[\widetilde{\xi },\widetilde{\xi }^{\prime
}]})\cap g_k\left( \cup _{i=1}^Jm_i(\widetilde{\Lambda }\cup \widetilde{%
\Lambda }^{\prime })\right) =\emptyset $ (true by remark \ref{remarkremark}%
), imply that $\widetilde{f}^{n_k}(\widetilde{\xi }^{\prime })$ also has to
cross $g_k.m_{i_1}(\widetilde{R})$ from outside $g_k(filled(\cup
_{i=1}^Jm_i( \widetilde{M}_{\widetilde{\Lambda }}\cup \widetilde{M}_{ 
\widetilde{\Lambda }^{\prime }}\cup \widetilde{R})))$ to inside. 
And this implies the existence of a \ ''rectangle'' as in the statement of
fact \ref{factfact} contained in 
$$
\widetilde{f}^{n_k}(Strip_{[\widetilde{\xi },\widetilde{\xi }^{\prime
}]})\cap g_km_{i_1}(\widetilde{R}). 
$$
So, $\widetilde{f}^{n_k+N}(m_{i_0}(\widetilde{R}))\cap g_k.m_{i_1}( 
\widetilde{R})$ contains such a \ ''rectangle'' and thus

$$
\widetilde{f}^{n_k+N}(\widetilde{R})\cap m_{i_0}^{-1}g_km_{i_1}(\widetilde{R}%
)\supset \widetilde{R}_1. 
$$
Clearly, $\widetilde{f}^{n_k+N}(\widetilde{R})\cap (\widetilde{R})\supset 
\widetilde{R}_0,$ by our choice of $\widetilde{\Lambda }$ and $\widetilde{%
\Lambda }^{\prime }.$\ 

\qed

\vskip 0.2truecm


So, finally we can build a ''topological horseshoe'': Arguing exactly as
when all crossings are $C^1$-transversal, it can be proved that for every
bi-infinite sequence in $\{0,1\}^{\mathbb{Z}}$, denoted $(a_n)_{n\in 
\mathbb{Z}}$, there is a compact set which realizes it (not necessarily a
point, as in the $C^1$-transverse case, see \cite{bw95}, and also \cite{jlm}
for a simpler application of the above construction).

If we denote by $M_k\subset R,$ the compact set associated with the sequence 
$(1)_{\mathbb{Z}}$ and $\widetilde{M}_k=\pi ^{-1}(M_k)\cap \widetilde{R},$
then by our construction $\widetilde{f}^{m(N+n_k)}(\widetilde{M}%
_k)=(m{}_{i_0}^{-1}g_km_{i_1})^m(\widetilde{M}_k),$ for all $m>0.$ In
particular, if $r\in M_k$ and $\widetilde{r}\in \pi ^{-1}(r)\cap \widetilde{M%
}_k,$ then $\widetilde{f}^{m(N+n_k)}(\widetilde{r})\in
(m{}_{i_0}^{-1}g_km_{i_1})^m(\widetilde{M}_k),$ for all $m>0.$

By our choice of $R,$ for all $m>0$ we can find $\beta _m^{\prime }$ a path
in $R$ joining $f^{m(N+n_k)}(r)$ to $r$ with $l(\beta _m^{\prime
})<2\epsilon .$ Thus, if $\widetilde I_{\widetilde r}^{m(N+n_k)}*\widetilde
\beta _m^{\prime }$ is the lift of $I_r^{m(N+n_k)}*\beta _m^{\prime }$ with
base point $\widetilde r$, then $\widetilde I_{\widetilde
r}^{m(N+n_k)}*\widetilde \beta _m^{\prime }$ is a path in $\mathbb{D}$
joining $\widetilde r$ to $(m_{i_0}^{-1}g_km_{i_1})^m(\widetilde r)$. In
particular 
$$
\frac{[I_p^{m(N+n_k)}*\beta _m^{\prime }]}{m(N+n_k)}=\frac{%
[(m_{i_0}^{-1}g_km_{i_1})^m]}{m(N+n_k)}=\frac{m[m_{i_0}^{-1}g_km_{i_1}]}{%
m(N+n_k)}=\frac{[m_{i_0}^{-1}]+[g_k]+[m_{i_1}]}{N+n_k}. 
$$

As $w=\lim _{k\to \infty }[g_k]/n_k,$ $N>0$ is fixed and there is just a
finite number of possibilities for $m_{i_0}$ and $m_{i_1},$ if $k>0$ is
large enough, then

$$
\frac{[m_{i_0}^{-1}]+[g_k]+[m_{i_1}]}{N+n_k}\text{ is as close as we want to 
}w. 
$$
So given an $error>0,$ if $k>0$ is sufficiently large, defining 
\begin{equation}
\label{defgwnw}g_w=m_{i_0}^{-1}g_km_{i_1}\text{ and }n_w=N+n_{k,} 
\end{equation}
we get that 
$$
\left\| \frac{[g_w]}{n_w}-w\right\| <error. 
$$

Using the above construction, we will show that $\rho _{mz}(f)=Conv(\rho
_{mz}(f))$. For this we need Steinitz's theorem \cite{hdk}. This theorem
says that if a point is interior to the convex hull of a set $X$ in $%
\mathbb{R}^n,$ it is interior to the convex hull of some set of $2n$ or
fewer points of $X.$

Since $\rho _{mz}(f)$ is a compact set, $Conv(\rho _{mz}(f))=Conv(Ext(\rho
_{mz}(f)),$ where $Ext(\rho _{mz}(f))$ is the set of all extremal points of $%
Conv(\rho _{mz}(f))$. Using Steinitz's theorem, any point in the interior of 
$Conv(\rho _{mz}(f))$ is a convex combination of at most $4g$ extremal
points.

Let $v$ be a point in $int(Conv(\rho _{mz}(f)))\cap \mathbb{Q}^{2g}$. By the
previous observation, there exists at most $4g$ extremal points (here
without loss of generality we will assume that exactly $4g$ extremal points
are used) $w_1,\ldots ,w_{4g}$ such that 
$$
v=\sum_{i=1}^{4g}\lambda _iw_i, 
$$
where $\lambda _i\in ]0,1[,$ for all $1\leq i\leq 4g$ and $\lambda _1+\ldots
+\lambda _{4g}=1$. By the previous construction for some general $w,$ choose
deck transformations $g_{w_1},\ldots ,g_{w_{4g}}$ and natural numbers $%
n_{w_1},\ldots ,n_{w_{4g}}$ (as in expression (\ref{defgwnw})), such that 
$$
v\in int\Big(Conv\Big(\frac{[g_{w_1}]}{n_{w_1}},\ldots ,\frac{[g_{w_{4g}}]}{%
n_{w_{4g}}}\Big)\Big). 
$$

This is always possible since $[g_{w_i}]/n_{w_i}$ can be chosen as close as
desired to $w_i$. As $[g_{w_i}]/n_{w_i}\in \mathbb{Q}^{2g}$ for all $1\leq
i\leq 4g$, and $v$ is a rational point in the interior of the convex hull of
these points, there exists $\lambda _1^{\prime },\ldots ,\lambda
_{4g}^{\prime }$, with $\lambda _i^{\prime }\in (0,1)\cap \mathbb{Q}$, $%
\lambda _1^{\prime }+\ldots +\lambda _{4g}^{\prime }=1$ such 
$$
v=\sum_{i=1}^{4g}\lambda _i^{\prime }\frac{[g_{w_i}]}{n_{w_i}}. 
$$

Thus, multiplying both sides of the previous equation by an appropriate
positive integer, we get positive integers $a_{Total},a_1,\ldots ,a_{4g}$
such that $a_{Total}=a_1+\ldots +a_{4g}$ and 
$$
a_{Total}v=\sum_{i=1}^{4g}a_i\frac{[g_{w_i}]}{n_{w_i}}, 
$$

For each $i\in \{1,2,...,4g\},$ $\widetilde f^{n_{w_i}}(\widetilde R)$
intersects $g_{w_i}(\widetilde R)$ in a vertical rectangle as in fact \ref
{factfact}. Since $\widetilde f$ commutes with every deck transformation, $%
\widetilde f^{n_{w_j}}(g_{w_i}(\widetilde R))$ intersects $%
g_{w_i}g_{w_j}(\widetilde R)$ in a similar rectangle. 

Let $N_{product}=n_{w_1}n_{w_2}\ldots n_{w_{4g}}$ and for all $1\leq i\leq
4g,$ let $u_i=N_{product}/n_{w_i}.$ By the previous definitions, $\widetilde{%
f}^{(a_iu_i)n_{w_i}}(\widetilde{R})=\widetilde{f}^{a_iN_{product}}(
\widetilde{R})$ satisfies: 
$$
\widetilde{f}^{(a_iu_i)n_{w_i}}(\widetilde{R})\cap g_{w_i}^{a_iu_i}(
\widetilde{R})\text{ contains a vertical rectangle as in fact \ref{factfact}.%
} 
$$
So, considering all iterates of this type for $1\leq i\leq 4g$ and composing
them we obtain that 
$$
\widetilde{f}^{a_{Total}N_{product}}(\widetilde{R})\cap h_v(\widetilde{R})
\text{ contains a vertical rectangle as in fact \ref{factfact}}, 
$$
where 
$$
h_v=g_{w_1}^{a_1u_1}\circ g_{w_2}^{a_2u_2}\circ \ldots \circ
g_{w_{4g}}^{a_{4g}u_{4g}}. 
$$

Clearly, as%
$$
\widetilde f^{a_{Total}N_{product}}(\widetilde R)\cap \widetilde R\text{
contains a vertical rectangle similar to }\widetilde{R}_0,\text{ just
thinner,} 
$$
we can consider the compact $f^{a_{Total}N_{product}}$-invariant subset $%
K_v\subset R$ of the topological horseshoe we just produced, associated with
the sequence $(1)_{\mathbb{Z}}.$ If $\widetilde{K}_v=\widetilde{R}\cap \pi
^{-1}(K_v)$, then 
$$
\widetilde f^{a_{Total}N_{product}}(\widetilde{K}_v)=h_v(\widetilde{K}_v). 
$$

\begin{figure}[!h]
	\centering
	\includegraphics[scale=0.2]{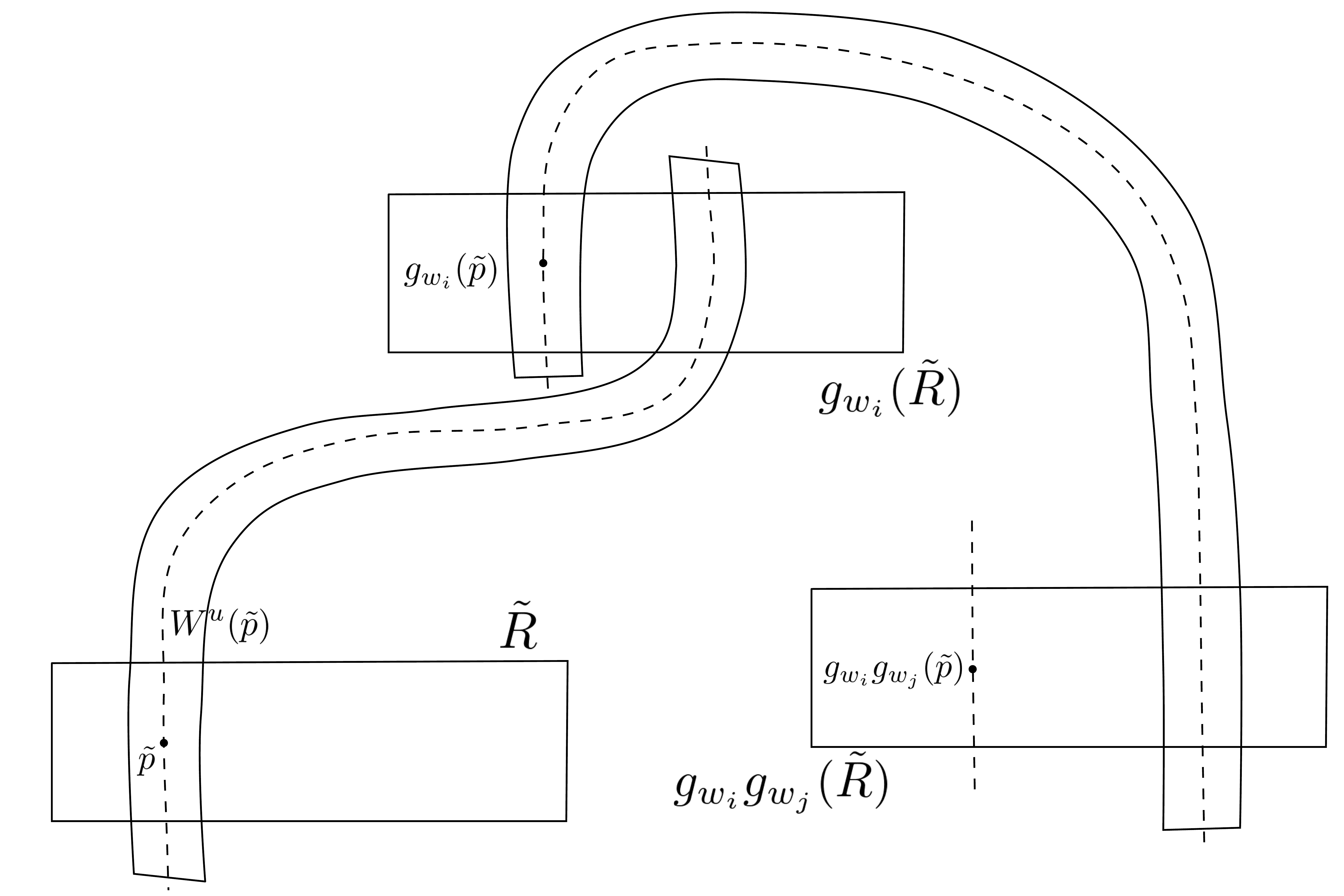}
	\caption{How to create intersections between iterates of 
		$\widetilde{R}$ and its translates.}
	\label{intersection}
\end{figure}

So, $h_v^{-1}\widetilde{f}^{a_{Total}N_{product}}(\widetilde{K}_v)=
\widetilde{K}_v,$ which implies using Brouwer's lemma on translation arcs 
\cite{franks}, that $h_v^{-1}\widetilde{f}^{a_{Total}N_{product}}$ has a
fixed point $\widetilde{z}_v.$ Since 
$$
\widetilde{f}^{a_{Total}N_{product}}(\widetilde{z}_v)=h_v(\widetilde{z}_v) 
$$
and 
$$
\frac{[h_v]}{a_{Total}N_{product}}=\frac{[g_{w_1}^{a_1u_1}\circ
g_{w_2}^{a_2u_2}\circ \ldots \circ g_{w_{4g}}^{a_{4g}u_{4g}}]}{%
a_{Total}N_{product}}=\sum_{i=1}^{4g}\frac{a_iu_i[g_{w_i}]}{%
a_{Total}N_{product}}=\frac 1{a_{Total}}\sum_{i=1}^{4g}a_i\frac{[g_{w_i}]}{%
n_{w_i}}=v, 
$$
we conclude that the $f$-periodic point $z_v=\pi (\widetilde{z}_v)$ has a
rotation vector $\rho (f,z_v)=v.$ This shows that $v\in \rho _{mz}(f).$
Since $\rho _{mz}(f)$ is compact, $\rho _{mz}(f)=Conv(\rho _{mz}(f)).$

Now let us consider the irrational case. For any $v\in (\mathbb{Q}%
^{2g})^c\cap int(\rho _{mz}(f)),$ exactly as in the rational case one can
find $4g$ rational points $w_1,\ldots ,w_{4g}$ in $\rho _{mz}(f)$ for which 
$$
v\in int(Conv(\{w_1,\ldots ,w_{4g}\})) 
$$
and such that for positive integers $n_{w_1},\ldots ,n_{w_{4g}},$ the
following holds:

$$
\begin{array}{c}
\widetilde f^{n_{w_i}}(\widetilde R)\cap g_{w_i}(\widetilde R) 
\text{ contains a vertical rectangle as in fact \ref{factfact},} \\ \text{%
for some }g_{w_i}\in Deck(\pi )\text{ such that }\left[ g_{w_i}\right]
/n_{w_i}=w_i. 
\end{array}
$$

As we did above, let $N_{product}=n_{w_1}\ldots n_{w_{4g}}$ and $%
u_i=N_{product}/n_{w_i}.$ Then, 
$$
\widetilde f^{N_{product}}(\widetilde R)\cap g_{w_i}^{u_i}(\widetilde R) 
\text{ also contains a vertical rectangle }\widetilde{R}_i\text{ as in fact 
\ref{factfact}.} 
$$

Clearly 
$$
\frac{\left[ g_{w_i}^{u_i}\right] }{N_{product}}=\frac{u_i.\left[
g_{w_i}\right] }{N_{product}}=w_i. 
$$

So going back to the surface $S,$ 
$$
f^{N_{product}}(R)\cap R\supseteq R_1\cup ...\cup R_{4g},\text{ where }%
R_i=\pi (\widetilde{R}_i). 
$$
And we claim that there exists an infinite sequence in $\{1,...,4g\}^{%
\mathbb{N}}$, denoted 
$$
a_1a_2...a_n... 
$$
such that for some constant $C^{*}>0,$ 
$$
\left\| \sum_{i=1}^n\left[ g_{w_{a_i}}^{u_{a_i}}\right]
-nN_{product}.v\right\| <C^{*},\text{ for all }n>0. 
$$
The existence of this kind of sequence is what is behind, in the torus case,
of the realization of irrational rotation vectors in the interior of the
rotation set by compact invariant sets with bounded mean motion in the
universal cover. This was done for relative pseudo-Anosov maps in lemma 3 of 
\cite{MZ} and extended to the original map using Handel's shadowing. See 
\cite{MZ} for details.

Now let $z\in R$ be any point which corresponds to the sequence $%
a_1a_2...a_n...,$ namely $f^{nN_{product}}(z)\in R_{a_n},$ for all $n\geq 1.$
Clearly for $\widetilde{z}\in \widetilde{R}\cap \pi ^{-1}(z)$ and any $n\geq
1,$%
$$
\widetilde{f}^{n.N_{product}}(\widetilde{z})\in
g_{w_{a_1}}^{u_{a_1}}.g_{w_{a_2}}^{u_{a_2}}...g_{w_{a_n}}^{u_{a_n}}(
\widetilde{R}), 
$$
so not only the rotation vector of $z$ is $v,$ but%
$$
\left\| [\alpha _z^l]-l.v\right\| <C^{*}+N_{product}.\left\| v \right\| +2\epsilon+\max \{d_{\mathbb{D}}(
\widetilde{f}^i(\widetilde{z}),\widetilde{z}):\widetilde{z}\in \mathbb{D}
\text{ and }0\leq i\leq N_{product}\}. 
$$
And this implies that the $\omega $-limit set of $z,$ denoted $K_v,$ has the
property we are looking for, because for any $z^{\prime }\in K_v,$ $\left\|
[\alpha _{z^{\prime }}^n]-n.v\right\| $ is smaller than some constant which is independent of $n$ and $z^{\prime}\in K_v.$

\qed

\vskip 0.2truecm

\section{Proof of theorem 4}

Here we just make use of the machinery developed in the proof of theorem 3.

Suppose, by contradiction, that for every $M>0,$ there exists $\omega \in
\partial \rho _{mz}(f),$ a supporting hyperplane $\omega \in H\subset 
\mathbb{R}^{2g},$ $z\in S$ and $n>0$ such that 
$$
\left( [\alpha _z^n]-n.\omega \right) .\overrightarrow{v_H}>M, 
$$
where $\overrightarrow{v_H}$ is the unitary normal vector to $H$ pointing
towards the connected component of $H^c$ which does not intersect $\rho
_{mz}(f).$

Fixed some fundamental domain of $S,$ denoted $\widetilde{Q}\subset 
\mathbb{D},$ there exists $\widetilde{z}=\pi ^{-1}(z)\cap \widetilde{Q}$
such that for some $g\in Deck(\pi ),$ 
$$
\widetilde f^n(\widetilde{z})\in g(\widetilde{Q})\text{ and }\left( \left[
g\right] -n.\omega \right) .\overrightarrow{v_H}>M-C_{\widetilde{Q}}, 
$$
where $C_{\widetilde{Q}}>0$ is a constant which depends only on the shape of 
$\widetilde{Q}.$ From the proof of the previous theorem, we know that there
are deck transformations $\left\{ m_1,m_2,...,m_J\right\} ,$ for some $J>0,$
which do not depend on the choices of

\begin{itemize}
\item  $M>0,$ $\omega \in \partial \rho _{mz}(f),$ the supporting hyperplane 
$\omega \in H\subset \mathbb{R}^{2g},$ $z\in S$ and $n>0,$
\end{itemize}

\noindent such that for some $i_0$ and $i_1$ in $\{1,...,J\},$ there exists
a compact subset $\widetilde{K}_M$ for which

$$
\widetilde{f}^{N+n}(\widetilde{K}_M)=m{}_{i_0}^{-1}gm_{i_1}(\widetilde{K}%
_M), 
$$
where $N>0$ is given in expression (\ref{defineNunif}). And thus for some
point $\widetilde{z}_M\in \mathbb{D},$ $\widetilde{f}^{N+n}(\widetilde{z}%
_M)=m{}_{i_0}^{-1}gm_{i_1}(\widetilde{z}_M).$

So, if $M>0$ is large enough so that 
$$
\left( \left[ m{}_{i_0}^{-1}gm_{i_1}\right] -(n+N).\omega \right) . 
\overrightarrow{v_H}>0, 
$$
we get a contradiction.

\qed

\vskip 0.2truecm

\section{Proof of theorem 5}

This proof is very similar to the proof of theorem 2 of \cite{jlm}. In
particular, the following lemma from that paper, which was proved in the
torus, holds without any modification under the hypotheses of the present
paper:

\begin{lemma}
(Adapted lemma 6 of \cite{jlm}) Suppose $f:S\rightarrow S$ is a $%
C^{1+\epsilon }$ diffeomorphism isotopic to the identity which has a fully
essential system of curves $\mathscr{C}.$ Let $\mu $ be a $f$-invariant
Borel probability measure such that its rotation vector $\rho (\mu )$
belongs to $\partial \rho _{mz}(f).$ Let $H$ be a supporting hyperplane at $%
\rho (\mu )$ and $\overrightarrow{v_H}$ be the unitary vector orthogonal to $%
H,$ pointing towards the connected component of $H^c$ which does not
intersect $\rho _{mz}(f).$ Then, if $x^{\prime }\in supp(\mu ),$ for any
integer $n>0,$

\begin{equation}
\label{fff}\left| \left( [\alpha _{x^{\prime }}^n]-n.\rho (\mu )\right) .
\overrightarrow{v_H}\right| \leq 2+M(f),
\end{equation}
where $M(f)>0$ comes from theorem 4.
\end{lemma}

Now the proof continues exactly as the proof of theorem 2 of \cite{jlm}. 
\qed

\vskip 0.2truecm

{\it Acknowledgements: }The first author is grateful to Juliana Xavier for
some conversations on this subject several years ago. After doing the torus
case I asked her on ways to generalize the results to other surfaces and she
explained several thing to me. Both of us are grateful to Andre de Carvalho
for several conversations and as we already said, to Alejandro Passeggi for
pointing out that theorem 3 was a consequence of theorem 2. 

\end{document}